# Polygonal surfaces in pseudo-hyperbolic spaces

Alex Moriani [1]

*Laboratoire J.A. Dieudonné, Université Côte d'Azur, France*

ABSTRACT

A polygonal surface in the pseudo-hyperbolic space $\mathbf{H}^{2,n}$ is a complete maximal surface bounded by a lightlike polygon in the Einstein universe $\mathbf{Ein}^{1,n}$ with finitely many vertices. In this article, we give several characterizations of them. Polygonal surfaces are characterized by finiteness of their total curvature and by asymptotic flatness. They have parabolic type and polynomial quartic differential. Our result relies on a comparison between three ideal boundaries associated with a maximal surface, corresponding to three distinct distances naturally defined on the maximal surface.

**Contents**

*E-mail address:* amoriani@proton.me.

[1] The author acknowledges funding by the European Research Council under ERC-Advanced grant 101095722.

https://doi.org/10.1016/j.aim.2025.110484

## 1. Introduction

### *1.1. General context*

The study of minimal surfaces is a broad and important part of differential geometry. A *minimal surface* is a critical point for the area functional on a Riemannian manifold. Some classical questions to ask are, given a boundary condition:

- Is there a minimal surface satisfying this boundary condition (the *Plateau problem*)?
- If a minimal surface exists, is it unique?
- What is the regularity one can expect for the minimal surface?
- What are the geometric properties of the surface that are imposed by the boundary condition?

When the Riemannian manifold $M$ considered has nonpositive curvature, it has an ideal boundary. We fix a simple closed curve, or a union of simple closed curves, in this boundary and ask for the existence and uniqueness of a minimal surface asymptotic to this curve. This is the *asymptotic Plateau problem.*

In the Riemannian case, already uniqueness is an issue. Surprisingly, the situation is much simpler in the pseudo-hyperbolic space $\mathbf{H}^{2,n}$.

Let us fix $(\mathbf{E}, \mathbf{q})$ a real vector space endowed with a quadratic form of non-degenerate signature $(2, n+1)$. We define $\mathbf{H}^{2,n}$ as being the set of $\mathbf{q}$-negative definite lines of $\mathbf{E}$. It is a complete pseudo-Riemannian manifold of signature $(2, n)$ and constant curvature $-1$. The Einstein universe $\mathbf{Ein}^{1,n}$, defined as being the frontier of $\mathbf{H}^{2,n}$ in $\mathbf{P}(\mathbf{E})$, is a natural boundary for the pseudo-hyperbolic space. Given a complete spacelike surface $\Sigma$ in $\mathbf{H}^{2,n}$, we define its *space boundary* $\partial^{\mathrm{s}}\Sigma$ as being the frontier of $\Sigma$ in $\mathbf{H}^{2,n} \cup \mathbf{Ein}^{1,n}$.

In [20], Labourie Toulisse and Wolf study the asymptotic Plateau problem in the pseudo-hyperbolic space $\mathbf{H}^{2,n}$ for a certain class of curves. They call *semi-positive loop* in $\mathbf{Ein}^{1,n}$ a loop that is the space boundary of a complete spacelike surface of $\mathbf{H}^{2,n}$. Their result is the following: a semi-positive loop in $\mathbf{Ein}^{1,n}$ bounds a unique complete *maximal* surface in $\mathbf{H}^{2,n}$. Maximal surfaces are supposed to be complete in our article, but the completeness of entire maximal surfaces seems to be an empty condition (see [31] for anti de-Sitter space $\mathbf{H}^{2,1}$).

### 1.2. Main theorem

If $\Sigma$ is a maximal surface in $\mathbf{H}^{2,n}$, then its sectional curvature takes values between $-1$ and $0$ (by [4]). Two extreme examples of maximal surfaces in $\mathbf{H}^{2,n}$ are the following.

- A totally geodesic copy of the hyperbolic plane, bounded by a projective circle,
- A Barbot surface, isometric to the Euclidean plane, bounded by a lightlike polygon with 4 vertices.

We study the problem of characterizing the maximal surfaces that have a behavior close to the behavior of Barbot surfaces.

Our main theorem is the following:

**Theorem 1.2.1** *(Main Theorem). Let $\Sigma$ be a maximal surface in $\mathbf{H}^{2,n}$, and denote by $g$ its induced metric. The following are equivalent:*

*(i) The space boundary of $\Sigma$ is a lightlike polygon with finitely many vertices,*
*(ii) The total curvature of $(\Sigma, g)$ is finite,*
*(iii) The surface $(\Sigma, g)$ is asymptotically flat.*

*Moreover, if the conditions are satisfied, and denoting by $N+4$ the number of vertices of $\partial^{\mathrm{s}}\Sigma$, the Riemann surface $(\Sigma, [g])$ is of parabolic type, has polynomial quartic differential of degree $N$, the total curvature of $\Sigma$ equals $-\frac{\pi}{2}N$ and the circles of radius $r$ have perimeter asymptotic to $(N+4)\frac{\pi}{2}r$ when $r \to +\infty$.*

The quartic differential is a holomorphic object associated with the second fundamental form of a maximal surface and will be defined below in Section 1.4.3.

Every maximal surface in $\mathbf{H}^{2,n}$ is a Hadamard surface. We will use techniques coming from geometry of non-positively curved manifolds to understand the geometry of maximal surfaces and its links with the space boundary.

There are three distance functions defined on a maximal surface $\Sigma$. The *induced distance function*, the *spacelike distance function*, and the *quartic distance function.*

- The induced distance function is the distance associated with the induced metric.
- The spacelike distance function between two points is the length of the $\mathbf{H}^{2,n}$-geodesic segment joining these two points.
- The quartic distance function is the distance associated with the flat metric with conical singularities coming from the quartic differential of $\Sigma$.

Remark that the first item of the Main Theorem deals with the spacelike distance function, that the second and third items deal with the induced distance function and that the fact that the quartic differential is polynomial deals with the quartic distance function.

The main idea of the proof is to assign an ideal boundary to each of the three distance functions in a natural manner, and to understand how they are related.

### *1.3. Related works*

We list here some works related with ours. They deal with minimal surfaces in homogeneous spaces, and with Higgs bundles in appropriate cases.

#### *1.3.1. Minimal surfaces in Euclidean spaces*

By a theorem of Osserman and Chern, [5, Theorem 1], a complete minimal surface $(\Sigma, g)$ in $\mathbb{R}^n$ has finite total curvature if and only if the following holds:

(i) the Riemann surface $X = (\Sigma, [g])$ is a compact Riemann surface with a finite number of punctures,
(ii) the minimal embedding corresponds to a collection $\{\varphi_1, \dots, \varphi_n\}$ of abelian differentials on $X$ that extend meromorphically at the punctures.

*Total curvature* The total curvature equals $-2\pi N$ where $N$ is a natural integer.

All finite total curvature complete minimal surfaces in the euclidean $n$-space are not classified. Some partial classification exists when we ask some topological restrictions (see [25, Chapter I, Section 3]). For example, if $\Sigma$ has genus 0 in $\mathbb{R}^3$, then it is the flat plane (total curvature 0) or the catenoid (total curvature $-2\pi$).

*Polygonal boundary* Let $S$ be a finite total curvature minimal surface in $\mathbb{R}^3$, corresponding to the compact Riemann surface $X$ with punctures $\{p_1, \dots, p_k\}$. Each end of $S$ is asymptotic to a plane covered a finite number of time as we explain now. Let $E_j$ be an end of $S$, corresponding to a sufficiently small neighborhood of $p_j$. Then ([15, Theorem 2.3]), denoting by $\mathbf{S}^2$ the unit sphere in $\mathbb{R}^3$, the set

$$S_j(r) = \{u \in \mathbf{S}^2 \mid ru \in E_j\}$$

converges, when $r$ goes to $+\infty$, to a great circle covered an integral number of times $d_j$ and the total curvature of $S$ can be computed to be:

$$\int_S K_S \mathrm{d}v_g = 2\pi \left( \chi(S) - \sum_{j=1}^{k} d_j \right) .$$

*1.3.2. Minimal surfaces in $\mathbf{H}^2 \times \mathbb{R}$*

In [17], Hauswirth and Rosenberg study complete minimal surfaces in $\mathbf{H}^2 \times \mathbb{R}$ with finite total curvature.

*Total curvature* By analyzing the harmonic equations, Hauswirth and Rosenberg show that the total curvature of $S$ is a mutiple of $-2\pi$, related to the orders of the zeros of $Q$.

*Polygonal boundary* The first examples of minimal surfaces in $\mathbf{H}^2 \times \mathbb{R}$ were *Scherk graphs* over ideal polygonal domains (see [6]).

In general [16, Theorem 2.1] the boundary of a minimal surface in $\mathbf{H}^2 \times \mathbb{R}$ is a polygonal curve consisting of geodesics in $\mathbf{H}^2 \times \{\pm\infty\}$ and vertical straight lines joining their endpoints.

In [26], Pyo and Rodriguez show that the Scherk graphs are not the only minimal surfaces with finite total curvature.

*1.3.3. Constant mean curvature cuts in Minkowski space*

Choi and Treibergs, in [30], [8], [7], worked on the link between CMC cuts in Minkowski space $\mathbf{R}^{n,1}$ and harmonic maps from $\mathbb{C}$ to the hyperbolic plane. Han, Tam, Treibergs and Wan [18] then used these results to establish a link between the growth of the Hopf differential of harmonic diffeomorphisms from $\mathbb{C}$ to $\mathbf{H}^2$ and the image of the diffeomorphism.

The Gauss map of a CMC cut $\Sigma$ can be thought of as being a map from $\Sigma$ to $\mathbf{H}^2$. Indeed the set of unit timelike future pointing vectors at a point of $\mathbb{R}^{2,1}$ is the hyperboloid model of $\mathbf{H}^2$. In [8, Theorem 4.8], it is shown that the image of the Gauss map is the interior of a closed convex set. We call the boundary of this convex set the *boundary* of the CMC cut.

We can now state their results as follows.

**Theorem** *(Choi–Treibergs, Han–Tam–Treibergs–Wan). Let $\bar{\Sigma}$ be a CMC cut of $\mathbb{R}^{2,1}$. The following are equivalent:*

*(i) The boundary of $\bar{\Sigma}$ is an ideal polygon with finitely many vertices,*
*(ii) The total curvature of $(\bar{\Sigma}, g)$ is finite,*

*(iii) The Riemann surface* $(\bar{\Sigma}, [g])$ *is of parabolic type and the Hopf differential of the Gauss map is polynomial.*

There is a correspondence between CMC cuts in $\mathbb{R}^{2,1}$ and maximal surfaces in $\mathbf{H}^{2,1}$. If $\Sigma$ is a maximal surface in $\mathbf{H}^{2,1}$ with induced metric $g$ and shape operator $B$, then the tensors $g$ and $B + \lambda \mathrm{id}$ satisfy the Gauss and Codazzi equations of a CMC cut of mean curvature $\lambda$ in $\mathbb{R}^{2,1}$.

The items of our main theorem are easily translated in terms of properties of associated CMC cuts. As a corollary of our main theorem, we obtain a new characterization of polygonal CMC cuts, being asymptotically flat:

**Corollary 1.3.1.** *Let* $\bar{\Sigma}$ *be a CMC cut of* $\mathbb{R}^{2,1}$*. The following are equivalent:*

*(i) The boundary of* $\bar{\Sigma}$ *is an ideal polygon with finitely many vertices,*
*(ii) The total curvature of* $(\bar{\Sigma}, g)$ *is finite,*
*(iii) The surface* $(\bar{\Sigma}, g)$ *is asymptotically flat,*
*(iv) The Riemann surface* $(\bar{\Sigma}, [g])$ *is of parabolic type and the Hopf differential of the Gauss map is polynomial.*

**Remark 1.3.2.** Our main theorem holds in any codimension, unlike the case of CMC cuts that are surfaces in a 3-dimensional space.

In the case of real split groups of rank 2, the situation has been extensively studied in the following works.

*1.3.4.* $\mathrm{SL}_3(\mathbb{R})$ *and affine spheres in* $\mathbb{R}^3$

There is a correspondence, called Calabi-Cheng-Yau correspondence, between properly convex open sets in $\mathbb{RP}^2$ and hyperbolic affine spheres of curvature $-1$ in $\mathbb{R}^3$. A hyperbolic affine sphere comes naturally with a complex structure and a holomorphic cubic differential. Dumas and Wolf, in [10], show that a hyperbolic affine sphere is of parabolic type with polynomial cubic differential if and only if the underlying set in the projective space is a polygon. Their proof is based on the resolution of the vortex equation and the study of the asymptotics of a solution. In their work, the easiest case is that of Ţiţeica surfaces, that play the role of Barbot surfaces in our setting.

*1.3.5.* $\mathrm{SL}_2(\mathbb{R}) \times \mathrm{SL}_2(\mathbb{R})$ *and lightlike polygons in* $\mathbf{Ein}^{1,1}$

In [29], Tamburelli studied polygonal surfaces in $\mathbf{H}^{2,1}$. He constructed a homeomorphism between the moduli space of polynomial quadratic differentials on the complex plane and the moduli space of lightlike polygons in $\mathbf{Ein}^{1,1}$. He interpreted the maximal surface equation as an instance of the vortex equation.

This is the counterpart of the work of Dumas and Wolf is the setting of maximal surfaces in $\mathbf{H}^{2,1}$.

*1.3.6.* $\mathrm{Sp}(4,\mathbb{R})$ *and lightlike polygons in* $\mathbf{Ein}^{1,2}$

Tamburelli and Wolf [33] also studied polygonal maximal surfaces in $\mathbf{H}^{2,2}$. They note that the space of polygons in $\mathbf{Ein}^{1,2}$ is not connected and describe a component of the moduli space of polygonal surfaces in $\mathbf{H}^{2,2}$ in terms of their holomorphic data. They do not deal with the other components.

*1.3.7.* $G_2$ *and almost complex curves in* $\mathbf{S}^{2,4}$

In another context, there is the work of Evans on polygonal surfaces in $\mathbf{S}^{2,4}$. In [11] he constructed a map from the moduli space of holomorphic polynomial sextic differentials of degree $k$ on the complex plane to the moduli space of annihilator polygons with $k+6$ vertices in $\mathbf{Ein}^{2,3}$.

This is close to our setting. The space $\mathbf{S}^{2,4}$ is anti-isometric to $\mathbf{H}^{4,2}$ and the surfaces that he is considering are maximal spacelike surfaces in $\mathbf{H}^{4,2}$.

*1.3.8. Quasiperiodic surfaces in* $\mathbf{H}^{2,n}$

In [19], Labourie and Toulisse investigate a subset of the set of semi-positive loops in $\mathbf{Ein}^{1,n}$: the *quasiperiodic* loops. They found interesting characterizations concerning the maximal surfaces bounded by a quasiperiodic loop, among which the fact of having curvature bounded by above by a negative constant.

Our surfaces are, in some sense, the opposite of quasiperiodic surfaces: polygonal surfaces are asymptotically flat.

## *1.4. Overview of the proof*

We explain now how goes the proof of the Main Theorem. The proof lies in a comparison between ideal boundaries associated with the three distances.

There are two main constructions for the ideal boundary of a Hadamard manifold: via rays and via horofunctions. We use both constructions to assign an ideal boundary to $\Sigma$ when endowed with these different distances.

*1.4.1. The ideal boundary*

A maximal surface $\Sigma$ in $\mathbf{H}^{2,n}$ together with its induced metric $g$ is a Hadamard surface. The ideal boundary $\Sigma(\infty)$ of $\Sigma$ is defined as being the set of equivalence classes of unit speed geodesic rays, up to bounded distance (Section 2.3.1). The ideal boundary strongly determines the geometry of the surface via the Tits distance. This distance is defined on the ideal boundary and the Tits distance has finite diameter on $\Sigma(\infty)$ if and only if $\Sigma$ has finite total curvature (see Section 4.8.1).

*1.4.2. The space boundary*

We have defined the space boundary $\partial^{\mathrm{s}}\Sigma$ of $\Sigma$ as being the set of limits of sequences of points in $\Sigma$ that converge in $\mathbf{Ein}^{1,n}$. As for the ideal boundary (Section 2.3.2), the space boundary can be realized as the frontier of an embedding of $\Sigma$ in the space $\mathcal{C}^*(\Sigma) = \mathcal{C}(\Sigma)/\mathbb{R}$ of continuous functions on $\Sigma$ modulo additive real constants.

The frontier of this embedding consists of a set of functions on $\Sigma$, defined up to an additive constant, called space-horofunctions, that are in bijection with the space boundary $\partial^{\mathrm{s}}\Sigma$ and have a totally explicit description (Section 2.3.3).

### *1.4.3. The quartic boundary*

We can define, using the second fundamental form on $\Sigma$, a holomorphic quartic differential on the underlying Riemann surface $(\Sigma, [g])$ as follows. We denote by $g_{\mathsf{N}}$ the induced metric on the normal bundle of $\Sigma$.

The second fundamental form $\mathbb{I}$ of a maximal surface $\Sigma$ is a section of the bundle $(\mathsf{T}^*\Sigma \odot \mathsf{T}^*\Sigma) \otimes \mathsf{N}\Sigma$, where $\odot$ denotes the symmetric tensor product. Applying the complex quadratic form $g_{\mathsf{N}}^{\mathbb{C}}$ to the $(2,0)$-part of the complexified tensor $\mathbb{I}^{\mathbb{C}}$ we obtain

$$\varphi_4 := g_{\mathsf{N}}^{\mathbb{C}}(2\mathbb{I}^{2,0}, 2\mathbb{I}^{2,0}) \ ,$$

a section of $(\mathsf{T}^{1,0}\Sigma)^{\otimes 4} \otimes \mathbb{C}$. By the Codazzi equation, this is a holomorphic quartic differential on $(\Sigma, [g])$ where $[g]$ denotes the conformal class of the induced metric $g$ on $\Sigma$ (see Section 5.1).

The holomorphic quartic differential defines a flat metric with conical singularities on $\Sigma$. Since the differential $\varphi_4$ is a quartic differential, a zero of order $k$ becomes a conical singularity of angle $2\pi + 4\pi/2$ of the flat metric. This metric makes $\Sigma$ a Hadamard space (in the sense of [1]) and we can define the ideal boundary via rays in the same way as for a Hadamard manifold (see Section 5.4).

We explain now how we compare these three ideal boundaries.

### *1.4.4. Barbot surfaces and renormalization*

The simplest case of polygonal surfaces is the case of maximal surfaces bounded by a 4-gon, called *Barbot surfaces* and *Barbot crowns*. A maximal surface $\Sigma$ has nonpositive curvature, and it is flat at a point if and only if it is flat everywhere if and only if it is a Barbot surface by [19, Proposition 5.3]. This is a starting point to describe asymptotically flat maximal surfaces.

An important tool is the *renormalization* process. Let $\mathsf{G}$ be the isometry group of $\mathbf{H}^{2,n}$. The group $\mathsf{G}$ is isomorphic to $\mathbf{PO}(2, n+1)$. The space

$$\{(\Sigma, x) \mid \Sigma \text{ maximal surface},\ x \in \Sigma\}$$

of pointed maximal surfaces has a natural action by $\mathsf{G}$ that is proper, continuous and cocompact ([20, Theorem 6.1] and [19, Theorem 4.8]). Hence, given any sequence $((\Sigma_k, x_k))_{k\in\mathbb{N}}$ of pointed maximal surfaces, there is a sequence $(g_k)_{k\in\mathbb{N}}$ of elements of $\mathsf{G}$ such that $(g_k \cdot (\Sigma_k, x_k))_{k\in\mathbb{N}}$ subconverges to a pointed maximal surface. We say that the sequence $((\Sigma_k, x_k))_{k\in\mathbb{N}}$ subconverges *up to renormalization.*

For example, take an asymptotically flat maximal surface $\Sigma$, and a diverging sequence $(x_k)_{k\in\mathbb{N}}$ of points of $\Sigma$. Then the sequence $((\Sigma, x_k))_{k\in\mathbb{N}}$ subconverges, up to renormalization, to a Barbot surface.

We prove in Section 4.3.4 the following.

**Theorem 1.4.1.** *Let* $\Sigma$ *be a polygonal surface in* $\mathbf{H}^{2,n}$*. Then for every diverging sequence* $(x_k)_{k\in\mathbb{N}}$ *of* $\Sigma$*, the sequence of pointed maximal surfaces* $((\Sigma, x_k))_{k\in\mathbb{N}}$ *subconverges to a pointed Barbot surface, up to renormalization.*

Hence polygonal surfaces are asymptotically flat, and this is $(i) \Rightarrow (iii)$ of the Main Theorem.

Barbot surfaces have an explicit description in terms of their space boundary. Also the quartic metric equals the induced metric of a Barbot surface. Hence we have a total correspondence between the space boundary, the ideal boundary and the quartic boundary of a Barbot surface.

*1.4.5. The space boundary and the ideal boundary*

We compare here the space boundary $\partial^{\mathrm{s}}\Sigma$ and the ideal boundary $\Sigma(\infty)$ of a maximal surface $\Sigma$. Recall that we identified $\partial^{\mathrm{s}}\Sigma$ with the set of space-horofunctions defined on $\Sigma$ up to an additive constant.

The space-horofunctions share some properties with intrinsic horofunctions. The most important for our purpose is that space-horofunctions are strictly quasiconvex on $\Sigma$ (Theorem 4.4.4). A *space-horoball* is a sublevel set of a space-horofunction. Hence the space-horoballs are strictly convex subsets of $\Sigma$. Given a semi-positive loop $\Lambda$ in $\mathbf{H}^{2,n}$, we define

- the *vertices*, points at the intersection of two segments of photon in $\Lambda$,
- the *edges*, maximal segments of photon contained in $\Lambda$ with non-empty interior,
- the *positive points*, points meeting no edge in $\Lambda$,
- the *positive intervals*, open intervals in $\Lambda$ containing no edge.

Using a description of the shape of horoballs, we define an increasing map in a subset of $\partial^{\mathrm{s}}\Sigma$:

$$\Phi : \{\text{positive points}\} \cup \{\text{vertices}\} \to \Sigma(\infty)\ .$$

**Remark 1.4.2.** The correspondence made between the space boundary and the ideal boundary is not as complete as the one in Barbot surfaces. We conjecture that the correspondence between the ideal boundary and the space boundary of polygonal surfaces is the same as in Barbot surfaces. In any case, the correspondence for vertices and positive points only is sufficient for our purpose.

Our construction is invariant by the action of $\mathsf{G}$, hence we can apply the renormalization process to understand the Tits distance in the ideal boundary of asymptotically flat maximal surfaces. The results are the following (proved in Theorem 4.7.1 and Theorem 4.7.2):

**Theorem 1.4.3.** *Let $\Sigma$ be an asymptotically flat maximal surface, and $v_1$, $v_2$ be adjacent vertices of $\partial^s\Sigma$. Then the Tits distance between the associated ideal points $\xi_{v_1}$ and $\xi_{v_2}$ equals $\pi/2$.*

**Theorem 1.4.4.** *Let $\Sigma$ be an asymptotically flat maximal surface, $I$ be a positive non-empty interval in $\partial^s\Sigma$, and $p, q$ be two distinct points of $I$. Then the Tits distance between the associated ideal points $\xi_p$ and $\xi_q$ equals $+\infty$.*

As a corollary, we understand the global geometry of the metric space $(\Sigma(\infty), \mathrm{Td})$. The following is proved in Theorem 4.8.5.

**Theorem 1.4.5.** *Let $\Sigma$ be a maximal surface. The Tits distance on the ideal boundary of $\Sigma$ has finite perimeter if and only if $\Sigma$ is a polygonal surface. A polygonal surface with $N+4$ vertices has Tits perimeter equal to $N\pi/2$.*

This gives the equivalence between items $(i)$ and $(ii)$ of the Main Theorem 1.2.1.

*1.4.6. Asymptotically flat maximal surfaces and the quartic metric*

Let now $\Sigma$ be an asymptotically flat maximal surface. Recall that the quartic metric and the induced metric are equal on Barbot surfaces. Hence for every positive number $\varepsilon$, there is a compact set $K$ in $\Sigma$ such that the induced metric $g$ and the quartic metric $g_4$ are $(1+\varepsilon)$-biLipschitz outside $K$ (Proposition 5.2.1). This makes $g_4$ a complete metric on $\Sigma$ with finitely many zeros. We prove in Theorem 5.2.4:

**Theorem 1.4.6.** *If a maximal surface $\Sigma$ is asymptotically flat, then the Riemann surface $(\Sigma, [g])$ is biholomorphic to the complex plane and the quartic differential $\varphi_4$ is a polynomial.*

We want to compare the quartic metric $g_4$ and the induced metric $g$. Since $g$ and $g_4$ are $(1+\varepsilon)$-biLipschitz outside a compact, they are $(1+\varepsilon, b)$-quasi isometric in $\Sigma$. This helps us to identify the $g$-circle $\mathbf{S}(r)$ and the $g_4$-circle $\mathbf{S}_4(r)$ of radii $r$, in a $(1+\varepsilon, b)$-quasi-isometric way. Taking the renormalized metric $d/r$ on $\mathbf{S}(r)$ and $d_4/r$ on $\mathbf{S}_4(r)$, and the limit when $r$ goes to $+\infty$, we obtain that the Tits boundary of $(\Sigma, g)$ and the Tits boundary of $(\Sigma, g_4)$ are isometric, yielding the following theorem (Theorem 5.5.1).

**Theorem 1.4.7.** *Let $\Sigma$ be an asymptotically flat maximal surface with polynomial quartic differential of degree $N$. Then $\Sigma$ is a polygonal surface with $N+4$ vertices.*

This is $(iii) \Rightarrow (i)$ of the Main Theorem, and concludes the scheme of the proof.

### 1.5. Organization of the paper

In Section 2, we recall the definitions and state some properties of $\mathbf{H}^{2,n}$ and $\mathbf{Ein}^{1,n}$. We define maximal surfaces and explain the work of Labourie Toulisse and Wolf about the asymptotic Plateau problem in $\mathbf{H}^{2,n}$. We also recall several definitions of the ideal boundary, and explain the construction of the space boundary as being the set of space-horofunctions on a maximal surface.

In Section 3 we define and study polygons in the Einstein universe. We emphasis the four vertices case, called Barbot crowns. We also describe a family of symmetric polygons. In the last part of this section, we prove the renormalization lemma, extensively used in what follows.

In Section 4, we study the link between the space boundary $\partial^{\mathrm{s}}\Sigma$ of a maximal surface $\Sigma$ and its ideal boundary $\Sigma(\infty)$. We begin with two examples: the Barbot surfaces and the maximal surfaces bounded by symmetric polygons. At the end we compute the Tits distance in terms of the space boundary allowing us to prove Theorem 4.8.5.

In Section 5, we compare the quartic metric associated with the holomorphic quartic differential and the induced metric. We prove that for an asymptotically flat maximal surface, the quartic differential is polynomial (Theorem 5.2.4) and the maximal surface is polygonal, relating the degree of the differential with the number of vertices (Theorem 5.5.1) by making a precise comparison between these two metrics.

*Acknowledgments*

I would like to thank Parker Evans, Enrico Trebeschi and Timothé Lemistre for useful discussions about polygons and maximal surfaces in pseudo-hyperbolic spaces, and Qiongling Li for discussions about finite total curvature and Higgs bundles. I thank Andrea Seppi for having pointed Corollary 1.3.1 as a consequence of the Main Theorem. I thank my advisors Jérémy Toulisse and François Labourie for having asked me the question at the origin of this article, and for their patience and encouragements during my research and the redaction of this article.

I also thank the anonymous referee for many useful comments that helped me to considerably improve the quality of the presentation.

## 2. Preliminaries

Once and for all, we fix $(\mathbf{E}, \mathbf{q})$ a real quadratic space of signature $(2, n+1)$ with a *complete orientation*, that is, a space orientation and a time orientation. The polar form of $\mathbf{q}$ is denoted by $\langle , \rangle$. We denote by

$$\mathbf{P} : \mathbf{E} \setminus \{0\} \to \mathbf{P}(\mathbf{E}),$$

or by $[.]$, the projectivization and by

$$[.]_+ : \mathbf{E} \setminus \{0\} \to \mathbf{E} \setminus \{0\}/\mathbb{R}_+^*$$

the quotient by the action of $\mathbb{R}_+^*$ by multiplication. We denote

$$\hat{\mathsf{G}} = \mathrm{O}(\mathbf{q})$$

the orthogonal group of $\mathbf{q}$ and

$$\mathsf{G} = \mathbf{P}\hat{\mathsf{G}} = \mathrm{O}(\mathbf{q})/\{\pm \mathrm{I}_{n+3}\}$$

its projectivized.

*2.1. The pseudo-hyperbolic space*

*2.1.1. Definition*

The pseudo-hyperbolic space is defined as follows

$$\mathbf{H}^{2,n} := \{x \in \mathbf{P}(\mathbf{E}) \mid \mathbf{q}_{|x} < 0\} \ ,$$

and its double cover

$$\widehat{\mathbf{H}}^{2,n} := \{x \in \mathbf{E} \mid \mathbf{q}(x) = -1\} \ .$$

Both spaces are endowed with a complete pseudo-riemannian metric of signature $(2, n)$ coming from the quadratic form $\mathbf{q}$ on $\mathbf{E}$ and are of constant sectional curvature -1.

The group $\hat{\mathsf{G}}$ acts naturally by isometries on $\widehat{\mathbf{H}}^{2,n}$ by multiplication on vectors. In the same manner, $\mathsf{G}$ acts isometrically on $\mathbf{H}^{2,n}$. Indeed, the groups $\mathsf{G}$ and $\hat{\mathsf{G}}$ are the full isometry groups of $\mathbf{H}^{2,n}$ and $\widehat{\mathbf{H}}^{2,n}$, respectively.

*2.1.2. Warped product structure*

An orthogonal decomposition $\mathbf{E} = P \oplus Q$ where $P$ is a positive definite plane and $Q$ is a negative definite $(n+1)$-subspace identifies $\widehat{\mathbf{H}}^{2,n}$ with a warped product as follows. Denote by $\mathbf{D}^2$ the unit disc of $(P, \mathbf{q}_{|P})$ and by $\mathbf{S}^n$ the unit sphere of $(Q, -\mathbf{q}_{|Q})$. By [9, Proposition 3.5], the map

$$\begin{aligned} \psi : \mathbf{D}^2 \times \mathbf{S}^n \to &\quad \widehat{\mathbf{H}}^{2,n} \\ (u, v) \mapsto &\ \frac{2}{1 - \mathbf{q}(u)} u \ + \ \frac{1 + \mathbf{q}(u)}{1 - \mathbf{q}(u)} v \end{aligned}$$

is a diffeomorphism onto and the pulled back metric $\psi^* \mathbf{q}$ equals $g_{\mathrm{hyp}} \oplus -f g_n$, where $g_{\mathrm{hyp}}$ is the hyperbolic metric on the Poincaré disc $\mathbf{D}^2$, $g_n$ is the round metric on the sphere $\mathbf{S}^n$ and $f : \mathbf{D}^2 \to \mathbb{R}$ is a smooth function defined by

$$f(u) = \left( \frac{1 + \mathbf{q}(u)}{1 - \mathbf{q}(u)} \right)^2 .$$

**Remark 2.1.1.** Denote by $\sigma$ the stereographic projection from the upper hemisphere $\mathbf{S}^2_+$ to the disc $\mathbf{D}^2$. We obtain a diffeomorphism $\phi = \psi(\sigma\cdot,\cdot)$ from $\mathbf{S}^2_+ \times \mathbf{S}^n$ to $\widehat{\mathbf{H}}^{2,n}$. The metric $\phi^*\mathbf{q}$ equals $(f\circ\sigma)\left(g_{\mathbf{S}^2_+} \oplus -g_n\right)$, where $g_{\mathbf{S}^2_+}$ denotes the round metric on the upper hemisphere. This is explained in [27, Section 3.2].

The image by $\psi$ of $\mathbf{D}^2 \times \{v_0\}$, for $v_0$ fixed, is totally geodesic and isometric to a hyperbolic plane. It corresponds to a connected component of $\widehat{\mathbf{H}}^{2,0}$ in the subspace $P + \mathbb{R}^+ v_0$ of signature $(2,1)$. The image by $\psi$ of $\{u_0\} \times \mathbf{S}^n$, for $u_0 \in \mathbf{D}^2$, is a totally umbilic submanifold called a *timelike sphere* and it is totally geodesic if and only if $u_0 = 0$.

*2.1.3. Warped projection*

The *warped projection* is, given a warped product $\psi : \mathbf{D}^2 \times \mathbf{S}^n \to \widehat{\mathbf{H}}^{2,n}$ and a base point $v_0$ in $\mathbf{S}^n$, the map

$$\begin{aligned} \pi_\psi : \widehat{\mathbf{H}}^{2,n} &\to \widehat{\mathbf{H}}^{2,n} \\ p &\mapsto \psi(\pi_1(\psi^{-1}(p))) \end{aligned}$$

where $\pi_1 : \mathbf{D}^2 \times \mathbf{S}^n \to \mathbf{D}^2 \times \mathbf{S}^n$ sends a point $(u,v)$ to $(u,v_0)$.

*2.1.4. Geodesics and the spacelike distance function*

Given two distinct points $x, y$ of $\mathbf{H}^{2,n}$, the set $\mathbf{P}(x \oplus y) \cap \mathbf{H}^{2,n}$ is an unparametrized geodesic. This geodesic is

- *spacelike* if the signature of $x \oplus y$ is $(1,1)$,
- *timelike* if the signature of $x \oplus y$ is $(0,2)$,
- *lightlike* if the signature of $x \oplus y$ is $(0,1)$.

A spacelike geodesic in $\widehat{\mathbf{H}}^{2,n}$ starting from a point $x$ and with direction $U$ a unit tangent vector is parametrized by arc length with the formula

$$t \mapsto \cosh(t)x + \sinh(t)U \ .$$

*2.1.5. Spacelike distance function*

**Definition 2.1.2** *(Spacelike distance).* Let $x$ and $y$ be two different points of $\mathbf{H}^{2,n}$ on a spacelike geodesic. Take $\hat{x}$ and $\hat{y}$ lifts of $x$ and $y$ to the double cover such that $\langle \hat{x}, \hat{y}\rangle < 0$. The *spacelike distance* between $x$ and $y$ is the number

$$\delta^{\mathrm{s}}(x,y) := \operatorname{arch}(-\langle \hat{x}, \hat{y}\rangle) \ .$$

**Remark 2.1.3.** The definition does not depend on the choices of lifts. Moreover, given two points $x, y$ on a spacelike geodesic, the spacelike distance $\delta^{\mathrm{s}}(x, y)$ equals the length of the $\mathbf{H}^{2,n}$-geodesic segment between them.

*2.1.6. The Einstein universe*

A reference for the Einstein universe is [2].

**Definition 2.1.4.** The *Einstein universe*, of signature $(1, n)$, is defined to be the projectivization of the isotropic cone of $(\mathbf{E}, \mathbf{q})$. We denote it $\mathbf{Ein}^{1,n}$. The double cover is defined as being the set of half lines in the isotropic cone of $(\mathbf{E}, \mathbf{q})$, and we denote it by $\widehat{\mathbf{Ein}}^{1,n}$.

The manifold $\mathbf{Ein}^{1,n}$ is the topological boundary of $\mathbf{H}^{2,n}$ seen as subsets of $\mathbf{P}(\mathbf{E})$.

The Einstein universe is endowed with a conformal class of Lorentzian metrics coming from $\mathbf{q}$. The manifold $\mathbf{Ein}^{1,n}$ is always space orientable, but time orientable (hence orientable) if and only if $n$ is odd. We fix once and for all a space orientation on $\mathbf{Ein}^{1,n}$.

The group $\hat{\mathsf{G}}$ acts naturally on $\widehat{\mathbf{Ein}}^{1,n}$ by multiplication on vectors $z$ belonging to the half lines $[z]_+$ in $\widehat{\mathbf{Ein}}^{1,n}$. The action of $\hat{\mathsf{G}}$ is conformal. In the same manner, $\mathsf{G}$ acts conformally on $\mathbf{Ein}^{1,n}$. Indeed, the groups $\mathsf{G}$ and $\hat{\mathsf{G}}$ are the groups of conformal automorphisms of $\mathbf{Ein}^{1,n}$ and $\widehat{\mathbf{Ein}}^{1,n}$, respectively.

*2.1.7. Extension of warped products*

A warped product decomposition of $\widehat{\mathbf{H}}^{2,n}$ as $\mathbf{S}^2_+ \times \mathbf{S}^n$ extends to a product description of $\widehat{\mathbf{Ein}}^{1,n}$ as being $\mathbf{S}^1 \times \mathbf{S}^n$ with the conformal class of Lorentzian metrics $[g_1 \oplus -g_n]$.

*2.1.8. Photons*

In general, two conformal metrics on a pseudo-riemannian manifold do not define the same geodesics. Even the images of the spacelike and timelike geodesics change, but not the images of lightlike geodesics [12, Proposition 2.131]. Lightlike geodesics are well-defined, as unparametrized curves, in $\mathbf{Ein}^{1,n}$, and we make the following definition.

**Definition 2.1.5.** A *photon* in $\mathbf{Ein}^{1,n}$ is a maximal lightlike geodesic (seen as an unparametrized curve).

**Proposition 2.1.6.** *A photon in* $\mathbf{Ein}^{1,n}$ *is the projectivization of an isotropic 2-plane.*

**Proof.** Identify $\widehat{\mathbf{Ein}}^{1,n}$ with $\mathbf{S}^1 \times \mathbf{S}^n$ by the warped product construction. Take the representative $g_1 \oplus -g_n$ of the conformal metric, and use the fact that a geodesic in a product projects to a geodesic in each factor. □

*2.2. Spacelike surfaces*

By *surface* we mean an embedded and connected submanifold of dimension 2.

*2.2.1. Graphs*

A surface $S$ in $\widehat{\mathbf{H}}^{2,n}$ is a *spacelike graph* if for any warped product decomposition $\mathbf{D}^2 \times \mathbf{S}^n \to \widehat{\mathbf{H}}^{2,n}$, the associated warped projections restrict to diffeomorphisms between $S$ and the hyperbolic plane. Equivalently, up to using the stereographic projection $\sigma$ from $\mathbf{S}^2_+$ to $\mathbf{D}^2$ (see Remark 2.1.1), a surface $S$ in $\widehat{\mathbf{H}}^{2,n}$ is a spacelike graph if and only if for any warped product decomposition $\mathbf{S}^2_+ \times \mathbf{S}^n \to \widehat{\mathbf{H}}^{2,n}$, the associated warped projections restrict to diffeomorphisms between $S$ and the upper hemisphere $\mathbf{S}^2_+$. We refer to [20, Propositions 3.6 and 3.10] and [27, Lemma 3.5] for the following proposition.

**Proposition 2.2.1.** *Let $\Sigma$ be a complete spacelike surface of $\widehat{\mathbf{H}}^{2,n}$. Then $\Sigma$ is a* spacelike graph *and in any warped product decomposition of $\widehat{\mathbf{H}}^{2,n}$ as $\mathbf{S}^2_+ \times \mathbf{S}^n$, the surface $\Sigma$ is the graph of a 1-Lipschitz map from $(\mathbf{S}^2_+, g_{\mathbf{S}^2_+})$ to $(\mathbf{S}^n, g_n)$.*

*Moreover, a complete spacelike surface of $\mathbf{H}^{2,n}$ has two disjoint diffeomorphic lifts in $\widehat{\mathbf{H}}^{2,n}$.*

*2.2.2. The topology of the set of complete spacelike surfaces*

We denote by $\mathcal{C}^\infty(\mathbf{D}^2, \mathbf{S}^n)$ the space of smooth functions from the unit disc to the $n$-sphere, endowed with the topology of uniform convergence on compact sets of all the derivatives.

By the previous proposition, given a warped product decomposition of $\widehat{\mathbf{H}}^{2,n}$, the set of complete spacelike surfaces of $\widehat{\mathbf{H}}^{2,n}$ embeds in $\mathcal{C}^\infty(\mathbf{D}^2, \mathbf{S}^n)$.

**Definition 2.2.2.** A sequence of complete spacelike surfaces $(\Sigma_k)_{k\in\mathbb{N}}$ of $\widehat{\mathbf{H}}^{2,n}$ *converges* to a complete spacelike surface $\Sigma$ if it holds under every warped product embedding in $\mathcal{C}^\infty(\mathbf{D}^2, \mathbf{S}^n)$. A sequence of complete spacelike surfaces of $\mathbf{H}^{2,n}$ *converges* if a sequence of lifts does.

*2.2.3. Semi-positive loops*

A *semi-positive loop*, as defined in [20, Definition 2.6], is a topological circle $\gamma$ in $\mathbf{Ein}^{1,n}$ satisfying the two following properties:

- Every triple of points in $\gamma$ spans a subspace of signature different than $(1, 2)$ in $\mathbf{E}$,
- There is at least one triple of points in $\gamma$ spanning a subspace of signature $(2, 1)$ in $\mathbf{E}$.

We use the same definition to define a semi-positive loop in $\widehat{\mathbf{Ein}}^{1,n}$.

**Proposition 2.2.3.**

(i) *[20, Proposition 2.11] A semi-positive loop of $\widehat{\mathbf{Ein}}^{1,n}$ is the graph of a 1-Lipschitz map from $(\mathbf{S}^1, g_1)$ to $(\mathbf{S}^n, g_n)$ in every warped product decomposition.*

(ii) *[20, Corollary 2.13] If a semi-positive loop contains two points that lie on a photon, the loop contains a segment of photon between them.*
(iii) *[20, Proposition 2.15] A semi-positive loop in $\mathbf{Ein}^{1,n}$ has two disjoint diffeomorphic lifts to $\widehat{\mathbf{Ein}}^{1,n}$. Moreover, if $\hat{p}$ and $\hat{q}$ are two points on the same lift of a semi-positive loop, then $\langle \hat{p}, \hat{q} \rangle \leq 0$.*

*2.2.4. The topology of the set of semi-positive loops*

We denote by $\mathrm{Lip}(\mathbf{S}^1, \mathbf{S}^n)$ the space of 1-Lipschitz functions from the unit circle to the $n$-sphere, endowed with the topology of uniform convergence.

By the previous proposition, given a warped product decomposition of $\widehat{\mathbf{H}}^{2,n}$, the set of semi-positive loops in $\widehat{\mathbf{Ein}}^{1,n}$ embeds in $\mathrm{Lip}(\mathbf{S}^1, \mathbf{S}^n)$.

**Definition 2.2.4.** A sequence of semi-positive loops in $\widehat{\mathbf{Ein}}^{1,n}$ *converges* if it holds under every warped product embedding in $\mathrm{Lip}(\mathbf{S}^1, \mathbf{S}^n)$. A sequence of semi-positive loops in $\mathbf{Ein}^{1,n}$ converges if a sequence of lifts does.

*2.2.5. Asymptotic Plateau problem*

The *space boundary* of a submanifold $M$ of $\mathbf{H}^{2,n}$ is the subset of $\mathbf{Ein}^{1,n}$ defined by

$$\partial^{\mathrm{s}} M = \left\{ p \in \mathbf{Ein}^{1,n} \mid p = \lim_{k\to\infty} x_k,\ (x_k)_{k\in\mathbb{N}} \in M^{\mathbb{N}} \right\} .$$

**Remark 2.2.5.** By Proposition 2.2.1, a complete spacelike surface in $\mathbf{H}^{2,n}$ has a space boundary homeomorphic to a circle. Moreover, by [20, Proposition 3.10], the space boundary of a complete spacelike surface is a semi-positive loop.

A complete spacelike surface in $\mathbf{H}^{2,n}$ is said to be *maximal* if its mean curvature vector vanishes. This is equivalent to being a local maximum for the area functional [20, Corollary 3.24].

In [20], Labourie, Toulisse and Wolf solve the Asymptotic Plateau problem in $\mathbf{H}^{2,n}$ as follows.

**Theorem 2.2.6** *(Asymptotic Plateau problem). A semi-positive loop in $\mathbf{Ein}^{1,n}$ is the space boundary of a unique maximal surface in $\mathbf{H}^{2,n}$.*

*2.2.6. Compactness theorems*

A *pointed maximal surface* is a pair $(\Sigma, x)$ where $\Sigma$ is a maximal surface and $x$ a point of $\Sigma$. The following is [20, Theorem 6.1] or [19, Theorem 4.8].

**Theorem 2.2.7.** *The group $\mathsf{G}$ acts properly, continuously and cocompactly on the space of pointed maximal surfaces.*

We refer to [20, Theorem 6.1] for the following.

**Theorem 2.2.8.** *Let $(\Sigma_k)_{k\in\mathbb{N}}$ be a sequence of maximal surfaces of $\mathbf{H}^{2,n}$. For each $k$, let $\Lambda_k$ be the space boundary of $\Sigma_k$. If the sequence $(\Lambda_k)_{k\in\mathbb{N}}$ converges to a semi-positive loop $\Lambda$, then the sequence $(\Sigma_k)_{k\in\mathbb{N}}$ converges to the maximal surface with space boundary $\Lambda$.*

## 2.3. The ideal boundary

We explain in this section what is the ideal boundary of a Hadamard manifold, using two equivalent definitions. We mimic the second definition to show how the space boundary can be thought of as an ideal boundary for the spacelike distance function.

Recall that a *Hadamard manifold* is a simply connected, complete, non-positively curved Riemannian manifold. A reference for the material about Hadamard manifolds and their ideal boundaries is [14].

### 2.3.1. The ideal boundary via rays

We define a *geodesic ray*, or simply a *ray*, in a complete Riemannian manifold as being a geodesic defined on $[0, +\infty)$ and parametrized by arclength.

**Definition 2.3.1.** Let $M$ be a Hadamard manifold. The *ideal boundary* of $M$, denoted $M(\infty)$, is defined as the quotient

$$M(\infty) := \{ \text{ rays in } M \ \}/\sim$$

where two rays $c_1$ and $c_2$ are equivalent ($c_1 \sim c_2$) if and only if there is a real number $r$ such that $d(c_1(t), c_2(t)) < r$ for all $t \geq 0$. If $c$ is a geodesic ray, we denote by $c(\infty)$ its class in $M(\infty)$.

If $M$ is a complete Riemannian manifold and $X$ is a tangent vector of $M$ at the point $x$, then we denote by $c_X$ the ray starting at $x$ with initial velocity $X$. The following is classical.

**Proposition 2.3.2.** *Let $M$ be a Hadamard manifold with ideal boundary $M(\infty)$. For every $\xi$ in $M(\infty)$ and every $x$ in $M$, there is a unique ray starting at $x$ representing the ideal boundary point $\xi$.*

So, given any point $x$ of $M$ there is a bijection from the unit tangent bundle $T^1_x M$ to the ideal boundary $M(\infty)$. We can give a topology to $M(\infty)$ via this bijection. It turns out that this topology does not depend on the chosen point $x$.

An ideal boundary point $\xi$ of a Hadamard manifold can be seen as a collection of rays, one for each point of the manifold. We denote by $\xi_x$ the representative of an ideal boundary point $\xi$ at a point $x$.

### 2.3.2. The ideal boundary via horofunctions

Another equivalent way to define the ideal boundary is via *horofunctions*. Let $M$ be a complete Riemannian manifold. Denote by $\mathcal{C}(M)$ the space of continuous functions on $M$, with the topology of uniform convergence on compact sets, that makes it a complete space. Denote by $\mathcal{C}^*(M)$ the quotient of this space by the equivalence relation $f \sim g$ if $f = g + r$ for $r$ a real number. Then $M$ embeds in $\mathcal{C}^*(M)$ via distance functions:

$$\begin{aligned}\psi : M &\to \mathcal{C}^*(M)\\ x &\mapsto [d(x,\cdot)].\end{aligned}$$

**Definition 2.3.3.** A *horofunction* on $M$ is an element of the frontier of $\psi(M)$ in $\mathcal{C}^*(M)$. The *horo boundary* of $M$, denoted $\partial^h M$, is the set of horofunctions.

It is shown in [14, 3.6] that, if $M$ is a Hadamard manifold, the ideal boundary $M(\infty)$ and the horo boundary $\partial^h M$ are homeomorphic, using the following identification.

Given an ideal boundary point $\xi$ in $M(\infty)$, we can look at a limit of renormalized distance functions as follows: take a base point $x_0$ in $M$, denote by $c$ the ray starting at $x_0$ with direction $\xi$, and define the curve

$$\begin{aligned}[0,+\infty) &\to \mathcal{C}(M)\\ t &\mapsto d(c(t),\cdot) - d(c(t),x_0).\end{aligned}$$

This curve has a limit when $t$ goes to $+\infty$, and this limit cannot be a distance function, hence it is a horofunction. One can show that the class of the limit function in $\mathcal{C}^*(M)$ does not depend on the choice of the base point, but only on the direction $\xi$. We call these functions *Busemann functions.* It turns out that, conversely, all horofunctions are Busemann functions for a certain ideal boundary point.

*2.3.3. The space boundary of complete spacelike surfaces*

Let $\Sigma$ be a complete spacelike surface of $\mathbf{H}^{2,n}$. Recall that its space boundary $\partial^{\mathrm{s}}\Sigma$ is defined as being the subset of $\mathbf{Ein}^{1,n}$ containing the limits of diverging sequences in $\Sigma$. We now present another point of view, imitating the second construction of the ideal boundary. As in the previous section, $\mathcal{C}(\Sigma)$ is the space of continuous functions on $\Sigma$ endowed with the topology of uniform convergence on compact sets and $\mathcal{C}^*(\Sigma)$ is the quotient of this space by constant functions.

Consider:

$$\begin{aligned}\psi^s : \Sigma &\to \mathcal{C}^*(\Sigma)\\ x &\mapsto [\delta^{\mathrm{s}}(x,\cdot)]\ ,\end{aligned}$$

where $\delta^{\mathrm{s}}$ is the spacelike distance function (recall Definition 2.1.2).

**Proposition 2.3.4.** *Let $\Sigma$ be a complete spacelike surface in $\mathbf{H}^{2,n}$ and $x$ be a point of $\Sigma$. Then the spacelike distance function $\delta^{\mathrm{s}}(x,\cdot)$ is well defined on $\Sigma$. In particular, $\delta^{\mathrm{s}}(x,\cdot)$ has a unique minimum reached at the point $x$ and this minimum is strict.*

**Proof.** By [20, Proposition 3.10], since $\Sigma$ is complete it is acausal (item $(i)$ and $(ii)$ of the proposition), hence for every point $y$ in $\Sigma$ the geodesic segment from $x$ to $y$, in $\mathbf{H}^{2,n}$, is spacelike. Moreover $\delta^{\mathrm{s}}(x,y)=0$ if and only if $x=y$. □

**Proposition 2.3.5.** *The map $\psi^s$ is a continuous embedding.*

**Proof.** Since $\Sigma$ is a spacelike surface, the map $\psi^s$ is well-defined. The map $\psi^s$ is also injective, since for every point $x$ of $\Sigma$ the minimum of each representative of $\psi^s(x)$ is obtained at $x$ only. We now take a lift $\hat{\Sigma}$ of $\Sigma$ to $\widehat{\mathbf{H}}^{2,n}$. The spacelike distance function $\delta^{\mathrm{s}}$ can be expressed in $\hat{\Sigma}$ as

$$\delta^{\mathrm{s}}(x_1,x_2)=\mathrm{arcch}(\langle x_1,x_2\rangle)\ .$$

Since $\hat{\Sigma}$ is a smooth submanifold of $\mathbf{E}$, and $\langle\cdot,\cdot\rangle$ a smooth function, $\psi^s$ is continuous. □

As for the ideal boundary, we want to characterize the frontier of this embedding. Before making computations, we recall that a complete spacelike surface of $\mathbf{H}^{2,n}$ admits two disjoint diffeomorphic lifts to the double cover (Proposition 2.2.3(iii)).

**Definition 2.3.6.** The *space-horofunction* associated to the point $p$ of $\widehat{\mathbf{Ein}}^{1,n}$ is the family of functions defined on $\widehat{\mathbf{H}}^{2,n}\setminus p^{\perp}$, up to an additive constant, by

$$h_p(x)=\log(|\langle z,x\rangle|)\ ,$$

for $z$ in $p$.

Here $p^{\perp}$ denotes the set of vectors orthogonal to $p$ for the quadratic form $\mathbf{q}$. Remark that $h_p=h_{-p}$, and $h_p(-x)=h_p(x)$ for any $p$ in $\widehat{\mathbf{Ein}}^{1,n}$ and $x$ in $\widehat{\mathbf{H}}^{2,n}$. Hence we also define the *space-horofunction* associated to a point $p$ in $\mathbf{Ein}^{1,n}$, and denote it by $h_p$, with the formula

$$\begin{aligned}h_p:\mathbf{H}^{2,n}\setminus L(p)&\to\mathbb{R}\\ x&\mapsto h_{\hat{p}}(\hat{x})\ ,\end{aligned}$$

where $\hat{p}$ and $\hat{x}$ are lifts of $p$ and $x$ to the double covers $\widehat{\mathbf{Ein}}^{1,n}$ and $\widehat{\mathbf{H}}^{2,n}$, respectively, and $L(p)$ is the set of points of $\mathbf{H}^{2,n}$ related to $p$ by a lightlike geodesic.

**Theorem 2.3.7.** *The frontier of $\psi^s(\Sigma)$ in $\mathcal{C}^*(\Sigma)$ is the set $\{h_p \mid p\in\partial^{\mathrm{s}}\Sigma\}$ of space-horofunctions associated with points of $\partial^{\mathrm{s}}\Sigma$. This defines a homeomorphism $\Sigma\cup\partial^{\mathrm{s}}\Sigma\to\overline{\psi^s(\Sigma)}$.*

We make the computations in the double cover.

**Proof.** The set $\{h_p \mid p \in \partial^{\mathrm{s}}\Sigma\}$ is contained in the frontier of $\psi^s(\Sigma)$ as follows.

Let $p$ be a point in $\partial^{\mathrm{s}}\Sigma$ and denote by $h_p$ the associated space-horofunction. Let $(x_k)_{k\in\mathbb{N}}$ be a sequence of points in $\Sigma$ that tends to $p$. Take a vector $z$ in $p$ and $\hat{x}_k$ a lift of $x_k$ to the double cover such that $\langle \hat{x}_k, z\rangle$ is negative for each integer $k$. Take a sequence $(\lambda_k)_{k\in\mathbb{N}}$ of positive real numbers such that $(\lambda_k\hat{x}_k)_{k\in\mathbb{N}}$ tends to $z$. Since the $\hat{x}_k$'s are in $\{\mathbf{q} = -1\}$, the sequence $(\lambda_k)_{k\in\mathbb{N}}$ tends to 0. By definition we have

$$f_k(x) := \delta^{\mathrm{s}}(x_k, \cdot) - \delta^{\mathrm{s}}(x_k, x_0) = \operatorname{arccosh}\left(-\langle \hat{x}_k, \hat{\cdot}\rangle\right) - \operatorname{arccosh}\left(-\langle \hat{x}_k, \hat{x}_0\rangle\right) \ .$$

Since $\operatorname{arccosh}(t) = \log(t + \sqrt{t^2 - 1})$ for every $t \geq 1$, we have

$$\begin{aligned}\exp(f_k(x)) = \exp(\delta^{\mathrm{s}}(x_k, x) - \delta^{\mathrm{s}}(x_k, x_0)) &= \frac{-\langle \hat{x}_k, \hat{x}\rangle + \sqrt{\langle \hat{x}_k, \hat{x}\rangle^2 - 1}}{-\langle \hat{x}_k, \hat{x}_0\rangle + \sqrt{\langle \hat{x}_k, \hat{x}_0\rangle^2 - 1}} \\ &= \frac{\langle \lambda_k\hat{x}_k, \hat{x}\rangle - \sqrt{\langle \lambda_k\hat{x}_k, \hat{x}\rangle^2 - \lambda_k^2}}{\langle \lambda_k\hat{x}_k, \hat{x}_0\rangle - \sqrt{\langle \lambda_k\hat{x}_k, \hat{x}_0\rangle^2 - \lambda_k^2}} \ .\end{aligned}$$

Hence, the sequence $(\exp(f_k))_{k\in\mathbb{N}}$ converges uniformly on the compact sets to the function

$$x \mapsto \frac{\langle z, \hat{x}\rangle}{\langle z, \hat{x}_0\rangle} \ .$$

That is,

$$\lim_{k\to+\infty} \psi^s(x_k) = h_p \ .$$

We extend $\psi^s$ to $\partial^{\mathrm{s}}\Sigma$ by $\psi^s(p) = h_p$.

We claim that $\psi^s$ is surjective from $\partial^{\mathrm{s}}\Sigma$ to $\mathrm{Fr}(\psi^s(\partial^{\mathrm{s}}\Sigma))$. Indeed, let $h$ be a function in $\mathrm{Fr}(\psi^s(\partial^{\mathrm{s}}\Sigma))$. There is a sequence $(x_k)_{k\in\mathbb{N}}$ of points in $\Sigma$ such that $(\psi^s(x_k))_{k\in\mathbb{N}}$ converges to $h$. Now, since $h$ is not an element of $\psi^s(\Sigma)$, the sequence $(x_k)_{k\in\mathbb{N}}$ diverges in $\Sigma$. The sequence $(x_k)_{k\in\mathbb{N}}$ converges to a point $p$ of $\partial^{\mathrm{s}}\Sigma$, hence $(\psi^s(x_k))_{k\in\mathbb{N}}$ subconverges to $\psi^s(p)$. Finally, $h = \psi^s(p)$ and $\psi^s$ is surjective.

For the last assertion, since $\delta^{\mathrm{s}}(x, \cdot)$ and $h_p$ are explicit, a computation as above shows that for every sequence $(x_k)_{k\in\mathbb{N}}$ in $\Sigma \cup \partial^{\mathrm{s}}\Sigma$ and $x$ in $\Sigma \cup \partial^{\mathrm{s}}\Sigma$, the sequence $(x_k)_{k\in\mathbb{N}}$ converges to $x$ if and only if the sequence $(\psi^s(x_k))_{k\in\mathbb{N}}$ converges to $\psi^s(x)$. □

**Remark 2.3.8.** This construction gives an extrinsic ideal boundary to $\Sigma$, even if the ambient space is not Riemannian. It can be generalized in various pseudo-Riemannian manifolds where complete non-compact acausal submanifolds can be endowed with a space boundary in the same manner.

## 3. Polygons

The aim of this section is to define the *lightlike polygons* in the Einstein universe, and study some of their properties.

### *3.1. Polygons in* $\mathbf{Ein}^{1,n}$

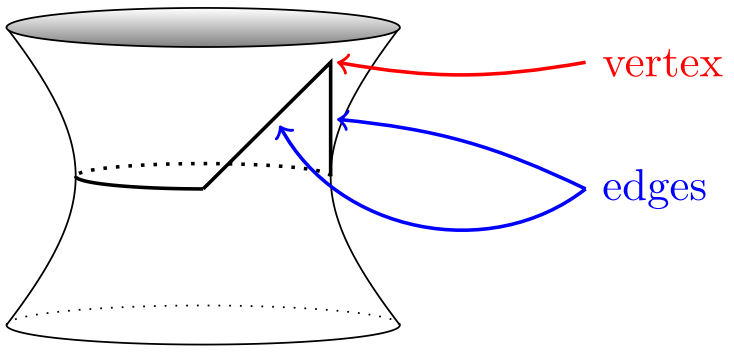

**Fig. 1.** A semi-positive loop in the boundary of $\mathbf{H}^{2,1}$ represented in an affine chart, with one vertex.

The notion of segment we have in $\mathbf{Ein}^{1,n}$ is that of segment of photons. Recall (Definition 2.1.5 and Proposition 2.1.6) that a photon is the image of a complete lightlike geodesic and is indeed the projectivization of an isotropic plane of $\mathbf{E}$.

**Definition 3.1.1** *(Lightlike polygons).* A *lightlike polygon* in $\mathbf{Ein}^{1,n}$ is a semi-positive loop consisting of finitely many segments of photons. We write a $N$-gon to designate a lightlike polygon with $N$ vertices.

**Definition 3.1.2.** Let $\Lambda$ be a semi-positive loop in $\mathbf{Ein}^{1,n}$. A *vertex* is a point of $\Lambda$ having a neighborhood, in $\Lambda$, consisting of the union of two segments of different photons. (See Fig. 1.)

**Definition 3.1.3.** Let $\Lambda$ be a semi-positive loop in $\mathbf{Ein}^{1,n}$. An *edge* of $\Lambda$ is a maximal segment of photon included in $\Lambda$ and containing at least two points. (See Fig. 1.)

**Proposition 3.1.4.** *Every lightlike polygon has at least four vertices.*

**Proof.** Given two points $p, q$ in $\mathbf{Ein}^{1,n}$, there is at most 1 photon passing through them. Hence there is no 2-gon. If a semi-positive loop in $\mathbf{Ein}^{1,n}$ contains three consecutive vertices $v_1$, $v_2$ and $v_3$, then the subspace of $\mathbf{E}$ spanned by the lines $v_1$, $v_2$ and $v_3$ is 3-dimensional. Now, if $v_3$ and $v_1$ are $\mathbf{q}$-orthogonal lines, the lines $v_1$, $v_2$ and $v_3$ span a 3-dimensional totally isotropic plane in $\mathbf{E}$ and this is impossible. $\square$

**Remark 3.1.5.** We can find loops in $\widehat{\mathbf{Ein}}^{1,n}$ that are the union of two segments of photon, but they are not semi-positive. We build an example in $\widehat{\mathbf{Ein}}^{1,1}$ that can be embedded in $\widehat{\mathbf{Ein}}^{1,n}$. Take $p$ in $\widehat{\mathbf{Ein}}^{1,1}$ and $q = -p$. There are four segments of photon between $p$ and $q$. Indeed there are two photons $\phi$ and $\psi$ passing through $[p]$ in $\mathbf{Ein}^{1,1}$ hence four lifts $\phi_+, \phi_-, \psi_+, \psi_-$ in the double cover. If $\phi_+$ and $\psi_+$ contain points with positive scalar product, then $\phi_+ \cup \psi_-$ is an example. Otherwise $\phi_+ \cup \psi_+$ is an example. Every loop consisting of two segments of photon in $\widehat{\mathbf{Ein}}^{1,n}$ is of this form. The projection of the obtained loop in $\mathbf{Ein}^{1,n}$ is the union of the two photons $\phi \cup \psi$.

The first interesting case is the object of the next section.

### *3.2. Barbot crowns*

Following [20], we introduce the following terminology.

**Definition 3.2.1.** A *Barbot crown* is a lightlike polygon with 4 vertices.

#### *3.2.1. The span of a Barbot crown*

A Barbot crown is always included inside a copy of $\mathbf{Ein}^{1,1}$. Indeed, take a Barbot crown $C$ with vertices $v_1$, $v_2$, $v_3$, $v_4$. Take a lift of $C$ in $\widehat{\mathbf{Ein}}^{1,n}$ with vertices $\hat{v}_i$, and let $z_i$ in $\hat{v}_i$ be non zero vectors. Two neighboring vertices lie on a photon, so $\langle z_i, z_{i+1}\rangle = 0$. Two vertices that are not adjacent are two non-orthogonal lines. Indeed, if two points of a semi-positive loop lie on a photon, the loop must contain a segment of photon between them (Proposition 2.2.3(ii)).

Hence, we have

$$(\langle z_i, z_j\rangle)_{1\le i,j\le 4} = \begin{pmatrix} 0 & 0 & \langle z_1, z_3\rangle & 0 \\ 0 & 0 & 0 & \langle z_2, z_4\rangle \\ \langle z_1, z_3\rangle & 0 & 0 & 0 \\ 0 & \langle z_2, z_4\rangle & 0 & 0 \end{pmatrix},$$

with determinant equal to $\langle z_1, z_3\rangle^2 \langle z_2, z_4\rangle^2 > 0$. Therefore, the vectors $\{z_i\}$ span a non-degenerate subspace of $\mathbf{E}$, containing totally isotropic planes. Hence $v_1 \oplus v_2 \oplus v_3 \oplus v_4$ is of signature $(2,2)$ and $C$ is included in the projectivization $\mathbf{P}(v_1 \oplus v_2 \oplus v_3 \oplus v_4) \cap \mathbf{Ein}^{1,n} \simeq \mathbf{Ein}^{1,1}$.

#### *3.2.2. Hyperbolic basis*

**Definition 3.2.2.** We say that a quadruple $(z_1, z_2, z_3, z_4)$ of distinct isotropic vectors of $\mathbf{E}$ is a *hyperbolic basis* if

- $\mathbb{R}z_i \oplus \mathbb{R}z_{i+1}$ is totally isotropic for each $1 \le i \le 4$ (with indices mod 4),
- $\langle z_i, z_{i+2}\rangle = -1/4$ for $i = 1$ and $i = 2$.

**Proposition 3.2.3.** *Given the vertices $v_1$, $v_2$, $v_3$ and $v_4$ of a Barbot crown, there exist nonzero vectors $z_1$, $z_2$, $z_3$ and $z_4$ in $v_1$, $v_2$, $v_3$ and $v_4$ respectively that form a hyperbolic basis.*

**Proof.** Take $w_1$, $w_2$, $w_3$ and $w_4$ nonzero vectors in $v_1$, $v_2$, $v_3$ and $v_4$, respectively. The first property is automatically satisfied and $\langle w_i, w_{i+2}\rangle \neq 0$. The vectors

$$z_1 = w_1,\ z_2 = w_2,\ z_3 = -\frac{1}{4\langle w_1, w_3\rangle} w_3 \text{ and } z_4 = -\frac{1}{4\langle w_2, w_4\rangle} w_4$$

then form a hyperbolic basis. □

**Proposition 3.2.4.** *A hyperbolic basis $(z_i)_{1\le i\le 4}$ determines in a unique way a Barbot crown in $\widehat{\mathbf{Ein}}^{1,n}$ passing through the vertices $\hat{v}_i = [z_i]_+$, hence a unique Barbot crown in $\mathbf{Ein}^{1,n}$ passing through the vertices $v_i = \mathbf{P}(z_i)$.*

**Proof.** For each $i$ in $\mathbb{Z}/4\mathbb{Z}$, the segment $\{(1-t)z_i + tz_{i+1} \mid t \in [0,1]\}$ gives rise to a segment of photon between $\hat{v}_i$ and $\hat{v}_{i+1}$, in $\widehat{\mathbf{Ein}}^{1,n}$. The union of these four segments is a topological circle $\hat{\Lambda}$, and is semi-positive. Indeed, the inequality

$$\langle (1-t)z_i + tz_{i+1}, (1-s)z_j + sz_{j+1}\rangle \le 0\ ,$$

for all $s, t$ in $(0,1)$ and $i, j$ in $\{1, \dots, 4\}$ proves that triples of points of $\Lambda$ cannot span a space of signature $(1,2)$ and that a triple of $\Lambda$ not containing two points on a photon span a space of signature $(2,1)$. Hence we have defined a Barbot crown $C$ in $\widehat{\mathbf{Ein}}^{1,n}$, that projects to a Barbot crown in $\mathbf{Ein}^{1,n}$. The Barbot crown $C$ in $\widehat{\mathbf{Ein}}^{1,n}$ is unique. Indeed, we have $\langle z_i, (1-t)z_j + tz_{j+1}\rangle > 0$ for $\{t < 0$ and $i = j-1\}$ and for $\{t > 1$ and $i = j+2\}$ (recall that two points $x, y$ on a semi-positive loop must satisfy $\langle x, y\rangle \le 0$ by Proposition 2.2.3(iii)). □

**Remark 3.2.5.** If $(z_i)_{1\le i\le 4}$ is a hyperbolic basis, then $(-z_i)_{1\le i\le 4}$ is also a hyperbolic basis. The Barbot crowns passing through $([z_i]_+)_{1\le i\le 4}$ and $(-[z_i]_+)_{1\le i\le 4}$, respectively, are in bijection via the map $-id_{|\mathbf{E}}$. Moreover, $(z_1, -z_2, z_3, -z_4)$ is also a hyperbolic basis, and the Barbot crown passing through $(z_1, -z_2, z_3, -z_4)$ is different than the previous Barbot crowns. Hence there are two distinct Barbot crowns in $\mathbf{Ein}^{1,n}$ with vertices $[z_1]$, $[z_2]$, $[z_3]$ and $[z_4]$: the projections of the Barbot crowns passing through $([z_1]_+, [z_2]_+, [z_3]_+, [z_4]_+)$ and $([z_1]_+, [-z_2]_+, [z_3]_+, [-z_4]_+)$ in $\widehat{\mathbf{Ein}}^{1,n}$, respectively.

### *3.2.3. Uniqueness of Barbot crowns under the action of* $\mathsf{G}$

**Proposition 3.2.6.** *The group $\hat{\mathsf{G}}$ acts transitively on the set of hyperbolic basis. The group $\mathsf{G}$ acts transitively on the set of Barbot crowns.*

**Proof.** Fix a basis $e_1, \dots e_{n+1}$ of $\mathbf{E}$ where $\mathbf{q}$ is written

$$\mathbf{q}(x) = -x^1x^3 - x^2x^4 - \sum_{5\le i\le n+3} (x^i)^2\ ,$$

for $x = \sum_i x^i e_i$. Let $(z_1, z_2, z_3, z_4)$ be a hyperbolic basis. The map sending each $z_i$ to $\frac{1}{2}e_i$ respects the quadratic form $\mathbf{q}$ and can therefor be extended to an isometry of $\mathbf{E}$, giving an element of $\hat{\mathsf{G}}$ sending the hyperbolic basis to $(\frac{1}{2}e_1, \frac{1}{2}e_2, \frac{1}{2}e_3, \frac{1}{2}e_4)$.

Let $v_1$, $v_2$, $v_3$ and $v_4$ be the ordered vertices of a Barbot crown $C$ in $\mathbf{Ein}^{1,n}$. Take a lift $\hat{C}$ of $C$ to the double cover $\widehat{\mathbf{Ein}}^{1,n}$ with $z_i$ in $\hat{v}_i$ non-zero vectors. Then, up to the multiplication of each $z_i$ by a positive constant, $(z_1, z_2, z_3, z_4)$ is a hyperbolic basis. Since $\hat{\mathsf{G}}$ acts transitively on the set of hyperbolic basis, we obtain the result. $\square$

*3.2.4. Determination of a Barbot crown*

We describe here a way to determine a unique Barbot crown, given three points in a certain configuration.

**Proposition 3.2.7.** *Given three points $v_1$, $v_2$ and $v_3$ of $\mathbf{Ein}^{1,n}$ satisfying the following properties:*

- $\mathbf{P}(v_1 \oplus v_2)$ *and* $\mathbf{P}(v_2 \oplus v_3)$ *are photons,*
- $v_1 \oplus v_3$ *is of signature* $(1, 1)$ *for* $\mathbf{q}$,

*there is a point $v_4$ in $\mathbf{Ein}^{1,n}$ such that the points $v_i$ are the vertices of a Barbot crown.*

If $n$ equals 1, the point $v_4$ is unique. If $n$ is bigger than 1, then the point is not unique, there is a $n-1$-dimensional manifold of possible choices.

**Proof.** We search a vector in $v_1^\perp \cap v_3^\perp$ non-orthogonal to $v_2$. The subspace

$$v_1^\perp \cap v_3^\perp = (v_1 \oplus v_3)^\perp$$

is of dimension $n+1$ and signature $(1, n)$ for $\mathbf{q}$, and contains $v_2$. Hence the result. $\square$

**Corollary 3.2.8.** *Let $\Lambda$ be a semi-positive loop containing a vertex $v$. Then there is a Barbot crown that coincides with $\Lambda$ along the edges passing through $v$.*

Recall that, given the vertices $v_1$, $v_2$, $v_3$ and $v_4$ of a Barbot crown, there are two distinct Barbot crowns of $\mathbf{Ein}^{1,n}$ passing through these points (Remark 3.2.5).

**Proof.** We denote by $v_2$ the point $v$ and by $v_1$, $v_3$ the two points at the end of the edges, in $\Lambda$, passing through $v_2$. Then we can apply Proposition 3.2.7. There exists $v_4$ in $\mathbf{Ein}^{1,n}$ such that the points $v_i$ are the vertices of a Barbot crown. Since there are two Barbot crowns passing through the points $v_i$ it remains to show that one of them works.

Take a lift of $\Lambda$ to the double cover with $[z_1]_+$, $[z_2]_+$ and $[z_3]_+$ lifting $v_1$, $v_2$ and $v_3$ respectively. Denote by $z_4$ a nonzero vector in $v_4$ such that the $[z_i]_+$ are the vertices of a lift $\hat{C}$ of $C$. Now, take the lift $[z_4]_+$ of $v_4$ such that $\langle z_2, z_4 \rangle < 0$. The Barbot crown in $\widehat{\mathbf{Ein}}^{1,n}$ with vertices $[z_1]_+$, $[z_2]_+$, $[z_3]_+$ and $[z_4]_+$ projects in $\mathbf{Ein}^{1,n}$ to a Barbot crown that coincides with $\Lambda$ along the edges passing through $v$. $\square$

*3.2.5. Barbot crowns and Cartan subgroup*

Take a Barbot crown $C$ whose vertices are $v_1$, $v_2$, $v_3$ and $v_4$. Take $z_i$ in $v_i$ for each $1 \leq i \leq 4$ such that $(z_1, z_2, z_3, z_4)$ is a hyperbolic basis. Define $\mathsf{A}$ as being the subgroup of $\hat{\mathsf{G}}$ consisting of elements stabilizing each vertex $v_i$ and stabilizing the orthogonal complement $(v_1 \oplus \ldots \oplus v_4)^\perp$ point wise. The group $\mathsf{A}$ can also be seen as a subgroup of $\mathsf{G}$. Let $(e_5, \ldots, e_{n+3})$ be an orthonormal basis of $(v_1 \oplus v_2 \oplus v_3 \oplus v_4)^\perp$. Hence $(z_1, \ldots, z_4, e_5, \ldots, e_{n+3})$ is a basis of $\mathbf{E}$ and we write elements of $\hat{\mathsf{G}}$ with respect to this basis.

**Definition 3.2.9.** We call $\mathsf{A}_0$ the *Cartan subgroup associated to* $C$ and denote its elements by

$$a(\lambda, \mu) := \mathrm{diag}(\lambda, \mu, \lambda^{-1}, \mu^{-1}, 1 \ldots, 1) \ ,$$

for $\lambda, \mu$ real numbers.

## 3.3. Symmetric polygons

We construct now symmetric polygons in $\widehat{\mathbf{Ein}}^{1,1}$. They are nontrivial examples of polygons and the symmetry will translate into geometric properties of the associated maximal surface. By *symmetric*, we mean that the isometry group of $\mathbf{H}^{2,1}$ acts as the dihedral group on the polygon.

*3.3.1. Construction*

Let $N$ be an even integer bigger than 4.

We use a warped product description of $\widehat{\mathbf{Ein}}^{1,1}$ as being $\mathbf{S}^1 \times \mathbf{S}^1$. Let $(u_k)_{k=1}^N$ and $(v_k)_{k=1}^N$ be two collections of points of the circle $\mathbf{S}^1$, defined by

$$u_k := \exp(2ik\pi/N + i\pi/N), \text{ and } v_k := \begin{cases} \exp(i\pi/N) & \text{if } k \text{ is odd,} \\ \exp(-i\pi/N) & \text{if } k \text{ is even.} \end{cases}$$

We obtain a collection of $N$ points $(p_k)_{k=1}^N = ((u_k, v_k))_{k=1}^N$ in $\widehat{\mathbf{Ein}}^{1,1}$. These points are the vertices of a lightlike polygon $P$ (see Fig. 2).

*3.3.2. Stability under the dihedral group*

The polygon $P$ is stabilized by the two elements $\mathbf{r}$ and $\mathbf{s}$, defined by their action on $\mathbf{S}^1 \times \mathbf{S}^1$ as follows

$$\mathbf{r} : (u, v) \mapsto (\exp(2i\pi/N)u, \bar{v}) \quad \text{and} \quad \mathbf{s} : (u, v) \mapsto (\bar{u}, \bar{v}) \ ,$$

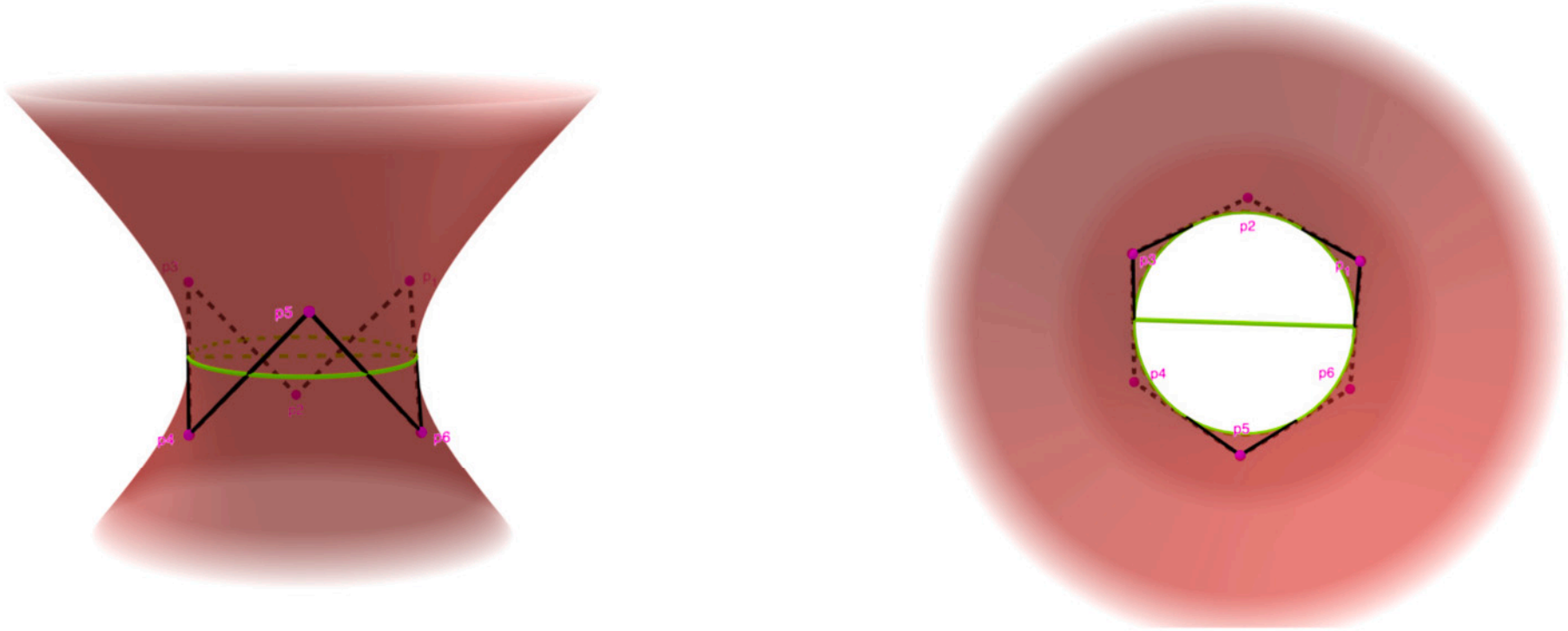

**Fig. 2.** The symmetric polygon $P$ with 6 vertices (in black) seen in an affine chart from the side (picture on the left) and from above (picture on the right). In this chart, the maps $\mathbf{r}$ and $\mathbf{s}$ are described as follows. The rotation $\mathbf{r}$ is the rotation of the green circle of angle $2\pi/6$ followed by a symmetry fixing the plane where the green circle lies. The symmetry $\mathbf{s}$ is the symmetry fixing the green segment. (For interpretation of the colors in the figure(s), the reader is referred to the web version of this article.)

where $\overline{\exp(i\theta)} := \exp(-i\theta)$ for every angle $\theta$. The subgroup of $\mathsf{O}(2,2)$ generated by $\mathbf{r}$ and $\mathbf{s}$ is isomorphic to the dihedral group with $2N$ elements and for every $k$ in $\{1,\dots,N\}$ we have

$$\mathbf{r}\cdot p_k = p_{k+1}, \text{ and } \mathbf{s}\cdot p_k = p_{-k} \ ,$$

where the indices are given modulo $N$. Again see Fig. 2 for a visual description of the maps $\mathbf{r}$ and $\mathbf{s}$.

### *3.3.3. Higher codimension*

The constructed symmetric polygons lie in $\mathbf{Ein}^{1,1}$. In $\mathbf{Ein}^{1,2}$, for an integer $N$ bigger than 4, the configuration space $\mathcal{C}_N$ of polygons with $N$ vertices is not connected. Together with Parker Evans, in an ongoing project, we classify the symmetric polygons of $\mathcal{C}_N$ in a similar way to the above example.

## *3.4. Renormalization of polygons*

We prove a lemma that we use to renormalize semi-positive loops with a vertex into Barbot crowns.

Let $C$ be a Barbot crown in $\mathbf{Ein}^{1,n}$ with vertices $v_1$, $v_2$, $v_3$ and $v_4$. Denote by $\phi$ the edge between $v_1$ and $v_2$ in $C$. Denote by $\psi$ the edge between $v_2$ and $v_3$ in $C$. Let $\mathsf{A}_0$ be the Cartan subgroup of $\mathsf{G}$ associated with the Barbot crown $C$ (see the previous paragraph for a definition of $\mathsf{A}_0$).

**Lemma 3.4.1** *(Renormalization lemma near a vertex). Let $\Lambda$ be a semi-positive loop in $\mathbf{Ein}^{1,n}$ having $\phi$ and $\psi$ as edges. Let $(a_k)_{k\in\mathbb{N}}$ be a sequence of $\mathsf{A}_0$, the Cartan subgroup associated to $C$, where $a_k = a(\lambda_k,\mu_k)$ satisfies the conditions*

*(1) the sequence* $(\mu_k)_{k\in\mathbb{N}}$ *tends to 0,*
*(2) the sequence* $(\lambda_k\mu_k)_{k\in\mathbb{N}}$ *is bounded from above,*
*(3) the sequence* $(\lambda_k^{-1}\mu_k)_{k\in\mathbb{N}}$ *is bounded from above.*

*Then the sequence of semi-positive loops* $(a_k\cdot\Lambda)_{k\in\mathbb{N}}$ *converges to* $C$*.*

Take $\hat{\Lambda}$ and $\hat{C}$ lifts of $\Lambda$ and $C$, respectively, that coincide above the edges $\phi$ and $\psi$. Consider $(z_1,z_2,z_3,z_4)$ a hyperbolic basis (recall Definition 3.2.2) such that $z_i$ is a nonzero vector in $v_i$ and the $[z_i]_+$ are the vertices of $\hat{C}$. Remark that the $[z_i]_+$, for $i$ in $\{1,2,3\}$, are lifts of $v_1$, $v_2$ and $v_3$ in $\hat{\Lambda}$, respectively.

We write every point $p$ of $\mathbf{Ein}^{1,n}$ as

$$p=\mathbf{P}(p_1z_1+p_2z_2+p_3z_3+p_4z_4+p_5)\ ,$$

for $p_i$ in $\mathbb{R}$ and $p_5$ in $(v_1\oplus v_2\oplus v_3\oplus v_4)^\perp\simeq\mathbb{R}^{0,n-1}$. We use the notation

$$p=[p_1:p_2:p_3:p_4:p_5]\ .$$

**Lemma 3.4.2.** *Let* $p$ *be a point in* $\Lambda\setminus(\phi\cup\psi)$*, and*

$$[p_1z_1+p_2z_2+p_3z_3+p_4z_4+p_5]_+$$

*be the corresponding lift in* $\hat{\Lambda}$*. The coordinates* $p_1$*,* $p_2$*,* $p_3$ *and* $p_4$ *satisfy*

$$p_4>0,\qquad p_1\geq 0\quad\textit{and}\quad p_3\geq 0\ .$$

*Moreover, if* $p_1>0$ *and* $p_3>0$*, then* $p_2<0$*.*

**Proof.** The loop $\Lambda$ is a semi-positive loop. By Proposition 2.2.3(iii) we have $\langle p,z_i\rangle\leq 0$ for every $i$ in $\{1,2,3\}$. Since $\langle z_i,z_i+2\rangle=-1/4$, we obtain

$$p_1=-4\langle p,z_3\rangle\geq 0\quad\text{and}\quad p_3=-4\langle p,z_1\rangle\geq 0\ .$$

Since $[p]$ is not in $\phi\cup\psi$, is it not in a photon with $v_2$, hence

$$p_4=-4\langle p,z_2\rangle>0\ .$$

For the last assertion, we compute

$$\begin{aligned}0=\mathbf{q}(p)&=\sum_{i,j=1}^{4}\langle p_iz_i,p_jz_j\rangle+\mathbf{q}(p_5)\\&=2\langle p_1z_1,p_3z_3\rangle+2\langle p_2z_2,p_4z_4\rangle+\mathbf{q}(p_5)\\&=-\frac{1}{2}\left(p_1p_3+p_2p_4\right)+\mathbf{q}(p_5)\ .\end{aligned}$$

Since $p_5$ is a vector in a negative definite subspace of $(\mathbf{E}, \mathbf{q})$, we have $\mathbf{q}(p_5) \leq 0$. By the previous computation,

$$p_1p_3 + p_2p_4 \leq 0 \ .$$

Recall that from the first part of the lemma $p_4 > 0$, hence if $p_1 > 0$ and $p_3 > 0$, then

$$p_2 < 0 \ . \quad \square$$

Now we are ready to prove the renormalization lemma. Recall that the action of an element $a(\lambda, \mu)$ on the point $p$ is

$$a(\lambda, \mu) \cdot p = [\lambda p_1 : \mu p_2 : \lambda^{-1} p_3 : \mu^{-1} p_4 : p_5] \ .$$

We want to show that every subsequence of $(a_k \cdot \Lambda)_{k\in\mathbb{N}}$ admits a subsequence tending to $C$. We will show that a sublimit of $(a_k \cdot \Lambda)_{k\in\mathbb{N}}$ contains $v_1$, $v_2$, $v_3$ and $v_4$.

**Proof of Lemma 3.4.1.** Remark that, up to extracting a subsequence, $\lambda_k \geq 1$ or $\lambda_k \leq 1$. The two cases are similar so we only consider the case $\lambda_k \geq 1$. In our case, the sequence $(\lambda_k^{-1}\mu_k)_{k\in\mathbb{N}}$ tends to 0. There are several possibilities for the sequence $(\lambda_k\mu_k)_{k\in\mathbb{N}}$:

(1) the sequence $(\lambda_k\mu_k)_{k\in\mathbb{N}}$ tends to 0,
(2) the sequence $(\lambda_k\mu_k)_{k\in\mathbb{N}}$ is bounded from below by a positive constant.

**Case (1).** Suppose that $(\lambda_k\mu_k)_{k\in\mathbb{N}}$ tends to 0. Take any point $p$ in $\Lambda \setminus (\phi \cup \psi)$. By assumptions, since $p$ is not in a photon with $v_2$ and we can write

$$p = [p_1 : p_2 : p_3 : 1 : p_5] \ ,$$

for $p_1, p_2, p_3$ in $\mathbb{R}$ and $p_5$ in $(v_1 \oplus v_2 \oplus v_3 \oplus v_4)^{\perp}$. Acting with $a_k$ we obtain

$$a_k \cdot p = [\lambda_k p_1 : \mu_k p_2 : \lambda_k^{-1} p_3 : \mu_k^{-1} : p_5] = [\lambda_k\mu_k p_1 : \mu_k^2 p_2 : \frac{\mu_k}{\lambda_k} p_3 : 1 : \mu_k p_5] \ .$$

It follows that the sequence $(a_k \cdot p)_{k\in\mathbb{N}}$ converges to $v_4$ since $(\mu_k)_{k\in\mathbb{N}}$, $(\lambda_k\mu_k)_{k\in\mathbb{N}}$ and $(\lambda_k^{-1}\mu_k)_{k\in\mathbb{N}}$ tend to 0.

**Case (2).** Suppose now that $(\lambda_k\mu_k)_{k\in\mathbb{N}}$ is bounded from below by a positive constant. Then the sequence $(\lambda_k)_{k\in\mathbb{N}}$ tends to $+\infty$. Up to extracting a subsequence, the sequence $(\lambda_k\mu_k)_{k\in\mathbb{N}}$ tends to a positive real number $r$. Take now $p$ in $\Lambda \setminus (\phi \cup \psi)$.

- If $p$ is in an edge with $v_3$ then $p_1 = 0$. By the same computation as above we conclude that $(a_k \cdot p)_{k\in\mathbb{N}}$ tends to $v_4$.

- If $p$ is neither in an edge with $v_3$ neither in an edge with $v_1$, we write $p = [p_1 : p_2 : p_3 : p_4 : p_5]$. By Lemma 3.4.2, we can assume $p_4 > 0$, $p_1 > 0$, $p_3 > 0$ and hence $p_2 < 0$. Looking at the sequence $(a_k \cdot p)_{k\in\mathbb{N}}$, we have

$$\lim_{k\to+\infty} a_k \cdot p = \lim_{k\to+\infty} \left[\frac{\lambda_k \mu_k p_1}{p_4} : \frac{\mu_k^2 p_2}{p_4} : \frac{\mu_k p_3}{\lambda_k p_4} : 1 : \frac{\mu_k p_5}{p_4}\right] = \left[r\frac{p_1}{p_4} : 0 : 0 : 1 : 0\right] .$$

Now by taking the point $p$ close to $v_3$ in $\Lambda \setminus (\phi \cup \psi)$, we can suppose that $p_1$, $p_2$, $p_4$ and $\mathbf{q}(p_5)$ are as close to 0 as we want. Since

$$\mathbf{q}(p_5) = \frac{1}{2}(p_1 p_3 + p_2 p_4) ,$$

the number $p_1/p_4$ is also as close to 0 as we want. Hence, up to extracting a subsequence, the limit of the sequence $(a_k \cdot \Lambda)_{k\in\mathbb{N}}$ contains points as close to $v_4$ as we want, hence it must contain $v_4$.

In every case, up to extracting a subsequence, the sequence $(a_k \cdot \Lambda)_{k\in\mathbb{N}}$ converges to a semi-positive loop containing $v_1$, $v_2$, $v_3$ (because each isometry $a_k$ preserves these points) and $v_4$ (by what precedes). Since a semi-positive loop contains a photon segment between two points as soon as two points lie on such a segment (Proposition II.7(ii)), any sublimit of $(a_k \cdot \Lambda)_{k\in\mathbb{N}}$ is the Barbot crown $C$ with vertices $v_1$, $v_2$, $v_3$ and $v_4$. □

## 4. The space boundary and the ideal boundary of maximal surfaces

In this section, we study the link between the ideal boundary of a maximal surface $\Sigma$ and its space boundary. As we have seen in Section 2.3.3, the space boundary $\partial^{\mathrm{s}}\Sigma$ can be identified with a set of space-horofunctions, allowing us to compare it with the ideal boundary $\Sigma(\infty)$.

We start by a study of Barbot surfaces in detail in Section 4.1. In Section 4.2 we explain a nontrivial example of polygonal surface: maximal surfaces bounded by symmetric polygons. In Section 4.3 we define a projection between maximal surfaces and explain the renormalization process for maximal surfaces. In Section 4.4, we study our most important tool, the *space-horofunctions* and their properties. We then link the ideal boundary and the space boundary in Section 4.5. Section 4.6 is devoted to the Tits distance on Hadamard manifolds, and how to compute it. With this toolbox, we compute the Tits distance of the ideal boundary of asymptotically flat maximal surfaces in Section 4.7 and link the total curvature with the properties of the space boundary in Section 4.8.

### *4.1. The case of Barbot surfaces*

A *Barbot surface* is a maximal surface having a Barbot crown as space boundary. The Barbot surfaces are the easiest polygonal surfaces to study: their boundaries have the least number of vertices (Proposition 3.1.4). The Barbot surfaces have flat induced metric and a totally explicit description. Given the four vertices of their space boundary, we describe the geodesics of Barbot surfaces, and in particular their limit points in $\mathbf{Ein}^{1,n}$. This allows us to understand the relation between the space boundary and the ideal boundary of Barbot surfaces. At the end of this section we discuss space-horofunctions in Barbot surfaces.

*4.1.1. Explicit parametrization*

Let $v_1$, $v_2$, $v_3$ and $v_4$ be the vertices of a Barbot crown $C$ in $\mathbf{H}^{2,n}$. Let $S$ be the maximal surface with space boundary $C$. Recall that the Cartan subgroup $\mathsf{A}$ in $\mathsf{G}$ corresponding to $C$ is defined to be the subgroup stabilizing the four vertices $v_i$ globally and $(v_1 \oplus v_2 \oplus v_3 \oplus v_4)^\perp$ pointwise. Let $(z_1, z_2, z_3, z_4)$ be a hyperbolic basis with $z_i$ in $v_i$ and define $x = [z_1 + z_2 + z_3 + z_4]$.

**Proposition 4.1.1.** *The Barbot surface $S$ is the orbit of $x$ under the action of $\mathsf{A}_0$ and has flat induced metric.*

**Proof.** We claim that the orbit $\mathsf{A}_0 \cdot x$ has vanishing mean curvature at $x$. Since $\mathsf{A}_0 \cdot x$ is included in $v_1 \oplus v_2 \oplus v_3 \oplus v_4$, we consider only the restriction of $\mathsf{A}$ to $v_1 \oplus v_2 \oplus v_3 \oplus v_4$. We do the computations in the double cover $\widehat{\mathbf{H}}^{2,n}$, where $[z_i]_+$ are the vertices of a lift of $C$ and where the lift of $x$ is still denoted $x$. One parameter subgroups of $\mathsf{A}_0$ are given by

$$\gamma_X(t) = \mathrm{diag}(\exp(tu), \exp(tv), \exp(-tu), \exp(-tv)) \ ,$$

in the basis $(z_1, z_2, z_3, z_4)$ and with $u, v$ in $\mathbb{R}$. We identify $\mathfrak{a} = \mathsf{T}_x(\mathsf{A}_0 \cdot x)$ with $\mathbb{R}^2$ via the isometric map that sends the matrix of $\mathfrak{a}$ with diagonal $(u, v, -u, -v)$ to the vector $(u, v)$. Hence a geodesic ray $\sigma_X$ of $\mathsf{A}_0 \cdot x$ starting from $x$ with initial velocity $X = (u, v)$ is parametrized as follows:

$$\sigma_X(t) = \gamma_X(t) \cdot x = \exp(tu)z_1 + \exp(tv)z_2 + \exp(-tu)z_3 + \exp(-tv)z_4. \qquad (1)$$

Now taking orthogonal unit tangent vectors $X = (1, 0)$ and $Y = (0, 1)$, the mean curvature vector $H$ of $\mathsf{A}_0 \cdot x$ is given by

$$H(x) = \mathbb{I}(\dot\sigma_X(0), \dot\sigma_X(0)) + \mathbb{I}(\dot\sigma_Y(0), \dot\sigma_Y(0)) = (z_1 + z_3 + z_2 + z_4)^\perp = 0 \ .$$

In the last equality, we take the part of the vector normal to $\mathsf{T}_x(\mathsf{A}_0 \cdot x)$. The orbit $\mathsf{A}_0 \cdot x$ has vanishing mean curvature at $x$, hence everywhere by homogeneity. Thus it is a maximal surface. It is flat and complete, being the orbit of the isometric action of $\mathsf{A}_0 \simeq \mathbb{R}^2$. It remains to show that the space boundary of $\mathsf{A}_0 \cdot x$ is $\hat{C}$. One can verify this by looking at limits of intrinsic geodesics of $\mathsf{A}_0 \cdot x$. $\square$

We can parametrize, given $(z_1, z_2, z_3, z_4)$ a hyperbolic basis, a Barbot surface $S$ associated to this hyperbolic basis by the following isometry (onto its image):

$$\begin{aligned} f : \mathbb{R}^2 &\to \widehat{\mathbf{H}}^{2,n} \\ (u, v) &\mapsto \exp(u)z_1 + \exp(v)z_2 + \exp(-u)z_3 + \exp(-v)z_4 \ . \end{aligned}$$

We can now describe the tangent and normal bundle of the Barbot surface $S$. Denote by $\mathbf{H}^{2,1}$ the projectivization of the set of negative lines of $v_1 \oplus v_2 \oplus v_3 \oplus v_4$. Let us define

- $\partial_u := \exp(u)z_1 - \exp(-u)z_3 = f_* \frac{\partial}{\partial_u}$,
- $\partial_v := \exp(v)z_2 - \exp(-v)z_4 = f_* \frac{\partial}{\partial_v}$,
- $\mathbf{n} := \exp(u)z_1 - \exp(v)z_2 + \exp(-u)z_3 - \exp(-v)z_4$.

**Proposition 4.1.2** *(Description of the tangent and normal bundles of a Barbot surface). At a point $x = f(u, v)$, the tangent bundle of the Barbot surface $S$ is spanned by the vectors $\partial_u$ and $\partial_v$. The normal bundle of $S$ splits in two parts, one in the tangent space of $\mathbf{H}^{2,1}$ spanned by the vector $\mathbf{n}$ and the other one in $\mathsf{T}\mathbf{H}^{2,1\perp}$. The second fundamental form $\mathbb{I}$ of $S$ satisfies, at the point $x = f(u, v)$,*

$$\mathbb{I}(\partial_u, \partial_u) = \nabla_{\partial_u}\partial_u = \frac{1}{2}\mathbf{n} \ ,$$

*where $\nabla$ is the Levi-Civita connection of $\mathbf{H}^{2,n}$. The shape operator $B_\mathbf{n}$ of $S$ satisfies, at the point $x = f(u, v)$ and relatively to the normal direction $\mathbf{n}$,*

$$B_\mathbf{n}(\partial_u) = -\partial_u, \quad B_\mathbf{n}(\partial_v) = \partial_v \ .$$

**Proof.** For the second fundamental form:

$$\begin{cases} \nabla_{\partial_u}\partial_u = \frac{1}{2}(\exp(u)z_1 - \exp(v)z_2 + \exp(-u)z_3 - \exp(-v)z_4) \ , \\ \nabla_{\partial_v}\partial_u = \nabla_{\partial_u}\partial_v = 0 \ , \\ \nabla_{\partial_v}\partial_v = -\nabla_{\partial_u}\partial_u \ . \end{cases}$$

For the shape operator:

$$\begin{cases} \langle B_\mathbf{n}(\partial_u), \partial_u\rangle = \langle \mathbb{I}(\partial_u, \partial_u), \mathbf{n}\rangle = -1/2 \ , \\ \langle B_\mathbf{n}(\partial_u), \partial_v\rangle = \langle B_\mathbf{n}(\partial_v), \partial_u\rangle = \langle \mathbb{I}(\partial_u, \partial_v), \mathbf{n}\rangle = 0 \ , \\ \langle B_\mathbf{n}(\partial_v), \partial_v\rangle = \langle \mathbb{I}(\partial_v, \partial_v), \mathbf{n}\rangle = 1/2 \ . \end{cases} \qquad \square$$

### *4.1.2. Geodesics*

Barbot surfaces are isometric to the Euclidean plane. Their intrinsic geometry is simple, and we investigate now their extrinsic geometry by the asymptotic behavior of geodesics in the space boundary. That is, we study the relation between the space boundary of a Barbot surface and its ideal boundary. Recall that $S$ is a Barbot surface with space boundary the Barbot crown $C = \partial^{\mathrm{s}} S$ and with ideal boundary $S(\infty)$.

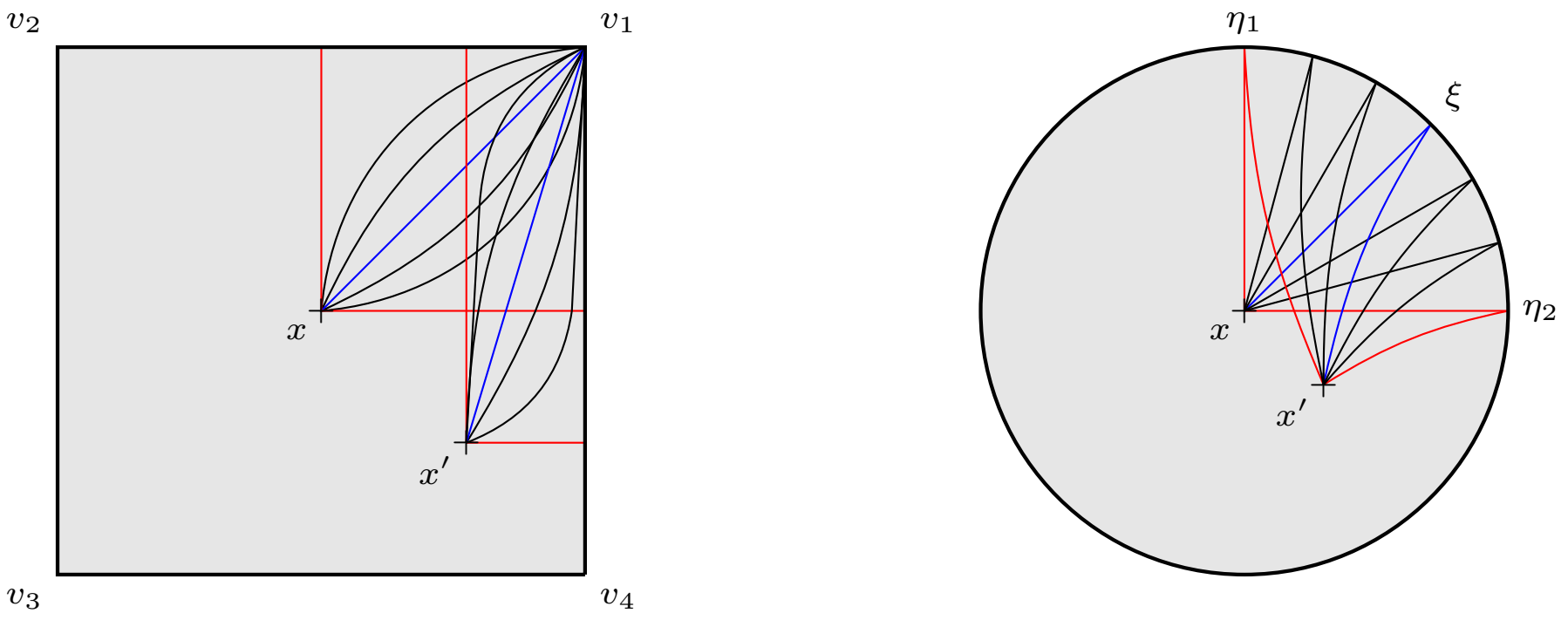

**Fig. 3.** A Barbot surface compactified with its space boundary (left) and its ideal boundary (right), respectively. In red, unstable geodesic rays corresponding to the edges $[v_1, v_2]$ and $[v_1, v_4]$ and to the ideal boundary points $\eta_1$ and $\eta_2$, respectively. In blue, stable geodesic rays corresponding to the vertex $v_1$ and to the ideal boundary point $\xi$.

**Proposition 4.1.3** *(Behavior of geodesics at infinity of a Barbot surface). There are four ideal points in $S(\infty)$ such that any representing ray of these ideal points converges, in $S \cup C$, to a point of the interior of an edge. Every other ray converges to a vertex (see Fig. 3).*

**Proof.** We first do the proof at the point $x = z_1 + z_2 + z_3 + z_4$. Let $X = (u, v)$. The geodesic ray of $S$ with initial velocity $X$ at $x$ is given by Equation (1). Suppose that $u = \pm v$. Then $\lim_{t \to +\infty} \exp(-|u|t)\sigma(t)$ equals

- $z_1 + z_2$ if $u = v > 0$,
- $z_2 + z_3$ if $u = -v$ and $u < 0$,
- $z_3 + z_4$ if $u = v < 0$,
- $z_4 + z_1$ if $u = -v$ and $u > 0$.

Suppose now that $u \neq \pm v$. Then $\lim_{t \to +\infty} \exp(-\max(|u|, |v|)t)\sigma(t)$ equals

- $z_1$ if $u > |v|$,
- $z_2$ if $v > |u|$,
- $z_3$ if $-u > |v|$,
- $z_4$ if $-v > |u|$.

Hence the four directions, at $x$, going to edges, are directed by the vectors $\pm(1, 1)$ and $\pm(-1, 1)$.

If $y$ is another point of $S$, then $y = a \cdot x$ for an element $a$ of $\mathsf{A}_0$. The differential of the map $a$, at $x$, sends the vector $(1, 1)$ to $(\partial_u + \partial_v)_y$ and sends the vector $(1, -1)$ to $(\partial_u - \partial_v)_y$. Since $a$ acts by isometries on $S$, stabilizes the vertices, sends points of an edge to points on the same edge transitively, this finishes the proof. □

*4.1.3. Stable and unstable directions*

We have computed in Proposition 4.1.2 the eigenspaces of the shape operator of a Barbot surface $S$. The vector fields $\partial_u$ and $\partial_v$ are eigenvectors of the shape operator (recall Section 4.1.1 for the definition of these vector fields), are unitary, and their integral curves are geodesics in $S$, hence they represent ideal boundary points in $S(\infty)$.

**Definition 4.1.4.** In a Barbot surface $S$, the ideal boundary points represented by the four eigenvectors $\partial_u$, $\partial_v$, $-\partial_u$ and $-\partial_v$, respectively, are called *stable* directions. The geodesic rays with velocity a stable direction are called *stable rays.*

**Remark 4.1.5.** A geodesic ray representing a stable direction will converge to a vertex of $C = \partial^{\mathrm{s}} S$. We use the terminology stable because given a stable direction, there is a neighborhood $\mathcal{U}$, in $S(\infty)$, of this direction such that for every direction $\xi$ in $\mathcal{U}$, every geodesic ray $c$ representing $\xi$ converges to the same vertex in $C$. See Fig. 3.

**Proposition 4.1.6.** *The angle distance between two consecutive stable ideal directions in a Barbot surface is* $\pi/2$.

**Remark 4.1.7.** The four vertices of the space boundary of a Barbot crown are in bijection with the four stable directions of the associated Barbot surface.

**Definition 4.1.8.** In a Barbot surface $S$, the four directions in $S(\infty)$ going to edges are called *unstable directions.* The geodesic rays with velocity an unstable direction are called *unstable rays.*

**Remark 4.1.9.** Let $\xi$ be an unstable direction of a Barbot surface $S$. Then the limit point, in $\partial^{\mathrm{s}} S$, of a geodesic ray of $S$ starting at $x$ with initial direction $\xi_x$ depends on $x$. The collection of these limit points, when $x$ moves through $S$, equals the interior of an edge. See Fig. 3.

*4.1.4. Space-horofunctions in Barbot surfaces*

Let $S$ be a Barbot surface in $\widehat{\mathbf{H}}^{2,n}$ with space boundary $C$ in $\widehat{\mathbf{Ein}}^{1,n}$ with vertices $v_1$, $v_2$, $v_3$ and $v_4$, and $(z_1, z_2, z_3, z_4)$ be a hyperbolic basis corresponding to these vertices. Recall that this data determines a parametrization of $S$ via the map

$$\begin{aligned} f : \mathbb{R}^2 &\to \widehat{\mathbf{H}}^{2,n} \\ (u, v) &\mapsto \exp(u) z_1 + \exp(v) z_2 + \exp(-u) z_3 + \exp(-v) z_4 \ . \end{aligned}$$

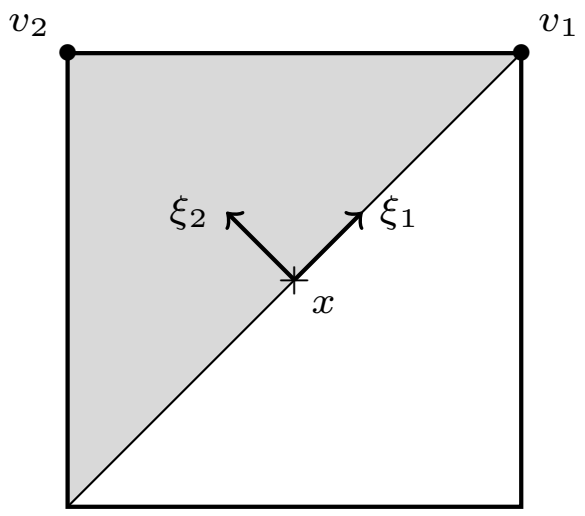

**Fig. 4.** A space-horoball (in gray) in a Barbot surface, associated to the vertex $v_2$. The stable directions $\xi_2$ and $\xi_1$ associated to $v_2$ and $v_1$, respectively, are represented at the point $x$. The direction $\xi_1$ is tangent to the space-horocycle associated to $v_2$ and passing through $x$.

See Section 4.1.1 for details. A direct computation using the formula defining space-horofunctions, see Definition 2.3.6, gives the following.

**Proposition 4.1.10.** *Let $S$ be a Barbot surface of $\widehat{\mathbf{H}}^{2,n}$ with space boundary $C$ and vertices $v_1$, $v_2$, $v_3$ and $v_4$. The space-horofunction associated to a vertex, say $v_1$, has a representative $h$ defined, under the parametrization $f$, by*

$$h \circ f((u,v)) = -u \ ,$$

*and the space-horofunction associated to a point on an edge, say $\mathbf{P}(tz_1 + (1-t)z_2)$, has a representative $\tilde{h}$ defined, under the parametrization $f$, by*

$$\tilde{h} \circ f((u,v)) = \log(t \exp(-u) + (1-t)\exp(-v)) \ .$$

**Remark 4.1.11.** In the case of vertices, in a Barbot surface $S$, the space-horofunction equals the intrinsic horofunction associated to the stable direction tending to the vertex. The horoballs associated to these horofunctions are half planes, and their boundaries are stable geodesics joining the two neighboring vertices. See Fig. 4. Hence a stable direction in $S(\infty)$, corresponding to the vertex $v_i$, can be determined by the space horofunctions of the neighboring vertices $v_{i-1}$ and $v_{i+1}$.

We will use this idea to define geodesic rays associated to a vertex in a maximal surface.

## *4.2. The case of symmetric polygonal surfaces in* $\mathbf{H}^{2,1}$

We present now a collection of examples of polygonal surfaces where we are able to describe the intrinsic and extrinsic geometry. These examples sit inside $\mathbf{H}^{2,1}$, when the polygon is symmetric under the action of a dihedral subgroup of the isometry group of $\mathbf{H}^{2,1}$. For the definition of the symmetric polygons, see Section 3.3. In this section we assume the results of the paper, in particular the main theorem, in order to show a nontrivial example of polygonal maximal surface to the reader before going to more abstract considerations.

Let $N$ be an even integer bigger than 4, and let $P$ be a symmetric polygon with $N$ vertices. Denote by $\Sigma$ the maximal surface with space boundary the polygon $P$. Recall that the polygon $P$ is constructed using a splitting of $\widehat{\mathbf{Ein}}^{1,1}$ as being $\mathbf{S}^1 \times \mathbf{S}^1$.

*4.2.1. Quartic differential*

We claim that the quartic differential $\varphi_4$ of $\Sigma$ is a monomial. Recall that, by Theorem 1.2.1, $\varphi_4$ is a polynomial of degree $N-4$ since the polygon $P$ has $N$ vertices.

**Proof.** Let $x$ be a point of $\Sigma$ being a zero of $\varphi_4$. The orbit of $x$ by the action of the subgroup generated by $\mathbf{r}$ consists only of zeros of $\varphi_4$. If $x$ is not fixed by $\mathbf{r}$, then we have $N$ pairwise distinct zeros of $\varphi_4$, which is impossible (the degree of $\varphi_4$ is $N-4$). As a consequence, $x$ is fixed by $\mathbf{r}$. But there is a unique point of $\Sigma$ fixed by $\mathbf{r}$, namely the point in the timelike sphere $\{0\} \times \mathbf{S}^1$ in the warped decomposition of $\widehat{\mathbf{H}}^{2,1}$ as being $\mathbf{D}^2 \times \mathbf{S}^1$. Finally, all the zeros of $\varphi_4$ are at the same point, that is $\varphi$ is a monomial. □

We chose a holomorphic chart $z$ in which the quartic differential is written

$$\varphi_4 = z^{N-4}\mathrm{d}z^4 \ ,$$

hence the point $z = 0$ is invariant by $\mathbf{r}$ and $\mathbf{s}$.

*4.2.2. Induced metric*

The maximal surface $\Sigma$ is defined using the unique solution to the following equation (see Equation (1) of [29]), where $u$ is a smooth map defined on $\mathbb{C}$:

$$u_{z\bar{z}} - \frac{1}{2}\exp(2u) + \frac{1}{2}\exp(-2u)|z^{N-4}|^2 = 0 \ . \tag{2}$$

Equation (2) being rotationally symmetric around 0, that is invariant under $z \mapsto \exp(i\theta)z$, the solution $u$ is rotationally symmetric. The induced metric of $\Sigma$ is recovered from $u$ (see [29, Paragraph 3.2]) as being

$$g = 2\exp(2u)|\mathrm{d}z|^2 \ ,$$

hence $g$ is rotationally symmetric, that is $g_z = g_{\exp(i\theta)z}$ for every real number $\theta$.

**Remark 4.2.1.** It should be noted that, although the second fundamental form has a finite group of symmetries, the induced metric has an entire circle of symmetries. In an intrinsic point of view, every direction looks the same at $z = 0$, but the extrinsic behavior of geodesics, as we shall see, is encoded by the second fundamental form.

*4.2.3. Ideal boundary*

Using the extrinsic symmetries of $P$, hence of the maximal surface $\Sigma$, we can identify several geodesics on $\Sigma$ passing through the point $z = 0$.

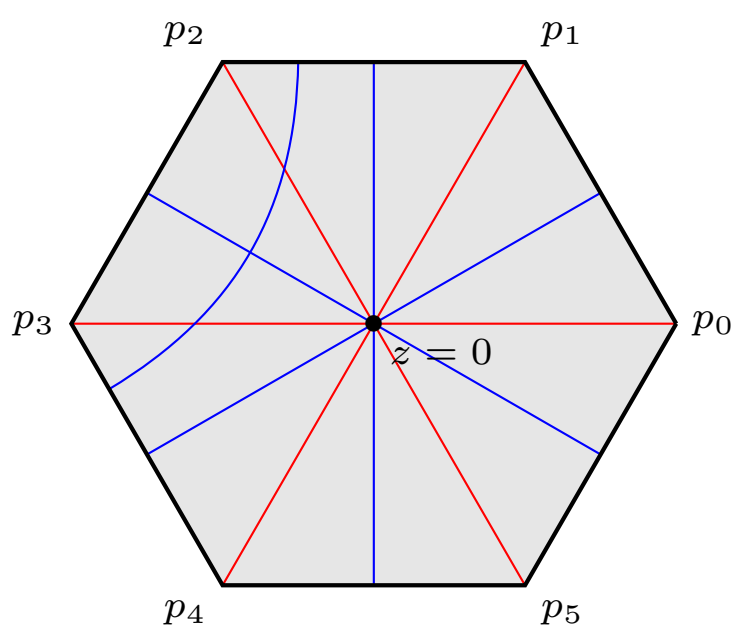

**Fig. 5.** A symmetric polygonal surface seen from above, with symmetry center $z = 0$. In red, the stable directions passing through the symmetry center. In blue, some vertical trajectories.

For every $k$ in $\{1, \ldots, N\}$, the symmetry $\mathbf{r}^k\mathbf{s}$ fixes pointwise a unique curve in $\Sigma$ passing through $z = 0$. This curve is hence a geodesic, denoted $c_k$, of $\Sigma$. In the cases where $k$ is even, this geodesic tends, at infinity, to the vertices $p_{k/2}$ and $p_{(N+k)/2}$ of $P$. In the cases where $k$ is odd, the geodesic $c_k$ has the following limit points in the space boundary:

$$\lim_{t \to \pm\infty} c_k(t) = \{(\exp(i\pi k/N), 1), (-\exp(i\pi k/N), 1)\} ,$$

which corresponds to two points of $P$, one in the interior of the edge joining $p_{\lfloor k/2 \rfloor}$ and $p_{\lceil k/2 \rceil}$, one in the interior of the edge joining $p_{\lfloor (N+k)/2 \rfloor}$ and $p_{\lceil (N+k)/2 \rceil}$, see Fig. 5.

### *4.2.4. Stable and unstable directions*

The previous section exhibits $2N$ points of the ideal boundary, namely the directions (forward and backward) of the geodesics $c_k$ at 0. These directions have an interpretation in terms of the differential $\varphi_4$: they are the "horizontal" and "vertical" directions of the differential.

**Definition 4.2.2.** A *horizontal direction* of a quartic differential $\varphi$, at a point $z$ of $\mathbb{C}$, is a tangent vector $u$ in $\mathsf{T}_z\mathbb{C}$ such that $\varphi((z, u)) > 0$. Similarly, a *vertical direction* of a quartic differential $\varphi$, at $z$, is a tangent vector $u$ in $\mathsf{T}_z\mathbb{C}$ such that $\varphi((z, u)) < 0$.

These definitions mimic the definition of horizontal and vertical directions of a quadratic differential. However, there are 4 horizontal directions (up to $\mathbb{R}_+^*$) at each regular point, and 4 vertical directions.

By symmetry and renormalization, we see that the horizontal directions passing through the symmetry center are stable directions of the maximal surface. The vertical directions can be thought of as generalization of unstable directions of the Barbot surface. But the horizontal and vertical trajectories that does not pass through the symmetry center are not geodesics, and hence are not exactly the stable and unstable directions.

The behavior of every horizontal or vertical trajectory is easy to guess, because they stay at finite distance from the ones passing through the symmetry center. An example is given in Fig. 5 for a vertical trajectory that does not pass through the symmetry center.

### 4.3. Renormalization

The aim of this Section is to show that polygonal surfaces are asymptotically flat. This is done using

- the *radial projection* between maximal surfaces (next section),
- the renormalization process (Section 4.3.4).

#### 4.3.1. The radial projection

Here we define the *radial projection* between two maximal surfaces. Take $\Sigma_1$ and $\Sigma_2$ two maximal surfaces of $\mathbf{H}^{2,n}$. We will project $\Sigma_1$ to $\Sigma_2$ following a direction in $\mathsf{N}_x\Sigma_1$ at each point $x$ of $\Sigma_1$.

**Definition 4.3.1.** The *radial projection* from $\Sigma_1$ to $\Sigma_2$ is the map that sends a point $x$ in $\Sigma_1$ to the unique point of $\Sigma_2$ in the timelike sphere orthogonal to $\mathsf{T}_x\Sigma_1$.

Recall that the timelike sphere passing through $x$ is $\exp_x(\mathsf{N}_x\Sigma)$ (see Section 2.1.2). Because complete spacelike surfaces are graphs, the radial projection is well defined (also see [19, Proposition 3.10 (iii)]).

We will always denote by $\rho$ the radial projection between two maximal surfaces. The timelike spheres orthogonal to $\mathsf{T}_x\Sigma_1$ vary smoothly with the point $\sigma$, so the projection $\rho$ is smooth.

#### 4.3.2. Local injectivity of the radial projection

In general, the radial projection between two maximal surfaces is not injective.

**Example 4.3.2.** Let $\mathbf{H}^2$ denote a totally geodesic copy of the hyperbolic plane in $\mathbf{H}^{2,1}$. This $\mathbf{H}^2$ corresponds to a subspace $F$ of signature $(2, 1)$ of $\mathbb{R}^{2,2}$. The orthogonal of $F$ is a negative definite line, hence a point of $\mathbf{H}^{2,1}$ denoted $x$. Then every timelike geodesic orthogonal to $\mathbf{H}^2$ pass through $x$. As a consequence, the radial projection from $\mathbf{H}^2$ to any maximal surface passing through $x$ is constant.

This example deals with two maximal surfaces far from each other. The situation is much better when two maximal surfaces are close. The closeness between two maximal surfaces is related to the following function.

**Definition 4.3.3** *(Timelike distance).* Let $x$ and $y$ be two different points of $\mathbf{H}^{2,n}$ in a timelike geodesic. Take $\hat{x}$ and $\hat{y}$ lifts of $x$ and $y$, respectively, to $\widehat{\mathbf{H}}^{2,n}$, such that $\langle \hat{x}, \hat{y} \rangle \leq 0$. The *timelike distance* between $x$ and $y$ is defined to be the number

$$\delta^{\mathrm{t}}(x,y) := \arccos(-\langle \hat{x}, \hat{y}\rangle) \ .$$

**Proposition 4.3.4.** *Let $\Sigma_1$ and $\Sigma_2$ be two maximal surfaces of $\mathbf{H}^{2,n}$ and $\rho$ denote the radial projection from $\Sigma_1$ to $\Sigma_2$. For every $x$ in $\Sigma_1$, if the timelike distance $\delta^{\mathrm{t}}(x,\rho(x))$ is smaller than $\arcsin(1/\sqrt{5})$, then there is a neighborhood $\mathcal{U}$ of $x$ in $\Sigma_1$ such that the restriction $\rho_{|\mathcal{U}}$ is a diffeomorphism onto its image.*

The proof is a careful analysis on the properties of some Jacobi fields in $\mathbf{H}^{2,n}$. We recall now some definitions and computations that can be found in [23, Chapter 8], before giving the proof. Let $\gamma$ be a geodesic of $\mathbf{H}^{2,n}$. A *Jacobi field along* $\gamma$ is a vector field $Y$ along $\gamma$ such that, denoting by $\dot{Y}$ and $\ddot{Y}$ the first and second covariant derivative of $Y$ along $\gamma$, we have the *Jacobi equation*:

$$\ddot{Y} + \mathsf{R}_{Y,\dot{\gamma}}\dot{\gamma} = 0 \ ,$$

where $\mathsf{R}$ denotes the Riemann curvature tensor of $\mathbf{H}^{2,n}$. The Jacobi fields along $\gamma$ are exactly the infinitesimal variation of one-parameter families of geodesics starting at $\gamma$. The Jacobi equation is an order 2 linear equation, hence a Jacobi field is determined by its value and derivative at a point. A Jacobi field is said to be *normal* if it is orthogonal to $\dot{\gamma}$.

If $\gamma$ is a timelike geodesic parametrized by arclength, and $Y$ a normal Jacobi field along $\gamma$, then the Jacobi equation reduces to

$$\ddot{Y} = -Y \ ,$$

because $\mathbf{H}^{2,n}$ has constant curvature $-1$. Moreover, if $U$ denotes the parallel transport of $Y(0)$ along $\gamma$ and $V$ denotes the parallel transport of $\dot{Y}(0)$ along $\gamma$, we have

$$Y(t) = \cos(t)U + \sin(t)V \ .$$

In a pseudo-Riemannian symmetric space, the parallel transport along a geodesic $\gamma$ corresponds to the differential of the translation along $\gamma$. Take a lift $\hat{\gamma}$ of $\gamma$ to $\widehat{\mathbf{H}}^{2,n}$. Let $X$ be a vector in $\mathsf{T}_{\hat{\gamma}(0)}\widehat{\mathbf{H}}^{2,n}$, orthogonal to the velocity of $\hat{\gamma}$ at 0. Then the parallel transport of $X$ along $\hat{\gamma}$ is constant equal to $X$. At the end, we obtain:

**Lemma 4.3.5.** *Let $\hat{\gamma}$ be a timelike geodesic of $\widehat{\mathbf{H}}^{2,n}$ parametrized by arclength and $Y$ be a normal Jacobi field along $\hat{\gamma}$. Then for every real number $t$, we have*

$$Y(t) = \cos(t)Y(0) + \sin(t)\dot{Y}(0) \ .$$

**Remark 4.3.6.** Let $x$ be a point in $\widehat{\mathbf{H}}^{2,n}$. The set $\{y \in \widehat{\mathbf{H}}^{2,n} \mid \langle x, y\rangle \leq 0\}$ is a fundamental domain for the covering $\widehat{\mathbf{H}}^{2,n} \to \mathbf{H}^{2,n}$. If $y$ is in a timelike geodesic starting from $x$, at timelike distance less than $\pi/2$, then $y$ is contained in the interior of this fundamental domain.

**Proof of Proposition 4.3.4.** Let $\Sigma_1$ and $\Sigma_2$ be two maximal surfaces of $\mathbf{H}^{2,n}$, $\rho$ be the radial projection from $\Sigma_1$ to $\Sigma_2$ and $x$ be a point of $\Sigma_1$ such that $\delta^{\mathrm{t}}(x, \rho(x)) < \arcsin(1/\sqrt{5})$. We want to show that $\mathrm{d}_x\rho$ is a linear isomorphism. Let $c$ denote a geodesic of $\Sigma_1$ such that $c(0) = x$ and with unit initial speed $c'(0)$. We claim that $(\rho \circ c)'(0)$ is a nonzero vector.

By continuity of the timelike distance function and of $\rho$, there is a neighborhood of $x$ in $\Sigma_1$ such that $\delta^{\mathrm{t}}(y, \rho(y)) < \arcsin(1/\sqrt{5})$ for every $y$ in this neighborhood. By the above remark, the curves $c$ and $\rho \circ c$ can be lifted consistently to $\widehat{\mathbf{H}}^{2,n}$, at least in a neighborhood of 0. We still denote such lifts $c$ and $\rho \circ c$.

We have now two spacelike curves $c$ and $\rho\circ c$ in $\widehat{\mathbf{H}}^{2,n}$, both defined in an interval $(-\varepsilon, \varepsilon)$, for a positive number $\varepsilon$. These curves satisfy the following property (up to shrink $\varepsilon$): for every $s$ in $(-\varepsilon, \varepsilon)$, the points $c(s)$ and $\rho \circ c(s)$ are related by a timelike geodesic of $\widehat{\mathbf{H}}^{2,n}$ of length smaller than $\arcsin(1/\sqrt{5})$. Denote by $\gamma_s$ the geodesic of $\widehat{\mathbf{H}}^{2,n}$ starting at $c(s)$ and passing at $\rho \circ c(s)$ at time $\delta^{\mathrm{t}}(c(s), \rho \circ c(s))$. The following map

$$\begin{aligned} V : (-\varepsilon, \varepsilon) \times \mathbb{R} &\to \widehat{\mathbf{H}}^{2,n} \\ (s,t) &\mapsto \gamma_s(t) \end{aligned}$$

is a one-parameter family of geodesics. Hence the vector field $J = \frac{\partial}{\partial s} \cdot V$ along $\gamma_0$ is a Jacobi field. We want to show that the vector $J$ at time $\delta^{\mathrm{s}}(x, \rho(x))$ is a nonzero vector. Denote by $Y$ the normal part of $J$. Hence $Y$ is a normal Jacobi field along $\gamma_0$ and by what precedes the proof we know that

$$Y(t) = \cos(t)Y(0) + \sin(t)\dot{Y}(0) \ .$$

Suppose that $Y$ vanishes at a certain positive time $t_0$. Then the normal Jacobi field $Y(\cdot + t_0)$ vanishes at 0, so it takes the form $Y(t + t_0) = \sin(t)\dot{Y}(t_0)$, hence

$$Y(t) = \sin(t - t_0)\dot{Y}(t_0) \ , \tag{3}$$

for every real number $t$.

**Claim:** the number $t_0$ is bigger than or equal to $\arcsin(1/\sqrt{5})$.

From Equation (3) we obtain the following equations:

$$Y(0) = -\sin(t_0)\dot{Y}(t_0) \ , \tag{4}$$

$$\dot{Y}(0) = \cos(t_0)\dot{Y}(t_0) \ . \tag{5}$$

Equation (4) implies that

$$1 = \mathbf{q}(Y(0)) = \sin(t_0)^2 \left(\dot{Y}(t_0)\right) ,$$

and Equations (4) and (5) imply that

$$\left(\dot{Y}(t_0)\right) = \mathbf{q}(Y(0)) + \mathbf{q}\left(\dot{Y}(0)\right) = 1 + \mathbf{q}\left(\dot{Y}(0)\right) .$$

Hence we have

$$\sin(t_0)^2 = \frac{1}{1 + \mathbf{q}\left(\dot{Y}(0)\right)} . \tag{6}$$

We want now to estimate the number $\mathbf{q}(\dot{Y}(0))$.

Since the Levi-Civita connection $\nabla$ of $\widehat{\mathbf{H}}^{2,n}$ is torsion free, and denoting by $S$ and $T$ the vector fields along $V$ defined by $S = \mathrm{d}V \cdot \frac{\partial}{\partial s}$ and $T = \mathrm{d}V \cdot \frac{\partial}{\partial t}$, respectively, we have:

$$\dot{J}(t) = (\nabla_T S)_{V(0,t)} = (\nabla_S T)_{V(0,t)} .$$

As a consequence, since $Y$ is the normal part of $J$, we have

$$\begin{aligned}
\langle \dot{Y}(0), Y(0) \rangle &= \langle \dot{Y}(0), c'(0) \rangle \\
&= \langle \dot{J}(0), c'(0) \rangle \\
&= \langle \nabla_S T, S \rangle_{V(0,0)} \\
&= \langle T, \nabla_S S \rangle_{V(0,0)} \\
&= \langle \dot{\gamma}_0(0), \mathbb{I}(c'(0), c'(0)) \rangle ,
\end{aligned}$$

where $\mathbb{I}$ denotes the second fundamental form of $\Sigma_1$. For a unitary vector $u$ at a point $x$ of $\Sigma_1$, we know by [19, Corollary 5.2] that $|\mathbf{q}(\mathbb{I}(u,u))| \leq 2$. Finally, by Cauchy-Schwarz inequality ($\dot{\gamma}_0(0)$ and $\mathbb{I}(c'(0), c'(0))$ live in a definite subspace: $\mathsf{N}_x\Sigma$),

$$|\langle \dot{Y}(0), Y(0) \rangle| \leq |\mathbf{q}(\dot{\gamma}_0(0))\mathbf{q}(\mathbb{I}(c'(0), c'(0)))| \leq 2 .$$

But $\dot{Y}(0) = -\frac{\cos(t_0)}{\sin(t_0)} Y(0)$, hence

$$\mathbf{q}(\dot{Y}(0)) = \left(\frac{\cos(t_0)}{\sin(t_0)}\right)^2 \leq 4 .$$

This last equation, together with Equation (6), gives

$$\sin(t_0) \geq \frac{1}{\sqrt{5}} ,$$

hence $\rho \circ c'(0) \neq 0$. We conclude with the inverse function theorem. □

#### *4.3.3. Asymptotics of a maximal surface near a vertex*

Let $\Sigma$ be a maximal surface with space boundary the semi-positive loop $\Lambda$. Suppose that $\Lambda$ has a vertex $v$ (recall Definition 3.1.2). By Corollary 3.2.8, there is a point $v_4$ in $\mathbf{Ein}^{1,n}$ such that $v_1$, $v_2$, $v_3$ and $v_4$ are the vertices of a Barbot crown $C$ coinciding with $\Lambda$ between $v_1$, $v_2$ and $v_3$, where $v_2 = v$ and $v_1$, $v_3$ are the points at the end of the two edges passing through $v$. We fix some notations (see Fig. 6):

- $S$ is the Barbot surface with space boundary $C$,
- $\gamma$ is the connected component of $\Lambda \setminus \{v_1, v_3\}$ containing $v_2$,
- $\bar{\Sigma}$ and $\bar{S}$ are the $\mathcal{C}^0$ surfaces with boundary $\Sigma \cup \gamma$ and $S \cup \gamma$, respectively,
- $\rho : S \to \Sigma$ denotes the radial projection (see Definition 4.3.1),
- $\mathsf{A}$ is the Cartan subgroup of $\mathsf{G}$ determined by the crown $C$ (see Section 3.2.5).

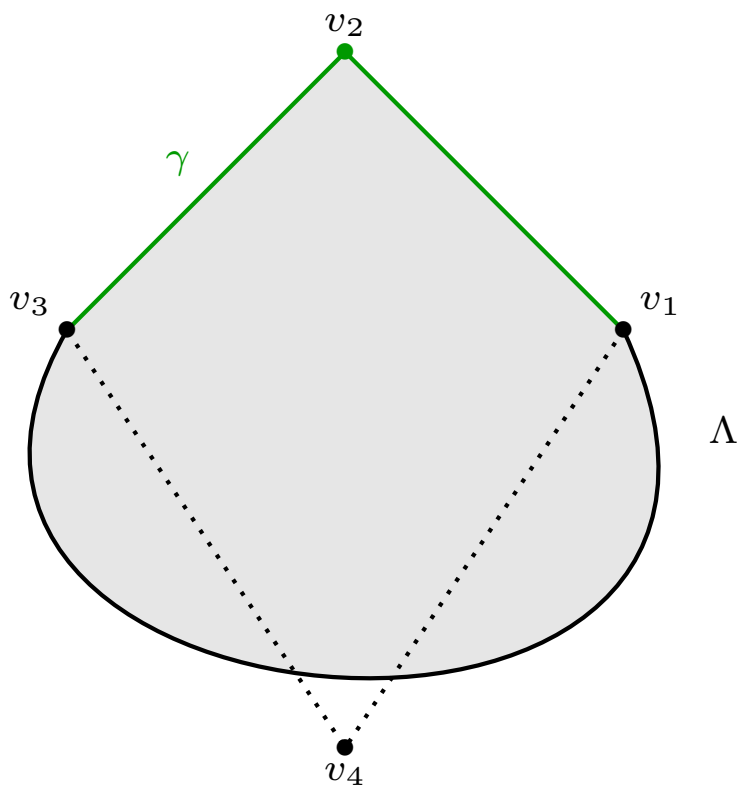

**Fig. 6.** A semi-positive loop $\Lambda$ with a Barbot crown coinciding with $\Lambda$ along an interval $\gamma$ containing the vertex $v_2$.

We now prove the following

**Proposition 4.3.7.** *The radial projection $\rho : S \to \Sigma$ extends as a continuous map $\bar{\rho} : \bar{S} \to \bar{\Sigma}$ that restricts to the identity on $\gamma$.*

**Proof.** Let $(x_k)_{k\in\mathbb{N}}$ be a sequence of points of $S$ tending to a point $p$ in $\gamma$. We fix a hyperbolic basis $(z_1, z_2, z_3, z_4)$ associated to the vertices $v_i$, and $x = [z_1 + z_2 + z_3 + z_4]$ in $S$. We have, for each integer $k$, real numbers $l_k, m_k$ such that

$$x_k = a(l_k, m_k) \cdot x = [l_k z_1 + m_k z_2 + l_k^{-1} z_3 + m_k^{-1} z_4] \ .$$

The sequence $(x_k)_{k\in\mathbb{N}}$ tends to a point $p$ in $\gamma$. That is:

$$[l_k : m_k : l_k^{-1} : m_k^{-1} : 0] = [l_k m_k^{-1} : 1 : l_k^{-1} m_k^{-1} : m_k^{-2} : 0] \to [a : 1 : b : 0 : 0] \ ,$$

for $a$, $b$ nonnegative real numbers. Hence we have

(1) the sequence $(m_k)_{k\in\mathbb{N}}$ tends to $+\infty$,
(2) the sequence $(l_k m_k)_{k\in\mathbb{N}}$ is bounded from below,

(3) the sequence $(l_k^{-1} m_k)_{k\in\mathbb{N}}$ is bounded from below.

Now we look at the sequence of points $(y_k)_{k\in\mathbb{N}}$ in $\Sigma$ defined by

$$y_k := \rho(x_k) \quad \forall k \in \mathbb{N} \ .$$

To conclude the proof of the proposition, it is enough to prove that the sequence $(y_k)_{k\in\mathbb{N}}$ converges to $p$.

**First step.** Claim: the timelike distance between $x_k$ and $y_k$ tends to 0 as $k$ goes to $+\infty$.

We look at the sequence $a_k = a(l_k^{-1}, m_k^{-1})$ of elements of $\mathsf{A}$. By definition, for all integers $k$, $a_k \cdot x_k = x$ and by the very definition of the radial projection, $a_k \cdot y_k$ is a point in the timelike sphere $\exp_x(\mathsf{N}_x S)$. By the renormalization lemma near a vertex (Lemma 3.4.1), up to extracting a subsequence, $(a_k \cdot \Lambda)_{k\in\mathbb{N}}$ converges to $C$ as $k$ goes to $+\infty$. Hence, by Theorem 2.2.8, the sequence of maximal surfaces $(a_k \cdot \Sigma)_{k\in\mathbb{N}}$ subconverges to $S$. Moreover, since $a_k \cdot y_k$ is in the timelike sphere defined by $x$ and $S$ for each integer $k$, we have

$$\lim_{k\to+\infty} a_k \cdot y_k = x \ .$$

Since $a_k$ is an isometry, the timelike distance between $x_k$ and $y_k$ equals the timelike distance between $x = a_k \cdot x_k$ and $a_k \cdot y_k$, hence tends to 0.

This concludes the first step.

To complete the proof of the proposition, we take coherent lifts $\hat{S}$, $\hat{\Sigma}$, $(\hat{x}_k)_{k\in\mathbb{N}}$ and $(\hat{y}_k)_{k\in\mathbb{N}}$ of $S$, $\Sigma$, $(x_k)_{k\in\mathbb{N}}$ and $(y_k)_{k\in\mathbb{N}}$ in $\widehat{\mathbf{H}}^{2,n}$, respectively.

**Second step.** Claim: the sequence $(\hat{y}_k)_{k\in\mathbb{N}}$, up to the multiplication by positive constants, tends to a vector in $p$.

For each integer $k$, we write $\hat{x}_k = \exp(l_k)z_1 + \exp(m_k)z_2 + \exp(-l_k)z_3 + \exp(-l_k)z_4$. The geodesic segment between $\hat{x}_k$ and $\hat{y}_k$ is parametrized by

$$\begin{aligned}[0, d_k] &\to \widehat{\mathbf{H}}^{2,n} \\ t &\mapsto \cos(t)\hat{x}_k + \sin(t) n_k,\end{aligned}$$

with $n_k$ a unit vector tangent to $\widehat{\mathbf{H}}^{2,n}$ at the point $\hat{x}_k$ and orthogonal to $\mathsf{T}_{\hat{x}_k}\hat{S}$. We know by the previous step that the number $d_k$, which equals the timelike distance between $\hat{x}_k$ and $\hat{y}_k$, tends to 0 as $k$ goes to $+\infty$.

The surface $\hat{S}$ lies inside a totally geodesic copy of $\widehat{\mathbf{H}}^{2,1}$, that we call abusively $\widehat{\mathbf{H}}^{2,1}$. This copy defines an orthogonal splitting $E = F \oplus G$ where $F$ is of signature $(2,2)$ and $G$ is of signature $(0, n-1)$. Let us decompose orthogonaly $n_k$ as

$$n_k = n_k^{\top} + n_k^{\perp} \ ,$$

where $n_k^\top$ is the part of $n_k$ in $F$ (tangent to $\widehat{\mathbf{H}}^{2,1}$) and $n_k^\perp$ the part of $n_k$ in $G$ (normal to $\widehat{\mathbf{H}}^{2,1}$). Since $-1 = \mathbf{q}(n_k) = \mathbf{q}(n_k^\top) + \mathbf{q}(n_k^\perp)$, we have

$$-1 \leq \mathbf{q}(n_k^\top) \leq 0 \quad \text{and} \quad -1 \leq \mathbf{q}(n_k^\perp) \leq 0 \ .$$

By continuity, and since $n_k^\perp$ lies in the unit $\mathbf{q}$-ball of $G$,

$$\lim_{k\to+\infty} \sin(d_k) n_k^\perp = 0 \ .$$

If $n_k^\top \neq 0$ for an infinite number of indices $k$ in $\mathbb{N}$, we extract a subsequence where $n_k^\top \neq 0$. By Proposition 4.1.2, the vector $n_k^\top$ satisfies

$$\frac{n_k^\top}{\mathbf{q}(n_k^\top)} = \pm(\exp(l_k)z_1 - \exp(m_k)z_2 + \exp(-l_k)z_3 - \exp(-l_k)z_4) \ .$$

Since $(x_k)_{k\in\mathbb{N}}$ tends to the point $p$ in $\gamma$, there exist positive real numbers $c_k$ for each integer $k$, and $z = az_1 + bz_2 + cz_3$ a vector in $p \setminus \{0\}$, such that

$$\lim_{k\to+\infty} c_k \hat{x}_k = z \ .$$

Then, we have the limit

$$\lim_{k\to+\infty} c_k \frac{n_k^\top}{\mathbf{q}(n_k^\top)} = \pm(az_1 - bz_2 + cz_3) \ .$$

In each case, by continuity,

$$c_k \hat{y}_k = c_k \cos(d_k)\hat{x}_k + c_k \sin(d_k) n_k^\top + c_k \sin(d_k) n_k^\perp \underset{k\to+\infty}{\longrightarrow} z$$

because the first term tends to $z$ and the other ones to zero. This concludes the second step, and the proof of the proposition. □

**Corollary 4.3.8.** *Keeping the same notations, the radial projection $\rho : S \to \Sigma$ realizes a homeomorphism from a neighborhood $\mathcal{U}$ of $v_2$ in $\bar{S}$ to a neighborhood $\mathcal{V}$ of $v_2$ in $\bar{\Sigma}$. The neighborhood $\mathcal{U}$ can be chosen to contain any given point $p$ of $\gamma$.*

**Proof.** By the proof of Proposition 4.3.7, we know that the timelike distance between a point $x$ in $S$ and $\rho(x)$ is as small as we want when approaching $\gamma$ in $\bar{S}$. Let $p$ be a point of $\gamma$. We chose a neighborhood $\mathcal{U}$ of $v_2$ in $\bar{S}$, homeomorphic to a disc and containing $p$, such that for every $x$ in $\mathcal{U} \cap S$

$$\delta^{\mathrm{t}}(x, \rho(x)) \leq \arcsin\left(\frac{1}{\sqrt{5}}\right) \ .$$

By Proposition 4.3.4, the radial projection $\rho$ is a local diffeomorphism from $\mathcal{U} \cap S$ to its image. Moreover, $\rho$ is proper in $\mathcal{U}$.

At the end, $\rho_{|\mathcal{U}}$ is a continuous and proper map, and is a local diffeomorphism in $\mathcal{U} \cap S$. Hence it is a homeomorphism from $\mathcal{U}$ to its image. □

*4.3.4. Renormalization*

We are in this paragraph again considering a maximal surface $\Sigma$ with boundary the semi-positive loop $\Lambda$ such that (see Fig. 6):

- the loop $\Lambda$ contains a vertex $v_1$,
- $C$ is a Barbot crown that coincides with $\Lambda$ along the two edges adjacent to $v_2$,
- $S$ is the Barbot surface with space boundary $C$,
- the vertices of $C$ are denoted $v_1$, $v_2$, $v_3$ and $v_4$,
- $\gamma$ is the connected component of $\Lambda \setminus \{v_1, v_3\}$ containing $v_2$,
- $\bar{\Sigma}$ and $\bar{S}$ are the $\mathcal{C}^0$ surfaces with boundary $\Sigma \cup \gamma$ and $S \cup \gamma$, respectively,
- $\rho : S \to \Sigma$ denotes the radial projection (see Definition 4.3.1),
- $\mathsf{A}$ is the Cartan subgroup of $\mathsf{G}$ determined by the crown $C$ (see Section 3.2.5).

**Lemma 4.3.9** *(Renormalization). Let $(y_k)_{k\in\mathbb{N}}$ be a sequence of elements of $\Sigma$ tending to a point $p$ in $\gamma$ and let $x$ be a point of $S$. Then there is a sequence $(a_k)_{k\in\mathbb{N}}$ of elements of $\mathsf{A}$ such that $((a_k \cdot \Sigma, a_k \cdot y_k))_{k\in\mathbb{N}}$ subconverges to the pointed Barbot surface $(S, x)$.*

By Corollary 4.3.8, the radial projection $\rho : S \to \Sigma$ realizes a homeomorphism from a neighborhood $\mathcal{U}$ of $v_2$ in $\bar{S}$ to a neighborhood $\mathcal{V}$ of $v_2$ in $\bar{\Sigma}$. Denote by $\bar{\rho}$ the obtained homeomorphism from $\mathcal{U}$ to $\mathcal{V}$. Let $p$ be a point of $\gamma$. Without loss of generality, we can suppose that $\mathcal{U}$ contains the point $p$.

**Proof.** Fix the point $x$ in $S$. This defines a hyperbolic basis $(z_1, z_2, z_3, z_4)$ such that $z_i$ are vectors in $v_i$ and $x = [z_1 + z_2 + z_3 + z_4]$. Since $(y_k)_{k\in\mathbb{N}}$ converges to $p$ in $\gamma$, there is an integer $K$ such that the point $y_k$ lies in $\mathcal{V}$ for all $k \geq K$. Up to considering $(y_{k+K})_{k\in\mathbb{N}}$, we assume that $K = 0$. Hence we can define for every integer $k$

$$x_k := \bar{\rho}^{-1}(y_k) \ .$$

For every integer $k$,

$$x_k = a(l_k, m_k) \cdot x \ ,$$

with $a(l_k, m_k)$ elements of the Cartan subgroup $\mathsf{A}$ written in the hyperbolic basis $(z_1, z_2, z_3, z_4)$. The point $p$ is in $\gamma$. As a consequence,

$$[l_k : m_k : l_k^{-1} : m_k^{-1} : 0] = [l_k m_k^{-1} : 1 : l_k^{-1} m_k^{-1} : m_k^{-2} : 0] \to [a : 1 : b : 0 : 0] \ ,$$

for $a$, $b$ nonnegative real numbers. Thus, $(l_k)_{k\in\mathbb{N}}$ and $(m_k)_{k\in\mathbb{N}}$ satisfy the following conditions

(1) the sequence $(m_k)_{k\in\mathbb{N}}$ tends to $+\infty$,
(2) the sequence $(l_k m_k)_{k\in\mathbb{N}}$ is bounded from below,
(3) the sequence $(l_k^{-1} m_k)_{k\in\mathbb{N}}$ is bounded from below.

Denote by $(a_k)_{k\in\mathbb{N}}$ the sequence of elements of $\mathsf{A}$ defined by $a_k := a(l_k^{-1}, m_k^{-1})$ for every integer $k$. By the renormalization lemma (Lemma 3.4.1), the sequence of polygons $(a_k \cdot \Lambda)_{k\in\mathbb{N}}$ subconverges to $C$. Moreover, since for all $k$ in $\mathbb{N}$, $a_k \cdot x_k = x$ and $y_k$ is the radial projection of $x_k$, the point $a_k \cdot y_k$ is in the timelike sphere $\exp(\mathsf{N}_x S)$. By Theorem 2.2.8, the sequence of maximal surfaces $(a_k \cdot \Sigma)_{k\in\mathbb{N}}$ subconverges to the maximal surface with space boundary $C$, that is $S$. The limit point $\lim a_k \cdot y_k$ is $x$ since $\exp(\mathsf{N}_x S)$ is closed and intersects $S$ in the unique point $x$. □

**Corollary 4.3.10.** *Polygonal surfaces are asymptotically flat.*

**Proof.** Let $\Sigma$ be polygonal. The Renormalization Lemma 4.3.9 shows that for every diverging sequence $(y_k)_{k\in\mathbb{N}}$ in $\Sigma$, the sequence $((\Sigma, y_k))_{k\in\mathbb{N}}$ subconverges, up to renormalization, to a pointed Barbot surface $(S, x)$. Hence,

$$\lim_{k\to+\infty} K_\Sigma(y_k) = K_S(x) = 0 . \quad \square$$

### *4.4. Space-horofunctions*

The results of this section are the followings. Let $\Sigma$ be a maximal surface of $\mathbf{H}^{2,n}$ with space boundary the semi-positive loop $\Lambda$. For each point $p$ in $\Lambda$, the space-horofunction associated with $p$ satisfies the following two properties:

- the function $h_p$ is quasiconvex, and strictly quasiconvex if $\Sigma$ is not a Barbot surface,
- the horocycles $\{h_p = C\}$ are curves that tend, in $\partial^{\mathrm{s}}\Sigma$, to the two points at the end of an edge with $p$ (see Fig. 7).

#### *4.4.1. Quasiconvexity*

Let $\Sigma$ be a Hadamard surface.

**Definition 4.4.1.** A function $h$ on $\Sigma$ is said to be *(strictly) quasiconvex* if its sublevel sets are (strictly) convex subsets of the surface.

Because of the definition of a convex set, the strict quasiconvexity of $h$ is equivalent to the following. For every geodesic segment $c$ defined in $[0, 1]$, for every $t$ in $(0, 1)$, we have the *strict quasiconvexity property*

$$h \circ c(t) < \max(h \circ c(0), h \circ c(1)) \ , \tag{7}$$

saying that $h \circ c$ is a strictly quasiconvex function on $\mathbb{R}$. Moreover, we have a property easiest to verify when $h$ is a $\mathcal{C}^2$ function.

**Proposition 4.4.2.** *Let $h$ be a function on $\Sigma$. If $h$ is $\mathcal{C}^2$, then $h$ is strictly quasiconvex if for every nonzero critical direction $X$, meaning $\mathrm{d}h \cdot X = 0$, the Hessian of $h$ satisfies $\nabla^2 h \cdot (X, X) > 0$.*

By the strict quasiconvexity property (7), the statement is equivalent to the following statement in the real line.

**Lemma 4.4.3.** *Let $f$ be a function on $\mathbb{R}$. If $f$ is $\mathcal{C}^2$, then $f$ is strictly quasiconvex if for every critical point $t$, meaning $f'(t) = 0$, the second derivative of $f$ satisfies $f''(t) > 0$.*

**Proof.** Let $f$ be a $\mathcal{C}^2$ function on $\mathbb{R}$. A point $t$ satisfying $f'(t) = 0$ and $f''(t) > 0$ is a strict local minimum. If $f$ satisfies the assumption

$$f'(t) = 0 \Rightarrow f''(t) > 0 \ ,$$

then every critical point of $f$ is a strict local minimum, hence $f$ has at most one critical point which is a strict minimum. This implies that $f$ is strictly monotone in $\mathbb{R}$ or in $(-\infty, t_0]$ and $[t_0, +\infty)$ where $t_0$ is a strict minimum. Hence $f$ is strictly quasiconvex. □

*4.4.2. Quasiconvexity of space-horofunctions*

**Theorem 4.4.4.** *Let $\Sigma$ be a maximal surface of $\mathbf{H}^{2,n}$. Take a point $p$ in the space boundary $\Lambda$ of $\Sigma$. Then any space-horofunction $(h_p)_{|\Sigma}$ is quasiconvex. If $\Sigma$ is not a Barbot surface, then any space-horofunction $(h_p)_{|\Sigma}$ is strictly quasiconvex.*

For the proof, we denote by $h$ the restriction of $h_p$ to $\Sigma$. We now compute the gradient and Hessian of the space-horofunction $h$ with respect to the induced metric on $\Sigma$.

**Proposition 4.4.5** *(Gradient and Hessian of space-horofunctions). Let $\Sigma$ be a spacelike surface of $\mathbf{H}^{2,n}$ and $p$ be a point of $\partial^{\mathrm{s}}\Sigma$. The gradient and the Hessian of the space-horofunction $h$ associated to $p$ and restricted to $\Sigma$ are at $x$, for vectors $X, Y$ in $\mathsf{T}_x\Sigma$:*

- $\mathrm{grad}(h)_x = z^\top / \langle x, z \rangle$,
- $\nabla^2_x h \cdot (X, Y) = \langle X, Y \rangle_x + \frac{\langle \mathrm{I\!I}(X,Y)_x, z \rangle}{\langle x, z \rangle} - (X \cdot h)_x (Y \cdot h)_x$,

*where $z$ is a nonzero vector in $p$, and $z^\top$ is the part of $z$ tangent to $\mathsf{T}_x\Sigma$ when we decompose orthogonaly $\mathbf{E} = x \oplus \mathsf{T}_x\Sigma \oplus \mathsf{N}_x\Sigma$.*

**Proof.** We fix a lift $\hat{\Sigma}$ of $\Sigma$ in $\widehat{\mathbf{H}}^{2,n}$ to make computations and chose coherently a nonzero vector $z$ in $p$. We fix the representative of $h$ that takes the value

$$h(x) = \log(-\langle x, z\rangle) \ ,$$

at any point $x$ in $\hat{\Sigma}$. Hence, if $X$ is a tangent vector of $\hat{\Sigma}$ at $x$,

$$\mathrm{d}_x h \cdot X = \frac{\langle X, z\rangle}{\langle x, z\rangle} \ ,$$

and the gradient of $h$ equals

$$\mathrm{grad}(h)_x = \frac{z^\top}{\langle x, z\rangle} \ .$$

Denote by $D$ the flat connection on $\mathbf{E}$. It splits, at a point $x$, as

$$D_X Y = \langle X, Y\rangle x + \nabla_X Y + \mathbb{I}(X, Y) \ ,$$

in the splitting $\mathbf{E} = x \oplus \mathsf{T}_x\hat{\Sigma} \oplus \mathsf{N}_x\hat{\Sigma}$, where $\nabla$ is the Levi-Civita connection of $\hat{\Sigma}$ and $\mathbb{I}$ the second fundamental form of $\hat{\Sigma}$ in $\widehat{\mathbf{H}}^{2,n}$.

To compute the Hessian, recall that for any $X$, $Y$ vector fields on $\hat{\Sigma}$,

$$\nabla^2 h \cdot (X, Y) = X \cdot (Y \cdot h) - \nabla_X Y \cdot h \ .$$

We first compute $X \cdot (Y \cdot h)$ for two vector fields $X, Y$ on $\hat{\Sigma}$. It gives

$$\begin{aligned} X \cdot (Y \cdot h)(x) =& X \cdot \frac{\langle Y, z\rangle}{\langle x, z\rangle} \\ =& \frac{\langle x, z\rangle\langle D_X Y, Z\rangle - \langle Y, z\rangle\langle X, z\rangle}{\langle x, z\rangle^2} \\ =& \frac{\langle\langle X, Y\rangle_x x + (\nabla_X Y)_x + \mathbb{I}_x(X, Y), z\rangle}{\langle x, z\rangle} - \frac{\langle Y, z\rangle_x\langle X, z\rangle_x}{\langle x, z\rangle^2}. \end{aligned}$$

Hence

$$\nabla^2_x h \cdot (X, Y) = \langle X, Y\rangle_x + \frac{\langle \mathbb{I}_x(X, Y), z\rangle}{\langle x, z\rangle} - (X \cdot h)_x(Y \cdot h)_x \ . \quad \square$$

**Proof of Theorem 4.4.4.** The space-horofunctions in Barbot surfaces have been computed in Section 4.1.4. They are either linear (when associated to a vertex) or strictly quasiconvex.

By [19, Lemma 5.12], the function

$$\beta : (x, X) \mapsto \frac{\langle \mathbb{I}_x(X, X), z\rangle}{\langle x, z\rangle} \ ,$$

defined on the unit tangent bundle of a maximal surface $\Sigma$, satisfies

$$\beta(x,\cdot)^2 \le (\mathbf{q}(\mathrm{grad}(h)_x) - 1)(1 + K_\Sigma(x)) \ .$$

By [19, Theorem 5.5], the norm of the gradient of a space-horofunction satisfies the inequalities:

$$0 < \mathbf{q}(\mathrm{grad}(h)) \le 2. \tag{8}$$

Since $\mathbf{q}(\mathrm{grad}(h)) \le 2$ and $K_\Sigma \le 0$, we have

$$\beta \ge -1 \ .$$

For a maximal surface that is not a Barbot surface, the strict inequality $K_\Sigma < 0$ gives

$$\beta > -1 \ .$$

Putting these inequalities in the formula of the Hessian of $h$ we obtain that if $X$ is a critical direction at a point $x$, on a Barbot surface

$$\nabla^2_x h \cdot (X, X) \ge 0 \ ,$$

and if $\Sigma$ is not a Barbot surface,

$$\nabla^2_x h \cdot (X, X) > 0 \ ,$$

making $h$ strictly quasiconvex by Proposition 4.4.2. □

*4.4.3. Space-horoballs and space-horocycles*

Let $\Sigma$ be a maximal surface with space boundary the semi-positive loop $\Lambda$. We study in this paragraph the shape of space-horoballs and space-horocycles associated to a point $p$ in $\Lambda$.

**Definition 4.4.6.** *Space-horoballs* and *space-horocycles* are sublevel sets and level sets of a space-horofunction, respectively, in a maximal surface of $\mathbf{H}^{2,n}$.

The geometry of space-horoballs is well understood on hyperbolic planes where they are exactly the intrinsic horoballs, and on Barbot surfaces (see Paragraph 4.1.4). In general, as an immediate consequence of Proposition 4.4.5, we have the following.

**Lemma 4.4.7.** *A space-horofunction has no critical point. Hence the space-horopheres are smooth simple and unbounded curves that foliate the maximal surface.*

Now we want to understand how the space-horocycles behave with respect to the space-boundary.

**Proposition 4.4.8.** *The points at infinity of a space-horocycle $\{h_p = C\}$ satisfy*

$$[q \in \partial^{\mathrm{s}}\{h_p = C\}] \Rightarrow [\langle q, p\rangle = 0] \ .$$

**Proof.** Denote $H = \{h_p = C\}$ the level set for the value $C$ in $\mathbb{R}$. The closure $\bar{H}$ of $H$ in $\bar{\Sigma} = \Sigma \cup \Lambda$ meets $\Lambda$. We want to show that the points in the frontier $\bar{H} \setminus H$ are in a photon with $p$. Let $(x_k)_{k\in\mathbb{N}}$ be a sequence of points of $H$ tending to a point $q$ at infinity. Take a lift $(\hat{x}_k)_{k\in\mathbb{N}}$ in $\widehat{\mathbf{H}}^{2,n}$. There exists a sequence of non-zero real numbers $(\lambda_k)_{k\in\mathbb{N}}$, tending to 0, such that $(\lambda_k\hat{x}_k)_{k\in\mathbb{N}}$ tends to $w$, a vector in $q \setminus \{0\}$, as $k \to +\infty$. Denote by $z$ a non-zero vector in $p$ such that $h(x) = \log(|\langle \hat{x}, z\rangle|)$ is a representative of $h_p$. Now we have

$$|\langle \lambda_k \hat{x}_k, z\rangle| = |\lambda_k||\langle \hat{x}_k, z\rangle| = |\lambda_k| \exp(C) \ .$$

Taking the limit $k \to +\infty$ we obtain $\langle w, z\rangle = 0$, so $q$ is in a photon with $p$. □

*4.4.4. The positive case*

As a corollary of the last proposition, we know the shape of space-horoballs in a particular case.

**Definition 4.4.9.** A point $p$ in a semi-positive loop $\Lambda$ is said to be a *positive point* if it does not belong to any edge in $\Lambda$.

**Corollary 4.4.10.** *Let $\Sigma$ be a maximal surface of $\mathbf{H}^{2,n}$, with space boundary the semi-positive loop $\Lambda$. If $p$ is a positive point in $\Lambda$, then the space-horoballs associated to $p$ are strictly convex subsets of $\Sigma$. The frontier of a space-horoball associated to $p$ in $\Sigma \cup \Lambda$ is a circle that intersect $\Lambda$ only at $p$ (see Fig. 7).*

*4.4.5. The vertex case*

We will need more work to understand the shape of space-horoballs associated to vertices. Recall that we call a vertex, in a semi-positive loop $\Lambda$, a point that is at the intersection of two distinct edges in $\Lambda$. If $v$ is a vertex, it is in the extremity of two different edges and we denote by $v_-$ and $v_+$ the two opposite extremities.

**Proposition 4.4.11.** *Let $\Sigma$ be a maximal surface of $\mathbf{H}^{2,n}$, with space boundary the semi-positive loop $\Lambda$. Suppose that $\Sigma$ is not a Barbot surface. If $v$ is a vertex of $\Lambda$, then the space-horoballs associated to $v$ are strictly convex subsets of $\Sigma$ and their boundaries, space-horocycles associated to $v$, have boundary $\{v_+, v_-\}$ at infinity.*

*Moreover, if $v_+$ is not a vertex (meaning that the only edge containing $v_+$ is the edge $[v_+, v]$), then the space-horoballs associated to $v_+$ are strictly convex subsets of $\Sigma$ and*

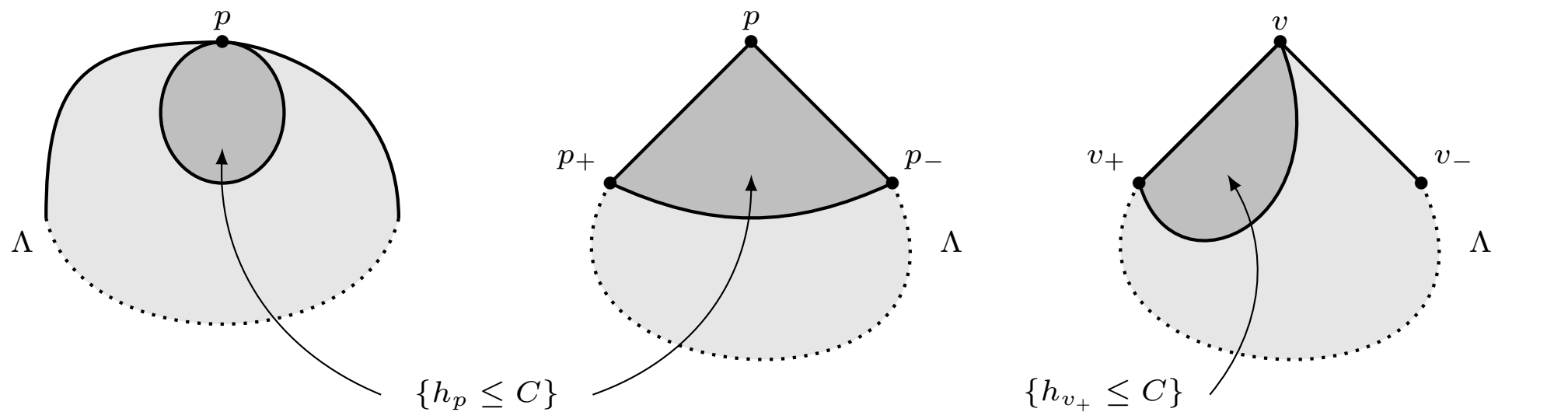

**Fig. 7.** Shape of space-horoballs in a maximal surface. On the left the situation where $p$ is a positive point. On the middle the situation where $p$ is a vertex. On the right the situation where $v_+$ is not a vertex but is contained in an edge with a vertex $v$.

*their boundaries, space-horocycles associated to $v_+$, have boundary $\{v_+, v\}$ at infinity (see Fig. 7).*

**Proof.** The claims about strict convexity of space-horoballs are the content of Theorem 4.4.4. We now study the space-horocycles.

By Corollary 3.2.8, there is a Barbot crown $C$ in $\mathbf{Ein}^{1,n}$ that coincide with $\Lambda$ along the two edges containing $v$. The vertices of $C$ are then $v_-$, $v$, $v_+$ and a fourth vertex $v_4$ of $\mathbf{Ein}^{1,n}$. Take a hyperbolic basis $(z_1, z_2, z_3, z_4)$ such that $z_1$ lies in $v_-$, $z_2$ lies in $v$, $z_3$ lies in $v_+$ and $z_4$ lies in $v_4$. This fixes a representative $h$ of $h_v$ that takes the value

$$h(y) = \log(-\langle \hat{y}, z_2 \rangle) \ ,$$

at points $y$ of $\Sigma$, where $\hat{y}$ is a lift of $y$ to the double cover $\widehat{\mathbf{H}}^{2,n}$ satisfying $\langle \hat{y}, z_2 \rangle < 0$. This also fixes a representative $h_+$ of $h_{v_+}$ that takes the value

$$h_+(y) = \log(-\langle \hat{y}, z_3 \rangle) \ ,$$

at points $y$ of $\Sigma$, where $\hat{y}$ is a lift of $y$ to the double cover $\widehat{\mathbf{H}}^{2,n}$ satisfying $\langle \hat{y}, z_3 \rangle < 0$.

Take a sequence of points $(y_k)_{k\in\mathbb{N}}$ in $\Sigma$ tending to a point $p = [tz_2 + (1-t)z_3]$ in the interior of the edge $[v_+, v]$, that is $t$ lies in $(0, 1)$. Denote by $S$ the Barbot surface with space boundary $C$ and by $x$ a point of $S$. By the Renormalization Lemma 4.3.9, there is a sequence of elements $(a_k)_{k\in\mathbb{N}}$ of the Cartan subgroup $\mathsf{A}_0$ associated with $C$ such that the sequence of pointed maximal surfaces $((a_k \cdot \Sigma, a_k \cdot y_k))_{k\in\mathbb{N}}$ converges to the pointed Barbot surface $(S, x)$. Now, by looking at the parametrization of $S$ using $(z_1, z_2, z_3, z_4)$ (see Section 4.1.1), we see that the sequence $(a_k)_{k\in\mathbb{N}}$ is written $a_k = a(\lambda_k, \mu_k)$ with

- the sequence $(\mu_k)_{k\in\mathbb{N}}$ tends to 0,
- the sequence $(\lambda_k)_{k\in\mathbb{N}}$ diverges to $+\infty$.

We compute

$$h(y_k) = \log\left(-\langle \hat{y}_k, z_2 \rangle\right)$$

$$
\begin{aligned}
&= \log\left(-\langle a_k \cdot \hat{y}_k, a_k \cdot z_2\rangle\right) && (a_k \text{ is an isometry})\\
&= \log\left(-\langle a_k \cdot \hat{y}_k, z_2\rangle\right) + \log(\mu_k) && (a_k \cdot z_2 = \mu_k z_2)\\
&= h(a_k \cdot y_k) + \log(\mu_k)\ ,
\end{aligned}
$$

and

$$
\begin{aligned}
h_+(y_k) &= \log\left(-\langle \hat{y}_k, z_3\rangle\right)\\
&= \log\left(-\langle a_k \cdot \hat{y}_k, a_k \cdot z_3\rangle\right) && (a_k \text{ is an isometry})\\
&= \log\left(-\langle a_k \cdot \hat{y}_k, z_3\rangle\right) - \log(\lambda_k) && (a_k \cdot z_3 = \lambda_k^{-1} z_3)\\
&= h_+(a_k \cdot y_k) - \log(\lambda_k)\ .
\end{aligned}
$$

By continuity of $h$ and $h_+$, the sequences $(h(a_k \cdot y_k))_{k\in\mathbb{N}}$ and $(h_+(a_k \cdot y_k))_{k\in\mathbb{N}}$ tend to the real numbers $h(x)$ and $h_+(x)$, respectively. As a consequence,

$$
\lim_{k\to+\infty} h(y_k) = -\infty \quad \text{and} \quad \lim_{k\to+\infty} h_+(y_k) = -\infty\ .
$$

By the same argument, for every sequence $(y_k)_{k\in\mathbb{N}}$ of points of $\Sigma$ tending to a point $p$ in the interior of the edge $[v, v_-]$ (hence the new $(\lambda_k)_{k\in\mathbb{N}}$ tends to 0), we have

$$
\lim_{k\to+\infty} h(y_k) = -\infty \quad \text{and} \quad \lim_{k\to+\infty} h_+(y_k) = +\infty\ .
$$

Since the endpoints of space-horocycles associated to $v$ must be in an edge with $v$ (Proposition 4.4.8), we obtain that the endpoints of space-horocycles associated to $v$ are $v_+$ and $v_-$, and that the endpoints of space-horocycles associated to $v_+$ are $v_+$ and $v$. □

We are now interested in comparing the space-horocycles associated to $v_+$ and $v_-$, near $v$.

**Corollary 4.4.12.** *Let $\Sigma$ be a maximal surface of $\mathbf{H}^{2,n}$ with space boundary the semi-positive loop $\Lambda$. Suppose that $\Sigma$ is not a Barbot surface. Suppose that $\Lambda$ contains a vertex $v$, and denote by $v_+$ and $v_-$ the two points at the end of edges with $v$. Then the space-horoballs associated to $v_+$ and $v_-$ are either disjoint or their intersection is a nonbounded strictly convex set containing only $v$ at infinity. Moreover, if we parametrize by $\gamma$ a space-horocycle associated to $v_-$, such that $\gamma(t)$ tends to $v$ when $t$ tends to $+\infty$, then the distance between $\gamma_{|[0,+\infty)}$ and any space-horocycle associated to $v_+$ is bounded by a constant depending only on the chosen space-horocycles.*

**Proof.** By Theorem 4.4.4, the space-horoballs associated to $v_-$ and $v_+$, respectively, are strictly convex subsets of $\Sigma$. If two space-horoballs intersect, the intersection is a strictly convex subset as being the intersection of two strictly convex subspaces. Let $H_+$ and $H_-$ be space-horoballs associated to $v_+$ and $v_-$, respectively.

Suppose by contradiction that $H_+ \cap H_-$ is a non-empty bounded subset of $\Sigma$. Up to looking at another space-horofunction associated to $v_-$, instead of $H_-$, we can suppose that the intersection is a unique point. Denote by $c$ a geodesic of $\Sigma$ intersecting $H_+ \cup H_-$ only at the point $c(0) = H_+ \cap H_-$. The geodesic $c$ is then tangent to $H_+$ and to $H_-$. By strict quasi-convexity of the space-horofunctions $h_+$ and $h_-$ associated to $v_+$ and $v_-$, respectively, the functions $h_+ \circ c$ and $h_- \circ c$ are increasing and diverge to $+\infty$. Up to changing orientation, the curve $c_{|[0,+\infty)}$ converges to $v$ at infinity, by the shape of $H_+$ and $H_-$. We want to renormalize along $c(t)$.

By Corollary 3.2.8, there is a Barbot crown $C$ in $\mathbf{Ein}^{1,n}$ that coincides with $\Lambda$ along the two edges $[v_+, v]$ and $[v, v_-]$ containing the point $v$. The vertices of $C$ are $v_-$, $v$, $v_+$ and a fourth point $v_4$ in $\mathbf{Ein}^{1,n}$. Denote by $S$ the Barbot surface with space boundary $C$, and let $x$ be a point of $S$. Also denote by $y_k$ the point $c(k)$ for every integer $k$. By the Renormalization Lemma 4.3.9, there is a sequence $(a_k)_{k\in\mathbb{N}}$ of the Cartan subgroup $\mathsf{A}_0$ associated with $C$ such that the sequence of pointed maximal surfaces $((a_k \cdot \Sigma, a_k \cdot y_k))_{k\in\mathbb{N}}$ converges to the pointed Barbot surface $(S, x)$. Now, take a hyperbolic basis $(z_1, z_2, z_3, z_4)$ corresponding to $C$. This gives us representatives $h_+$ and $h_-$ of the space-horofunctions associated to $v_+$ and $v_-$, respectively. The parametrization of $S$ associated with the hyperbolic basis $(z_1, z_2, z_3, z_4)$ (see Section 4.1.1) shows that the sequence $(\mu_k)_{k\in\mathbb{N}}$ tends to 0. We compute

$$\begin{aligned} h_+(y_k) &= \log\left(-\langle \hat{y}_k, z_3\rangle\right) \\ &= \log\left(-\langle a_k \cdot \hat{y}_k, a_k \cdot z_3\rangle\right) && (a_k \text{ is an isometry}) \\ &= \log\left(-\langle a_k \cdot \hat{y}_k, z_3\rangle\right) - \log(\lambda_k) && (a_k \cdot z_3 = \lambda_k^{-1} z_3) \\ &= h_+(a_k \cdot y_k) - \log(\lambda_k)\ . \end{aligned}$$

Similarly,

$$h_-(y_k) = h_-(a_k \cdot y_k) + \log(\lambda_k)\ . \tag{9}$$

We know that $(h_+(y_k))_{k\in\mathbb{N}}$ and $(h_-(y_k))_{k\in\mathbb{N}}$ tend to $+\infty$, the first computation implies that $(\lambda_k)_{k\in\mathbb{N}}$ tends to 0 and the second implies that $(\lambda_k)_{k\in\mathbb{N}}$ tends to $+\infty$, by continuity of space-horofunctions. This is a contradiction.

Hence, the intersection between $H_+$ and $H_-$ is either empty or unbounded. In the second case, the closure of $H_+ \cap H_-$ in $\Sigma \cup \Lambda$ contains only $v$ because of the shape of $H_+$ and $H_-$ (see Proposition 4.4.11).

We now prove the last assertion of the Corollary. Let $\gamma : \mathbb{R} \to \Sigma$ denote a continuous parametrization of the space-horocycle $\partial H_-$, such that $\gamma(t)$ tends to $v$ when $t$ tends to $+\infty$. We define a sequence of points $(y_k)_{k\in\mathbb{N}}$ by $y_k := \gamma(t_k)$ for an increasing sequence $(t_k)_{k\in\mathbb{N}}$ of real numbers tending to $+\infty$. Since $y_k$ tends to $v$, we can make the same computations as above. Let $(a_k)_{k\in\mathbb{N}}$ be a sequence of elements of $\mathsf{A}_0$ such that $((a_k \cdot \Sigma, a_k \cdot y_k))_{k\in\mathbb{N}}$ converges to the pointed Barbot surface $(S, x)$. Since $h_-(y_k)$ is constant and since $h_-(a_k \cdot y_k)$ tends to $h_-(x)$ by continuity, we have that $a_k = a(\lambda_k, \mu_k)$ for a sequence $(\lambda_k)_{k\in\mathbb{N}}$ converging to $\exp(h_-(y_0) - h_-(x))$ (see Equation (9)). But we also have

$$h_+(y_k) = h_+(a_k \cdot y_k) - \log(\lambda_k) \ ,$$

hence the sequence $(h_+(y_k))_{k\in\mathbb{N}}$ tends to a constant. Since the norm of the gradient of $h_+$ is bounded (by [19, Theorem 5.5]), the sequence $(y_k)_{k\in\mathbb{N}}$ stays at finite distance from the space-horocycle $H_+$. This finishes the proof. □

### *4.5. The link between the space boundary and the ideal boundary of a maximal surface*

In this section, we provide a comparison between the space boundary and the ideal boundary of maximal surfaces. Let $\Sigma$ be a maximal surface of $\mathbf{H}^{2,n}$ with space boundary the semi-positive loop $\Lambda$. To each vertex of $\Lambda$ we associate a *stable direction*, that is a point of $\Sigma(\infty)$, in a natural way. We also associate to each positive point of $\Lambda$ a point of $\Sigma(\infty)$, called *regular*, in a natural way. Here natural means that

- a geodesic ray following a stable direction or a regular direction with respect to the point $p$ of $\Lambda$ converges to $p$ in $\Sigma \cup \Lambda$,
- the construction is $\mathsf{G}$-invariant.

#### *4.5.1. Polygonal case*

Let $v$ be a vertex of $\Lambda$ and $v_-$, $v_+$ be the two points of $\Lambda$ at the end of the two edges that meet $v$. We denote by $h_-$, $h$ and $h_+$ representatives of space-horofunctions associated to $v_-$, $v$ and $v_+$, respectively.

**Definition 4.5.1.** A geodesic ray $c$ of $\Sigma$ is said to be *stable associated to $v$* if $t \mapsto \mathrm{d}(c(t), B)$ is bounded for every space-horoball $B$ associated to $v_-$ or $v_+$.

**Theorem 4.5.2.** *Let $\Sigma$ be a maximal surface with space boundary the semi-positive loop $\Lambda$ and let $v$ be a vertex of $\Lambda$. For each point $x$ in $\Sigma$, there is a unique stable geodesic ray associated to $v$ starting from $x$. Moreover, two stable geodesic rays associated to $v$ are asymptotic.*

**Definition 4.5.3.** The collection of stable rays associated to a vertex $v$ defines, by the previous theorem, a point in $\Sigma(\infty)$, called the *stable ideal point associated to $v$*.

**Proof.** If $\Sigma$ is a Barbot surface, the result follows from our description of space-horoballs in Barbot surfaces (see Remark 4.1.11 and Fig. 4). Suppose now that $\Sigma$ is not a Barbot surface.

By Proposition 4.4.11 and Corollary 4.4.12, the space-horoballs associated to the vertices $v_+$ and $v_-$ are strictly convex subsets of $\Sigma$ and their intersection is either empty or is a nonbounded convex set containing only $v$ at infinity with boundary two curves at finite distance towards $v$. This ensures that two space-horocycles associated to $v_+$ and $v_-$ are not tangent at a point.

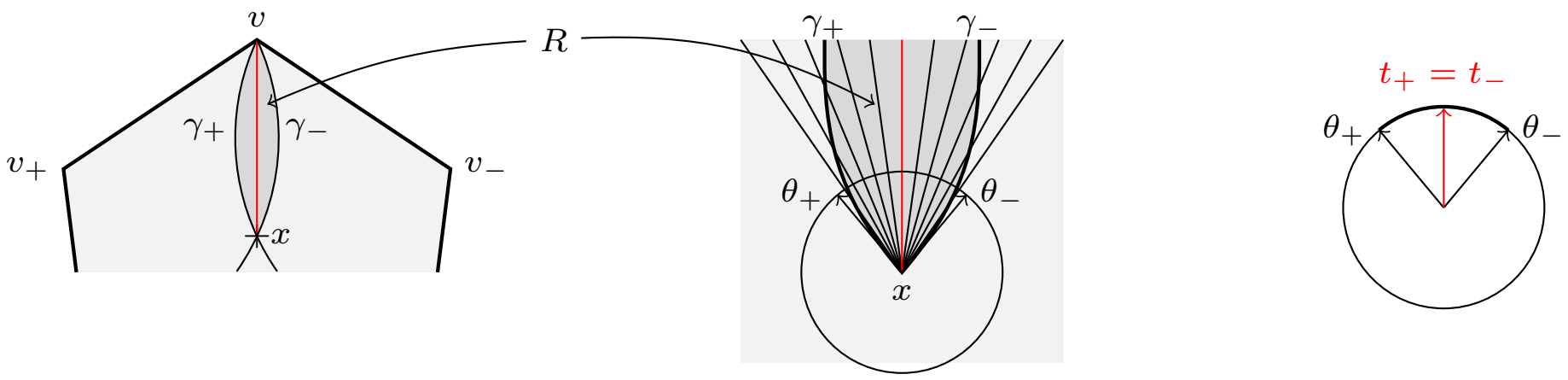

**Fig. 8.** Construction of a stable geodesic ray emanating from $x$.

We now prove the existence of stable geodesic rays associated to $v$. Take $x$ in $\Sigma$. Denote by $\gamma_-$ and $\gamma_+$ the space-horocycles passing through $x$ and associated to $v_-$ and $v_+$, respectively. We restrict $\gamma_+$ and $\gamma_-$ to the part between $x$ and $v$. These curves bound a convex nonbounded region denoted $R$ that is not empty because $\gamma_+$ and $\gamma_-$ are not tangent at $x$. Denote by $\mathbf{S}^1$ the unit circle of $(\mathsf{T}_x\Sigma, \mathbf{q}_{|\mathsf{T}_x\Sigma})$, and parametrize it by arclength $\theta$ in $[0, 2\pi]$. Denote by $\theta_-$ the direction of $\gamma_-$ at $x$ and $\theta_+$ the direction of $\gamma_+$ at $x$ such that $\theta_- \leq \theta_+$ and every direction $\theta$ in $(\theta_-, \theta_+)$ points inside $R$. The angles $\theta_-$ and $\theta_+$ are different because $\gamma_-$ and $\gamma_+$ are not tangent at $x$ (by strict convexity of the horoballs). The set of geodesics rays starting from $x$, pointing inside $R$ and crossing $\gamma_- \setminus \{x\}$ corresponds to an interval $(\theta_-, t_-)$ (it is defined by an open condition). Similarly, the set of geodesic rays starting from $x$, pointing inside $R$ and crossing $\gamma_+ \setminus \{x\}$ corresponds to an interval $(t_+, \theta_+)$. See Fig. 8 for a picture of the construction.

**Claim**: $t_- \leq t_+$.

Indeed, if $\theta_- < \theta < t_-$, the corresponding ray $c_\theta$ crosses $\gamma_-$, so it goes out of the sublevel set $\{h_- \leq h_-(x)\}$ and never goes inside after, by convexity. So $c_\theta$ cannot cross $\gamma_+$, and $\theta < t_+$. Hence $t_- \leq t_+$.

Every geodesic ray emanating from $x$ with a direction in $[t_-, t_+]$ stays inside $R$, so converges to $v$ at infinity, by staying at finite distance of $\gamma_-$ and $\gamma_+$, hence is stable.

Now, since $\gamma_-$ and $\gamma_+$ are at bounded distance and $\Sigma$ is negatively curved, $t_- = t_+$. This proves existence and uniqueness.

To prove that the collection of stable rays defines an ideal point, first remark that if we take a geodesic ray $c$ with initial direction $\dot{c}(0)$ a stable direction, then the velocity of the geodesic ray $\dot{c}(t)$ is a stable direction for every $t > 0$, by definition. Finally, two stable geodesic rays tending to the same vertex $v$ stay at finite distance, by definition, so they are asymptotic and the collection of stable geodesic rays tending to $v$ is a point of the ideal boundary $\Sigma(\infty)$. □

**Remark 4.5.4.** The above construction is $\mathsf{G}$ invariant because $\mathsf{G}$ sends space-horoballs to space-horoballs and geodesics to geodesics hence stable geodesic rays to stable geodesic rays.

*4.5.2. Positive case*

For this paragraph, $\Sigma$ is a maximal surface with boundary the semi-positive loop $\Lambda$. Suppose moreover that $\Lambda$ contains *positive* points, that is, points meeting no edge in $\Lambda$.

We can now assign to every positive point $p$ in $\Lambda$ a collection of ideal points in a natural way, similar to the construction of stable geodesic rays associated to a vertex.

**Definition 4.5.5.** Let $\Sigma$ be a maximal surface in $\mathbf{H}^{2,n}$ and $p$ be a positive point of $\partial^{\mathrm{s}}\Sigma$. A geodesic ray $c$ of $\Sigma$ is said to be *regular* with respect to $p$ if it stays inside the space-horoball $\{h_p \leq h_p(c(0))\}$. An ideal boundary point $\xi$ of $\Sigma(\infty)$ is said to be *regular* with respect to $p$ if every geodesic ray of $\Sigma$ representing $\xi$ is regular with respect to $p$.

**Theorem 4.5.6.** *Let $\Sigma$ be a maximal surface in $\mathbf{H}^{2,n}$. For every positive point $p$ in $\partial^{\mathrm{s}}\Sigma$, there is an ideal point $\xi$ in $\Sigma(\infty)$ that is regular with respect to $p$.*

**Proof.** Suppose that $\partial^{\mathrm{s}}\Sigma$ contains some positive points. In particular, $\Sigma$ is not a Barbot surface. By Corollary 4.4.10, the space-horoballs associated to $p$ are strictly convex non bounded subsets of $\Sigma$ containing only $p$ at infinity.

Fix a point $x_0$ of $\Sigma$, and look at the set of geodesic rays starting from $x_0$, pointing inside the space-horoball $\{h_p \leq h_p(x_0)\}$ and crossing the space-horocycle $\{h_p = h_p(x_0)\}$. The collection of such rays is the union of two intervals in the unit tangent space of $\Sigma$ at $x_0$. Between these intervals is the direction of at least one ray $c_0$ staying inside the space-horoball (for every positive $t$, $c_0(t)$ is in the space-horoball $\{h_p \leq h_p(x_0)\}$). Denote by $\xi$ the associated ideal point ($\xi := c_0(\infty)$). By Corollary 4.4.10, the closure of the space-horoball $\{h_p \leq h_p(x_0)\}$ in $\mathbf{H}^{2,n} \cap \mathbf{Ein}^{1,n}$ intersects $\partial^{\mathrm{s}}\Sigma$ only at the point $p$. Hence the ray $c_0$ converges to $p$ at infinity.

Now we prove that every ray with direction $\xi$ is regular with respect to $p$, that is, stays inside space-horoballs associated to $p$. Suppose that $c$ is a geodesic ray of $\Sigma$. If $c$ goes out of the space-horoball $\{h_p \leq h_p(c(0))\}$, the distance between $c(t)$ and $\{h_p \leq h_p(c(0))\}$ goes to $+\infty$ when $t$ goes to $+\infty$. Indeed, denote by $s$ the biggest time for which $c(s)$ lies in $\{h_p = h_p(c(0))\}$. Then $c$ crosses a support geodesic, at $c(s)$, of the strictly convex subset $\{h_p \leq h_p(c(0))\}$. Hence $c$ diverges at least linearly from the space-horoball $\{h_p \leq h_p(c(0))\}$, by nonpositivity of the sectional curvature of $\Sigma$.

Let $c$ be a geodesic ray starting at a point $x$ of $\Sigma$ with direction $\xi$. If $c$ goes out of the space-horoball $\{h_p \leq h_p(x)\}$, the distance between $c(t)$ and the strictly convex subset $\{h_p \leq h_p(x)\}$ diverges to $+\infty$ with $t$ by nonpositivity of the sectional curvature of $\Sigma$. Two space-horoballs associated with the same point are at finite distance (the norm of $\mathrm{grad}(h_p)$ is bounded by above). As a consequence, if $c$ goes out of $\{h_p \leq h_p(x)\}$ then the distance between $c(t)$ and $\{h_p \leq h_p(x_0)\}$ also diverges. But $c_0$ is included in $\{h_p \leq h_p(x_0)\}$, and $c$ and $c_0$ are asymptotic rays. Finally, $c$ cannot go out of the space-horoball $\{h_p \leq h_p(x)\}$ and $\xi$ is regular with respect to $p$. □

**Remark 4.5.7.** Let $\gamma$ be an isometry of $\mathbf{H}^{2,n}$. If $c$ is a geodesic ray regular with respect to the positive point $p$, then $\gamma \cdot c$ will be a geodesic ray of $\gamma \cdot \Sigma$ regular with respect to $\gamma \cdot p$. Indeed, the isometry group sends space-horofunctions to space horofunctions and geodesic rays to geodesic rays.

**Remark 4.5.8.** We do not show that there is a unique ideal boundary point regular with respect to $p$.

### *4.6. The Tits distance*

Let $M$ be a Hadamard manifold, that is, a complete and simply connected Riemannian manifold of nonpositive sectional curvature. We explain in this section the definition of the Tits distance in the ideal boundary $M(\infty)$, and how to compute it in practice. The main reference for this section is [14].

#### *4.6.1. The angle distance*

**Definition 4.6.1.** Let $\xi, \eta$ be two ideal boundary points of $M$. The *angle distance* between $\xi$ and $\eta$, denoted $\measuredangle(\xi, \eta)$, is defined by

$$\measuredangle(\xi, \eta) := \sup_{x \in M} \measuredangle_x(\xi, \eta) \ ,$$

where the angle $\measuredangle_x(\xi, \eta)$ is defined to be the usual Riemannian angle between the two rays $\xi_x, \eta_x$ corresponding to $\xi, \eta$, at $x$, respectively.

The angle distance is indeed a distance in $M(\infty)$.

**Remark 4.6.2.** The angle distance on $M(\infty)$ equals the usual round distance on the unit sphere if and only if $M$ is the euclidean space. In the hyperbolic space, the angle distance is $\pi$ between any two different points of the ideal boundary: indeed two ideal points can always be joined by a geodesic in this case. In general, the angle distance helps to detect flatness or hyperbolicity of Hadamard manifolds.

#### *4.6.2. The l-distance*

Another interesting distance can be defined on $M(\infty)$. The idea is to look at the renormalized distance function restricted to spheres of growing radius. Take two rays $c_1, c_2$ starting from the same point in $M$. Using the triangular inequality, we obtain that for every $0 < s < t$ we have

$$0 \leq \frac{d(c_1(s), c_2(s))}{s} \leq \frac{d(c_1(t), c_2(t))}{t} \leq 2 \ .$$

**Definition 4.6.3.** Let $\xi_1, \xi_2$ be two ideal boundary points of $M$, and $c_1, c_2$ be two rays starting at the same point representing $\xi_1$ and $\xi_2$, respectively. The $l$-distance between $\xi_1$ and $\xi_2$, denoted $l(\xi_1, \xi_2)$, is defined by

$$l(\xi_1, \xi_2) = \lim_{t \to +\infty} \frac{d(c_1(t), c_2(t))}{t} .$$

This definition does not depend on the chosen point, the function $l$ is a distance and is related to the angle distance by the following formula ([14, Lemma 4.4]).

**Proposition 4.6.4.** *If $\xi$ and $\eta$ are two ideal boundary points of $M$, then*

$$l(\xi, \eta) = 2 \sin \left( \frac{\measuredangle(\xi, \eta)}{2} \right) .$$

This formula is the same that links the restriction of the Euclidean distance on a sphere of radius 1 to the spherical distance.

*4.6.3. Inner distances*

We make a brief aside to explain the construction of the *inner distance* associated to a metric space ([14, 4.6]). Let $(A, d)$ be a metric space. If $c : [0, 1] \to A$ is a continuous curve, we can define the *length* of $c$ as being the supremum over all subdivisions $0 = t_0 \leq t_1 \leq \ldots \leq t_k = 1$ of the sums

$$\sum_{i=1}^{k} d(c(t_{i-1}), c(t_i)) .$$

If $x, y$ are in $A$, we define $d_i(x, y)$ the *inner distance* between $x$ and $y$ by

$$d_i(x, y) = \inf\{L(c) \mid c \text{ is a continuous curve joining } x \text{ and } y\} .$$

There could be no continuous curve between $x$ and $y$ or every existing curve between $x$ and $y$ may be of infinite length, in which case $d_i(x, y) = \infty$. A distance $d$ is said to be *inner* if the corresponding inner distance $d_i$ equals $d$.

**Remark 4.6.5.** The distance of a Riemannian manifold is an inner distance. Also in general the inner distance $d_i$ of any distance $d$ is an inner distance ($(d_i)_i = d_i$).

*4.6.4. The Tits distance*

**Definition 4.6.6.** The *Tits distance* on $M(\infty)$, denoted Td, is defined to be the inner distance of the angle distance,

$$\mathrm{Td} = \measuredangle_i .$$

It turns out that the Tits distance is also the inner distance associated with $l$ ([14, Lemma 4.8]). Hence they are two ways to understand the Tits distance on the ideal boundary of a Hadamard manifold, compute the angle distance or compute the $l$-distance. In this article, we use both the angle distance and the $l$-distance, each being more appropriate in some context.

### 4.6.5. Computations

We will need several properties of the Tits distance to make computations.

**Lemma 4.6.7** *(Lemma 4.2 in [14]). Let $\xi$ and $\eta$ be two ideal boundary points of $M$. Then we can compute the angle distance between $\xi$ and $\eta$ as follows:*

$$\measuredangle(\xi,\eta) = \lim_{t\to+\infty} \measuredangle_{c(t)}(\xi,\eta) \ ,$$

*where $c : [0, +\infty)$ is a geodesic ray of direction $\xi$.*

**Lemma 4.6.8** *(Lemma 4.7 in [14]). If $\xi$ and $\eta$ are two ideal boundary points of $M$ such that there is no geodesic $c$ of $M$ satisfying $c(+\infty) = \xi$ and $c(-\infty) = \eta$, then $\mathrm{Td}(\xi,\eta) = \measuredangle(\xi,\eta)$. In particular, if $\measuredangle(\xi,\eta) < \pi$, then $\measuredangle(\xi,\eta) = \mathrm{Td}(\xi,\eta)$.*

**Lemma 4.6.9** *(Lemma 4.10 in [14]). If $\mathrm{Td}(\xi,\eta)$ is larger than $\pi$, then there is a geodesic $c$ in $M$ with $c(+\infty) = \xi$ and $c(-\infty) = \eta$.*

Since $(M(\infty), \mathrm{Td})$ is a metric space, we can ask for the existence and uniqueness of geodesics between two points. In fact, where Td is finite, $(M(\infty), \mathrm{Td})$ is a geodesic space:

**Lemma 4.6.10** *(Lemma 4.11 in [14]). Any two points in $M(\infty)$ with finite Tits distance can be joined by a minimal geodesic of $M(\infty)$. Moreover, the geodesic is unique if the Tits distance is smaller than $\pi$.*

## 4.7. Tits distance in the ideal boundary of a maximal surface

Let $\Sigma$ be a maximal surface with space boundary the semi-positive loop $\Lambda$. To each vertex $v$ of $\Lambda$, we associated a "stable" ideal boundary point. Also to each positive point $p$ of $\Lambda$, we showed the existence of some ideal boundary points regular with respect to $p$. We now compute the Tits distance between these identified ideal boundary points. More precisely, we compute

- the Tits distance between the stable directions associated with two neighboring vertices,
- the Tits distance between two regular rays in a positive interval.

### 4.7.1. The Tits distance between two neighboring vertices

**Theorem 4.7.1** *(Tits distance between neighboring vertices). Let $\Sigma$ be a maximal surface of $\mathbf{H}^{2,n}$ with space boundary the semi-positive loop $\Lambda$. Suppose that $\Lambda$ has two adjacent vertices $v_2$ and $v_3$. The Tits distance between the stable points associated to $v_2$ and $v_3$, respectively, equals $\pi/2$.*

**Proof.** Suppose that $v_2$ and $v_3$ are two adjacent vertices of $\Lambda$. Denote by $v_1$ the point at the end of the edge containing $v_2$ but not $v_3$. Denote by $\xi_2$ and $\xi_3$ the stable points of $\Sigma(\infty)$ associated to $v_2$ and $v_3$, respectively. Take a stable geodesic ray $c$ associated to $v_2$. By the properties of the angle distance (see Lemma 4.6.7), we can compute the angle distance between $\xi_2$ and $\xi_3$ as follows:

$$\measuredangle(\xi_2, \xi_3) = \lim_{t \to +\infty} \measuredangle_{c(t)}(\xi_2, \xi_3) \ .$$

To this aim, we will renormalize along the geodesic $c$. By Corollary 3.2.8, we can find a point $v_4$ of $\mathbf{Ein}^{1,n}$ such that $v_1$, $v_2$, $v_3$ and $v_4$ are the vertices of a Barbot crown $C$ that intersect $\Lambda$ precisely along the two edges adjacent to $v_2$. Denote by $S$ the Barbot surface with space boundary $C$ and by $\mathsf{A}_0$ the Cartan subgroup associated with $C$. By Corollary 4.3.8, the radial projection $\rho$ from $S$ to $\Sigma$ provides a homeomorphism $\bar{\rho}$ from a neighborhood of $v_2$ in $S \cup C$ to a neighborhood of $v_2$ in $\Sigma \cup \Lambda$, and this homeomorphism restricts as the identity on $C$.

Let $(t_k)_{k\in\mathbb{N}}$ be an increasing sequence of positive real numbers tending to $+\infty$. For $k$ big enough, the point $c(t_k)$ is in the range of $\bar{\rho}$. Without loss of generality we suppose that $c(t_k)$ is in the range of $\bar{\rho}$ for every integer $k$. We define $x_k := \bar{\rho}^{-1}(c(t_k))$ for every integer $k$. We will now renormalize along the sequence $(x_k)_{k\in\mathbb{N}}$, that converges to $v_2$ by construction. This is done using the Renormalization Lemma 4.3.9: fix a point $x$ of $S$, then there is a sequence $(a_k)_{k\in\mathbb{N}}$ of elements of $\mathsf{A}_0$ such that the sequence of pointed maximal surfaces $((a_k \cdot \Sigma, a_k \cdot c(t_k)))_{k\in\mathbb{N}}$ converges to $(S, x)$.

Recall that every element of the sequence $(a_k)_{k\in\mathbb{N}}$ stabilizes $v_1$, $v_2$ and $v_3$. At each step the stable ray $c(\cdot + t_k)$, of $\Sigma$, starting from $c(t_k)$ and associated to $v_2$ is sent to the stable geodesic ray denoted $c_k$, of $a_k \cdot \Sigma$, starting at $c_k(0) = a_k \cdot c(t_k)$ and associated to $v_2$. Also at each step, the stable ray starting from $c(t_k)$ and associated to $v_3$ is sent to the stable geodesic ray denoted $\gamma_k$ starting at $\gamma(0) = a_k \cdot c(t_k)$ and associated to $v_3$.

Since $a_k \cdot c(t_k)$ converges to $x$, the sequence of initial velocities $(\dot{c}_k(0))_{k\in\mathbb{N}}$ subconverges to a unit vector $\dot{c}_\infty(0)$ in $\mathsf{T}_x S$. For every integer $k$, the unit vector $\dot{c}_k(0)$ is contained in the cone of $\mathsf{T}_{a_k \cdot c(t_k)} a_k \cdot \Sigma$ of vectors pointing inside the intersection between the space-horoballs $\{h_1 \leq h_1(c_k(0))\}$ and $\{h_3 \leq h_3(\gamma_k(0))\}$ associated with the vertices $v_1$ and $v_3$, respectively. Now, by the very definition of space-horofunctions, a space-horoball in a maximal surface is the intersection of a closed set of $\mathbf{H}^{2,n}$ with this surface. Hence, at the limit, the unit vector $\dot{c}_\infty(0)$ in $\mathsf{T}_x S$ points toward the intersection between the space-horoballs $\{h_1 \leq h_1(x)\}$ and $\{h_3 \leq h_3(x)\}$, in $S$, associated with the vertices $v_1$ and $v_3$, respectively. But this intersection is a stable geodesic of $S$ (see Remark 4.1.11 and Fig. 4). Hence $\dot{c}_\infty(0)$ is the tangent vector, at $x$, of the stable geodesic ray $c_\infty$ associated with $v_2$. By the same argument, the sequence of initial velocities $(\dot{\gamma}_k(0))_{k\in\mathbb{N}}$ converges to the stable direction $\dot{\gamma}_\infty(0)$, at $x$, associated with $v_3$.

By Proposition 4.1.6, $\measuredangle_x(\dot{c}_\infty(0), \dot{\gamma}_\infty(0)) = \pi/2$, hence

$$\measuredangle(\xi_2, \xi_3) = \pi/2 \ .$$

By Lemma 4.6.8, and since the angle distance is smaller than $\pi$, the Tits distance between $\xi_2$ and $\xi_3$ is

$$\mathrm{Td}(\xi_2, \xi_3) = \measuredangle(\xi_2, \xi_3) = \pi/2 \ . \quad \square$$

*4.7.2. The Tits distance in a positive interval*

We say that an interval $I$ in a semi-positive loop $\Lambda$ is *positive* if it contains no edge. Equivalently, every point $p$ in $I$ is a positive point of $\Lambda$.

**Theorem 4.7.2** *(Tits distance between positive points). Let $\Sigma$ be a maximal surface with space boundary the semi-positive loop $\Lambda$. Let $p$, $q$ be two distinct points in an open and positive interval $I$ of $\Lambda$. Then the Tits distance between two ideal points regular with respect to $p$ and $q$, respectively, is infinite.*

**Proof.** Let $I$ be an open and positive interval, and $p$, $q$ be two points of $I$. Denote by $\xi_p$, $\xi_q$ some ideal points that are regular with respect to $p$, $q$, respectively. We will first show that the angle distance between $\xi_p$ and $\xi_q$ is $\pi$.

Let $r$ and $s$ be two positive points of $\Lambda$, one in each connected component of $\Lambda \setminus \{p, q\}$. Let $\mathbf{H}^{2,0}$ be a totally geodesic hyperbolic plane having $r$ and $s$ in its space boundary. Denote by $\tilde{\gamma} : \mathbb{R} \to \mathbf{H}^{2,0}$ the geodesic of $\mathbf{H}^{2,0}$ going from $s$ to $r$ and by $g_t$ the isometry of $\mathbf{H}^{2,n}$ that restricts, in $\mathbf{H}^{2,0}$, to the hyperbolic translation sending $\tilde{\gamma}(t)$ to $\tilde{\gamma}(0)$.

We use the radial projection $\rho$, Definition 4.3.1, between $\mathbf{H}^{2,0}$ and $\Sigma$, to project $\tilde{\gamma}$ to a curve $\gamma := \rho(\tilde{\gamma})$ on $\Sigma$. We claim that $\lim\limits_{t \to +\infty} \gamma(t) = r$. Indeed, take any sequence $(t_k)_{k \in \mathbb{N}}$ of real numbers tending to $+\infty$. For each integer $k$, the geodesic of $\mathbf{H}^{2,n}$ between $\gamma(t_k)$ and $\tilde{\gamma}(t_k)$ is timelike, making a sequence of timelike geodesics indexed by $k$. Recall that a timelike geodesic in $\mathbf{H}^{2,n}$ is the projectivization of a plane of signature $(0, 2)$. Up to a subsequence, this sequence of planes converges to a plane that cannot be of signature $(1, 1)$ (such a plane would project to a spacelike geodesic of $\mathbf{H}^{2,n}$). Since $(\tilde{\gamma}(t_k))_{k \in \mathbb{N}}$ tends to $r$, the only possibility is that $(\gamma(t_k))_{k \in \mathbb{N}}$ converges to $r$. Similarly, $\lim\limits_{t \to -\infty} \gamma(t) = s$.

Since $p$ and $q$ are positive points, they are not in a photon of $\mathbf{Ein}^{1,n}$ containing $s$ or $r$. Hence

$$\lim_{t \to +\infty} g_t \cdot p = \lim_{t \to +\infty} g_t \cdot q = s, \text{ and } \lim_{t \to -\infty} g_t \cdot p = \lim_{t \to -\infty} g_t \cdot q = r \ .$$

Denote by $h_p$ and $h_q$ representatives of the space-horofunctions associated with $p$ and $q$, respectively (recall Definition 2.3.6). We have

$$\begin{aligned} g_t \cdot \{x \in \Sigma \mid h_p(x) \le h_p(\gamma(t)\} &= \{g_t \cdot x \in g_t \cdot \Sigma \mid h_p(x) \le h_p(\gamma(t))\} \\ &= \{y \in g_t \cdot \Sigma \mid h_p(g_t^{-1} \cdot y) \le h_p(\gamma(t))\} \\ &= \{y \in g_t \cdot \Sigma \mid h_{g_t \cdot p}(y) \le h_{g_t \cdot p}(g_t \cdot \gamma(t))\} \ . \end{aligned}$$

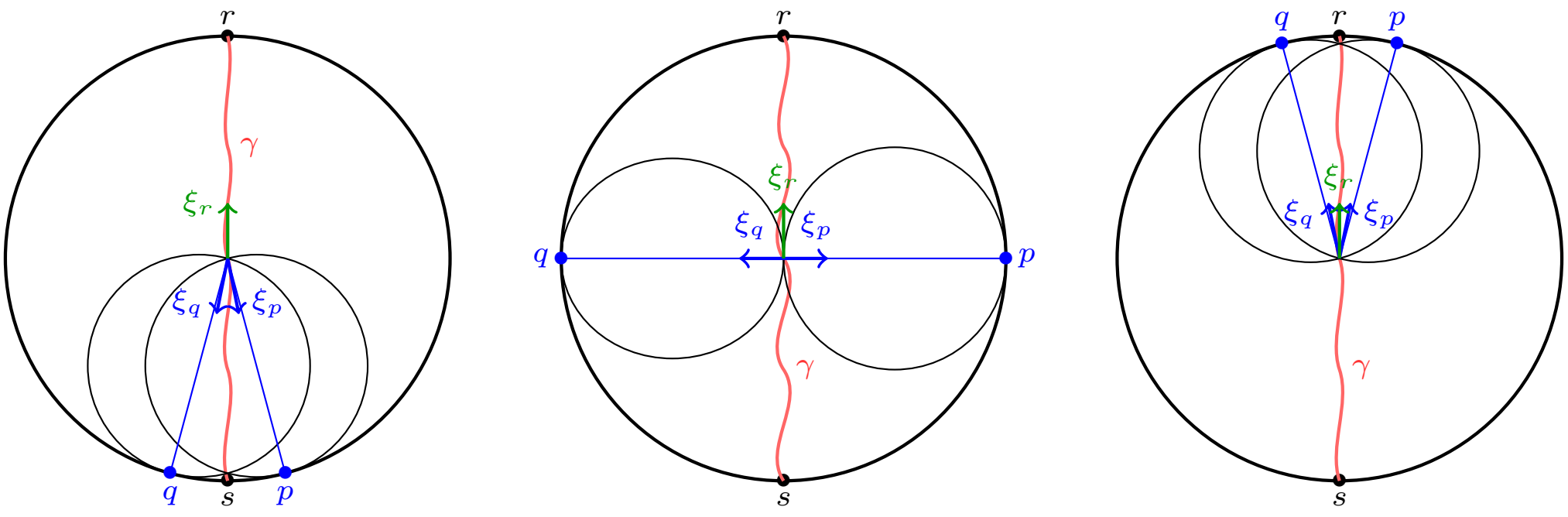

**Fig. 9.** The angle between ideal boundary points $\xi_p$ and $\xi_q$ that are regular with respect to two positive points $p$ and $q$, respectively. On the left, the case where $p$ and $q$ are "close to each other" on the side of $s$. On the right, the case where $p$ and $q$ are "close to each other" on the side of $r$. In between, there is a moment where $\xi_p$ and $\xi_q$ are opposite to each other.

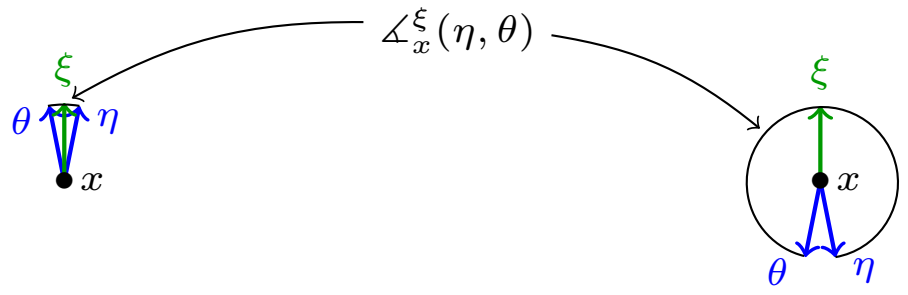

**Fig. 10.** Illustration of the measure $\measuredangle^\xi(\eta, \theta)$ at the point $x$. On the left, the measure is close to 0, on the right, the measure is close to $2\pi$.

Similarly,

$$g_t \cdot \{x \in \Sigma \mid h_q(x) \leq h_q(\gamma(t))\} = \{y \in g_t \cdot \Sigma \mid h_{g_t \cdot q}(y) \leq h_{g_t \cdot q}(g_t \cdot \gamma(t))\} .$$

When $t$ goes to $+\infty$ (respectively $-\infty$), the sequences of space-horofunctions $(h_{g_t \cdot p})_{k \in \mathbb{N}}$ and $(h_{g_t \cdot q})_{k \in \mathbb{N}}$ converge smoothly and uniformly on compact sets, and modulo an additive constant, to a representative of the space-horofunction $h_s$ (respectively $h_r$). As a consequence, the space-horocycles

$$\{y \in g_t \cdot \Sigma \mid h_{g_t \cdot p}(y) = h_{g_t \cdot p}(g_t \cdot \gamma(t))\} \text{ and } \{y \in g_t \cdot \Sigma \mid h_{g_t \cdot q}(y) = h_{g_t \cdot q}(g_t \cdot \gamma(t))\}$$

are more and more tangent at $g_t \cdot \gamma(t)$ (see the left and right drawings in Fig. 9), when $t$ goes to $+\infty$ or $-\infty$. Moreover, when $t$ goes to $+\infty$ (respectively $-\infty$), the direction $\xi_r$ points outside (respectively inside) the limit of the space-horoballs associated to $p$ and $q$. That is:

$$\limsup_{t \to +\infty} \measuredangle^{\xi_r}_{g_t \cdot \gamma(t)}(\xi_p, \xi_q) \leq \pi, \text{ and } \liminf_{t \to +\infty} \measuredangle^{\xi_r}_{g_t \cdot \gamma(t)}(\xi_p, \xi_q) \geq \pi \ ,$$

where $\measuredangle^\xi_x(\eta, \theta)$ denotes the angle, with values in $[0, 2\pi)$, between $\eta$ and $\theta$ *in the side of* $\xi$. More precisely, $\measuredangle^\xi_x(\eta, \theta)$ is the distance between $\eta_x$ and $\theta_x$ in $\mathsf{T}^1_x S \setminus \{-\xi_x\}$ if $\xi$, $\eta$ and $\theta$ are vector fields of the maximal surface $S$ (see Fig. 10).

Finally, either there is a value $t$ for which

$$\measuredangle_{g_t\cdot\gamma(t)}(\xi_p, \xi_q) = \pi \ ,$$

or there is a sign $\varepsilon$ in $\{-1, 1\}$ such that

$$\lim_{t\to\varepsilon\infty} \measuredangle_{g_t\cdot\gamma(t)}(\xi_p, \xi_q) = \pi \ .$$

That is,

$$\measuredangle(\xi_p, \xi_q) = \pi \ .$$

This is true for every pair of ideal points $(p, q)$ in $I^2$, with $p$ different from $q$, and every ideal boundary points $\xi_p$, $\xi_q$ regular with respect to $p$, $q$, hence the result. □

**Remark 4.7.3.** Let $p$ be a positive point of $\partial^s\Sigma$, where $\Sigma$ is a complete maximal surface of $\mathbf{H}^{2,n}$. We did not prove that an ideal ray regular with respect to $p$ is unique. If several exist, they may have finite Tits distance.

### *4.8. The total curvature of maximal surfaces*

We show in this section that a maximal surface is polygonal if and only if it has finite Tits perimeter.

#### *4.8.1. A formula for the total curvature of a Hadamard surface in terms of the Tits distance*

For a Hadamard surface, having finite Tits perimeter is equivalent to having finite total curvature, by a result of Ohtsuka [21] that we explain now.

Let $(\Sigma, g)$ be a Hadamard manifold of dimension 2.

We precise, given two ideal boundary points $\xi, \eta$ of $\Sigma$, the angular domain delimited by $\xi, \eta$ and a point $x$ of $\Sigma$. Given $x$ in $M$ denote by $c_\xi$ and $c_\eta$ the geodesic rays starting at $x$ with direction $\xi$ and $\eta$, respectively. One might want to define the angular domain determined by $\xi, \eta$ and $x$ by

$$\{\sigma_t([0,1]) \mid \sigma_t : [0,1] \to \Sigma \text{ geodesic ray from } c_\xi(t) \text{ to } c_\eta(t),\ t \in [0, +\infty)\} \ .$$

This is not totally satisfying. Indeed, by Lemma 4.6.9, if $\mathrm{Td}(\xi, \eta) > \pi$, the above set is an infinite triangle with two ideal points. But if $\mathrm{Td}(\xi, \eta) < \pi$, the above set is what we want. Also, we will cut the domain in several subdomains and we need to make a consistent choice, unambiguous when cuting the domain. Hence we make the following definition, following [21].

**Definition 4.8.1.** Let $\xi, \eta$ be two ideal boundary points of $\Sigma$.

If $\mathrm{Td}(\xi,\eta) < +\infty$, the *angular domain* delimited by $\xi, \eta$ and $x$ is defined as being the connected component of $\Sigma \setminus (c_\xi([0,+\infty]) \cup c_\eta([0,+\infty]))$ containing the ray from $x$ to a midpoint between $\xi$ and $\eta$ in $(\Sigma(\infty), \mathrm{Td})$.

If $\mathrm{Td}(\xi,\eta) = +\infty$, the *angular domain* delimited by $\xi, \eta$ and $x$ is any connected component of $\Sigma \setminus (c_\xi([0,+\infty]) \cup c_\eta([0,+\infty]))$.

**Lemma 4.8.2** *(Theorem in [21]). Let $(\Sigma, g)$ be a Hadamard manifold of dimension 2. Take $\xi$ and $\eta$ two ideal points, $x$ in $\Sigma$ and denote by $A$ the angular domain determined by $\xi, \eta$ and $x$. Denote by $K_y$ the Gaussian curvature of $\Sigma$ at the point $y$ of $\Sigma$.*

*If* $\mathrm{Td}(\xi,\eta) < +\infty$, *then*

$$\int_A K_y \mathrm{d}v_g(y) = \measuredangle_x(\xi,\eta) - \mathrm{Td}(\xi,\eta) \ .$$

*If* $\mathrm{Td}(\xi,\eta) = +\infty$, *then*

$$\int_A K_y \mathrm{d}v_g(y) = -\infty \ .$$

The lemma is a generalized version of the Gauss Bonnet formula. We recall the Gauss Bonnet formula in the case of a geodesic triangle. Let $(\Sigma, g)$ be a Riemannian 2-dimensional manifold. Let $c_1$, $c_2$ and $c_3$ be three geodesics defined in $[0,1]$ such that $c_i(1)$ equals $c_{i+1}(0)$ for $i$ in $\{1,2,3\}$ and the concatenation of the three curves bounds a region homeomorphic to a disc. Denote by $T$ the corresponding geodesic triangle (the compact region bounded by the three segments). Denote by $\theta_1$, $\theta_2$ and $\theta_3$ the interior angles of the triangle. Then we have the *Gauss-Bonnet formula*:

$$\int_T K_y \mathrm{d}v_g(y) = -\pi + \sum_{i=1}^{3} \theta_i \ .$$

The idea of the proof of Lemma 4.8.2 is to take geodesic rays $c_\xi$ and $c_\eta$ starting at $x$ and to consider the triangle with vertices $x$, $c_\xi(s)$ and $c_\eta(t)$. The Gauss-Bonnet formula gives us the total curvature inside this triangle, depending on the angles at each vertex. Now making $s$ and $t$ go to infinity gives us the result, using properties of the angle distance, in the case $\measuredangle(\xi,\eta)$ is smaller than $\pi$. If $\mathrm{Td}(\xi,\eta)$ is bigger than $\pi$, then by 4.6.10 we can find several ideal boundary points $\xi_1, \ldots, \xi_k$ such that $\mathrm{Td}(\xi_i, \xi_{i+1})$ is smaller than $\pi$ for every $i$ and this yields the result for $\mathrm{Td} < +\infty$.

When $\mathrm{Td} = \infty$, we can divide the segment between $\xi$ and $\eta$ in $M(\infty)$ in an arbitrary large number of segments of Tits length bigger than $\pi$. We conclude that the total curvature cannot be bounded by below.

This beautiful formula implies a strong consequence about the total curvature of Riemannian surfaces with finite diameter Tits boundary.

**Theorem 4.8.3.** *Let $(\Sigma, g)$ be a Hadamard surface. The Tits distance in the ideal boundary $\Sigma(\infty)$ is that of a circle of perimeter $P$ in $[2\pi, +\infty)$ if and only if the total curvature of $(\Sigma, g)$ satisfies*

$$\int_\Sigma K = 2\pi - P \ .$$

**Proof.** Suppose that the Tits distance on the ideal boundary is that of a circle of perimeter $P$ in $[1, +\infty[$. For $k = \lfloor P/\pi \rfloor + 1$, there are $k$ points $\xi_1, \ldots, \xi_k$ of $\Sigma(\infty)$ such that

- The distance $\mathrm{Td}(\xi_i, \xi_{i+1})$ is less than $\pi$, for every $i$,
- The distance $\mathrm{Td}(\xi_i, \xi_j)$ is bigger than $\pi$, for every $|i-j| > 1$,
- The sum of the distances $\mathrm{Td}(\xi_i, \xi_{i+1})$, for $i$ through $\{1, \ldots, k\}$ equals $P$.

Here the indices are modulo $k$. We can apply the formula of Ohtsuka, Lemma 4.8.2, for a point $x$ in $\Sigma$ and the domains bounded by $\xi_i$, $\xi_{i+1}$ for $i$ through $\{1, \ldots, n\}$. Now, the sum of the angles $\measuredangle_x(\xi_i, \xi_{i+1})$ at $x$ equals $2\pi$, and the sum of the distances $\measuredangle(\xi_i, \xi_{i+1}) = \mathrm{Td}(\xi_i, \xi_{i+1})$ equals $P$, hence the formula for the total curvature.

Now if there are two points $\xi$ and $\eta$ in the ideal boundary $\Sigma(\infty)$ such that $\mathrm{Td}(\xi, \eta) = +\infty$, then the formula of Ohtsuka, Lemma 4.8.2, says that the total curvature of $\Sigma$ is $-\infty$. □

**Remark 4.8.4.** The proof is done by finding a $\pi$-net in $(\Sigma(\infty), \mathrm{Td})$, describing the Tits distance completely. Nets will be used in Section 5.5 to compare different Tits boundaries in the same surface.

### *4.8.2. Application to maximal surfaces*

**Theorem 4.8.5.** *Let $\Sigma$ be a maximal surface in $\mathbf{H}^{2,n}$. The Tits distance on $\Sigma(\infty)$ has finite perimeter if and only if $\Sigma$ is a polygonal surface. A polygonal surface with $N+4$ vertices has Tits perimeter equal to $(N+4)\pi/2$, hence total curvature equal to $-\frac{\pi}{2}N$.*

**Proof.** Suppose that $\Lambda$, the space boundary of $\Sigma$, is a lightlike polygon with $N+4$ vertices. As shown in Theorem 4.5.2 and Theorem 4.7.1, they are $N+4$ distinct ideal points $\xi_1, \ldots, \xi_{N+4}$, corresponding to vertices, such that $\mathrm{Td}(\xi_i, \xi_{i+1}) = \pi/2$ for each $i$. Then the Tits distance $(\Sigma(\infty), \mathrm{Td})$ is that of a circle of perimeter $(N+4)\pi/2$. Using the formula in Theorem 4.8.3 we obtain that

$$\int_\Sigma K = 2\pi - (N+4)\frac{\pi}{2} = -\frac{\pi}{2}N \ .$$

Suppose now that $\Lambda$ is not a lightlike polygon. If $\Lambda$ contains an infinite number of adjacent segment of photons, the arguments used in the polygonal case show that the total curvature of $\Sigma$ is $-\infty$. Suppose that $\Lambda$ does not contain an infinite number of adjacent photons. Then $\Lambda$ contains a positive open non-empty interval, denoted $I$.

By Theorem 4.7.2, for each point $p$ of $I$ there is an associated ideal point $\xi$, called regular with respect to $p$, such that the Tits distance between each pair of these regular ideal point is infinite. Hence, by Theorem 4.8.3 we have

$$\int_{\Sigma} K = -\infty \ ,$$

and this finishes the proof. □

By [3], a simply connected non-compact Riemannian surface with finite total curvature must be of parabolic type. Hence we obtain the following:

**Corollary 4.8.6.** *A polygonal surface is of parabolic type.*

## 5. The quartic boundary of maximal surfaces

We study now the quartic boundary of a maximal surface, and its links with the space boundary and the ideal boundary.

In Section 5.1, we define the quartic differential and its associated flat metric. In Section 5.2, we show that an asymptotically flat maximal surface uniformizes with the complex plane and has polynomial quartic differential. Section 5.3 is devoted to a characterization of compact Tits boundaries, allowing us to compare two Hadamard distances on the same surface. In Section 5.4, we explain how a maximal surface with the quartic metric can be thought of as a Hadamard surface, and explain how the Tits distance construction extends to this setting. We use this results in Section 5.5 to show that the quartic boundary and the ideal boundary of a polygonal maximal surface are isometric.

### *5.1. The quartic differential associated to maximal surfaces*

We denote now by $g$ the induced metric on the tangent space of a spacelike surface $\Sigma$ and by $g_{\mathsf{N}}$ the induced metric on its normal bundle $\mathsf{N}\Sigma$. If $\Sigma$ is a maximal surface of $\mathbf{H}^{2,n}$, there is a holomorphic object on $(\Sigma, [g])$ one can build using the second fundamental form: the *quartic differential*. This differential can be turned into a flat metric with conical singularities, that will be very close to the induced metric in the polygonal case.

#### *5.1.1. Definition of the quartic differential*

The second fundamental form $\mathbb{I}$ of $\Sigma$ is a symmetric 2 tensor with values in the normal bundle $\mathsf{N}\Sigma$. By complexifying, we obtain a 2-tensor

$$\mathbb{I}^{\mathbb{C}} : \mathsf{T}^{\mathbb{C}}\Sigma \times \mathsf{T}^{\mathbb{C}}\Sigma \to \mathsf{N}^{\mathbb{C}}\Sigma \ ,$$

that we can decompose under types using the complex structure coming from $g$ ($T^{\mathbb{C}}\Sigma = \mathsf{T}^{0,1}\Sigma \oplus \mathsf{T}^{1,0}\Sigma$). The $(1,1)$ part of $\mathbb{I}^{\mathbb{C}}$ is null because it equals the $g$-trace of $\mathbb{I}$ and $\Sigma$ is maximal. The bundle

$$E = (\mathsf{T}^{1,0}\Sigma)^{\vee} \otimes (\mathsf{T}^{1,0}\Sigma)^{\vee} \otimes \mathsf{N}^{\mathbb{C}}\Sigma$$

is endowed with the connection $\nabla^E$ coming from the complexification of the Levi-Civita connection of $\Sigma$ on the tangent bundle and from the normal connection of $\Sigma$ on the normal bundle. The $(0,1)$ part $\bar{\partial}^E$ of the associated exterior derivative $\mathrm{d}_{\nabla^E}$ is a Dolbeault operator. The $(2,0)$ part of $\mathbb{I}^{\mathbb{C}}$ is $\bar{\partial}^E$ holomorphic (see [20, Lemma 5.6] for a proof).

**Definition 5.1.1.** The *holomorphic quartic differential* associated to $(\Sigma, [g])$ is the holomorphic section of the fourth power of the cotangent bundle $\mathcal{K}^4$ defined by

$$\varphi_4 := g_{\mathsf{N}}^{\mathbb{C}} \left(2\mathbb{I}^{2,0}, 2\mathbb{I}^{2,0}\right) \ .$$

*5.1.2. Expression in a chart*

First, denote by $J$ the complex structure associated with the induced metric $g$. Let $X, Y$ be complex vector fields on $\Sigma$. They can be decomposed in $(1,0)$ and $(0,1)$ type and we have

$$\begin{aligned} 2\mathbb{I}^{2,0}(X,Y) &= 2\mathbb{I}^{\mathbb{C}}(X^{1,0}, Y^{1,0}) \\ &= \frac{1}{2}\mathbb{I}^{\mathbb{C}}(X - iJX, Y - iJY) \\ &= \mathbb{I}(X,Y) - i\mathbb{I}(JX,Y) \ , \end{aligned}$$

where for the last equality we used the fact that $\mathbb{I}(J\cdot, J\cdot) = -\mathbb{I}(\cdot,\cdot)$.

Secondly, we compute $\varphi_4$ in a holomorphic chart. Take $z = x + iy$ a holomorphic coordinate centered at a point $p$ of $\Sigma$. The metric $g$ reads in this chart

$$g_z = \exp(2f(z))(\mathrm{d}x^2 + \mathrm{d}y^2) \ ,$$

with $f$ a smooth function on $\mathbb{C}$. We now have

$$2\mathbb{I}^{2,0} = 2\mathbb{I}^{2,0}(\partial_z, \partial_z)\mathrm{d}z^2 \ ,$$

and

$$2\mathbb{I}^{2,0}(\partial_z, \partial_z) = \mathbb{I}(\partial_x, \partial_x) - i\mathbb{I}(\partial_x, \partial_y) \ .$$

Hence, taking the complex extension of $g_{\mathsf{N}}$, we obtain in the chart $z$

$$\varphi_4 = \varphi_4(z)\mathrm{d}z^4 \ ,$$

with

$$\varphi_4(z) = g_{\mathsf{N}}(\mathbb{I}(\partial_x, \partial_x), \mathbb{I}(\partial_x, \partial_x)) - g_{\mathsf{N}}(\mathbb{I}(\partial_x, \partial_y), \mathbb{I}(\partial_x, \partial_y)) - 2i g_{\mathsf{N}}(\mathbb{I}(\partial_x, \partial_x), \mathbb{I}(\partial_x, \partial_y)) \ . \tag{10}$$

Denote by $f$ the homomorphism $\mathbb{I}(\partial_x, \cdot) : (\mathsf{T}\Sigma, g) \to (\mathsf{N}\Sigma, -g_{\mathsf{N}})$. The last equation says that at a point $z$ where $\varphi_4(z) = 0$ the map $f$ is conformal:

$$\begin{cases} g_{\mathsf{N}}(f(\partial_x), f(\partial_x)) = g_{\mathsf{N}}(f(\partial_y), f(\partial_y)) \ , \\ g_{\mathsf{N}}(f(\partial_x), f(\partial_y)) = 0 \ . \end{cases}$$

*5.1.3. The metric associated to a quartic differential*

A holomorphic quartic differential $\varphi$ on $\mathbb{C}$ is written in the $z$-parameter

$$\varphi = \varphi(z)\mathrm{d}z^4 \ ,$$

where $z \mapsto \varphi(z)$ is a holomorphic function on $\mathbb{C}$.

It defines a smooth flat metric $g := |\varphi|^{1/2}$ on $\mathbb{C} \setminus \{ \text{ zeros of } \varphi\}$. More precisely, at every point $z$:

$$g_z = |\varphi(z)|^{\frac{1}{2}} |\mathrm{d}z|^2 \ .$$

This formula is invariant by coordinate change.

In the neighborhood of a zero, it is possible ([28, Chapter III, Paragraph 6]) to find a local parameter $w$ such that the quartic differential $\varphi$ is written in this coordinate

$$\varphi = \varphi(w)\mathrm{d}w^4 = w^k \mathrm{d}w^4 \ ,$$

where $k$ is the order of the zero. This parameter is unique up to the multiplication by a $(k+4)$-th root of unity. In this parameter, the metric $g$ is written

$$g_w = |w|^{\frac{k}{2}} |\mathrm{d}w|^2 \ ,$$

and is isometric to a Euclidean cone of angle $\theta = \frac{\pi}{2}k + 2\pi$ ([32, Proposition 1]). At a nonzero point $z$, the metric $g$ has vanishing curvature. This is the reason why we talk about *flat metrics with conical singularities*.

**Definition 5.1.2.** Let $\Sigma$ be a maximal surface in $\mathbf{H}^{2,n}$ and $\varphi_4$ be the associated quartic differential. We define

$$g_4 := |\varphi_4|^{1/2} \ ,$$

and call $g_4$ the *quartic metric* of $\Sigma$.

### 5.2. Asymptotically flat maximal surface

We know that polygonal surfaces are of parabolic type. The aim of this section is to prove that asymptotically flat maximal surfaces (not supposed to be polygonal) are parabolic and have polynomial quartic differential.

**Proposition 5.2.1.** *Let $\Sigma$ be an asymptotically flat maximal surface of $\mathbf{H}^{2,n}$ with induced metric $g$ and quartic metric $g_4$. For every $\varepsilon > 0$ there is a compact set $K$ in $\Sigma$ such that $g$ and $g_4$ are $(1+\varepsilon)$-biLipschitz outside $K$.*

**Proof.** Let $\Sigma$ be an asymptotically flat maximal surface. Let $(x_k)_{k\in\mathbb{N}}$ be a diverging sequence of points of $\Sigma$. By Theorem 2.2.7, there is a sequence $(\gamma_k)_{k\in\mathbb{N}}$ of isometries of $\mathbf{H}^{2,n}$ such that the sequence $((g_n, \cdot\Sigma, g_n \cdot x_n))_{k\in\mathbb{N}}$ subconverges to a pointed maximal surface $(S, y)$. Now

$$K^S(y) = \lim_{k\to+\infty} K^\Sigma(x_k) = 0 \ ,$$

hence, by [19, Proposition 5.3], $S$ is a Barbot surface.

In a Barbot surface, the quartic differential equals $\mathrm{d}z^4$, as one can see using Equation (10) and the computations in Proposition 4.1.2. Hence, in this case, the quartic metric equals the flat induced metric.

Recall that the convergence of pointed maximal surfaces is uniform and smooth on compact sets (see Section 2.2.2). As a consequence, for every positive number $\varepsilon$, there is an integer $K$ such that for every integer $k$ bigger than $K$ the metrics $g$ and $g_4$ are $(1+\varepsilon)$-biLipschitz at $x_k$.

By contradiction, suppose that there is a positive number $\varepsilon$ such that $\Sigma$ has no compact subset outside of which $g$ and $g_4$ are $(1+\varepsilon)$-biLipschitz. In this case we can find a diverging sequence of points $(x_k)_{k\in\mathbb{N}}$ of $\Sigma$ such that $g$ and $g_4$ are not $(1+\varepsilon)$-biLipschitz along this sequence. This contradicts asymptotic flatness by the previous arguments. □

We now follow the arguments of [10, Section 7] to prove that an asymptotically flat maximal surface is of parabolic type and has polynomial quartic differential. We put the proof in our context and notations.

**Lemma 5.2.2** *(Lemma 7.3 in [10]). Let $\Sigma$ be an asymptotically flat maximal surface in $\mathbf{H}^{2,n}$ with quartic differential $\varphi_4$ and quartic metric $g_4$. Then $\varphi_4$ has a finite number of zeros, and $g_4$ is complete.*

**Proof.** The previous lemma shows that the zero set of $\varphi_4$ lies in a compact set. It is discrete because $\varphi_4 \neq 0$, hence finite.

Outside a compact $K$, the induced metric $g$ and the quartic metric $g_4$ are bi-Lipschitz. If $c$ is a $g_4$-geodesic, it is a $g$-quasi geodesic outside $K$, hence is complete. □

**Lemma 5.2.3** *(Lemma 7.5 in [10]). Let $\Sigma$ be an asymptotically flat maximal surface in $\mathbf{H}^{2,n}$ with induced metric $g$. The Riemann surface $(\Sigma,[g])$ is of parabolic type.*

**Proof.** If not, then $(\Sigma,[g])$ is bi-holomorphic to the unit disc. We write $\varphi_4$ in the $z$-coordinate:

$$\varphi_4 = P(z)\exp(f(z))\mathrm{d}z^4,$$

with $P$ a polynomial and $f$ holomorphic since $\varphi_4$ has a finite number of zeroes.

The metric $|\exp(f)|^{1/2}|\mathrm{d}z|^2$ is a flat metric, conformal to the hyperbolic metric of the unit disc. Hence it cannot be complete because the unit disc is not bi-holomorphic to the complex plane. As a consequence, there is a divergent path $\gamma$ of finite $|\exp(f)|^{1/2}|\mathrm{d}z|^2$-length.

Now the $g_4$-length of $\gamma$ is also finite because $P$ is bounded on the unit disc. Hence $g_4$ is not complete, contradicting the previous lemma. □

**Theorem 5.2.4** *(Lemma 9.6 in [24]). Let $\Sigma$ be an asymptotically flat maximal surface in $\mathbf{H}^{2,n}$ with quartic differential $\varphi_4$. Then $\varphi_4$ is polynomial.*

**Proof.** We identify $(\Sigma,[g])$ with the complex plane with $z$-coordinate. So $\varphi_4 = \varphi_4(z)\mathrm{d}z^4$. Since $\varphi_4$ has a finite number of zeros,

$$\varphi_4 = P(z)\exp(f(z))\mathrm{d}z^4 \ ,$$

with $f$ a holomorphic function on $\mathbb{C}$. We will prove that $f$ is constant. If $P$ is of degree $N$ and $M$ is an integer bigger than $N/4$, then

$$|\varphi_4(z)|^{1/4} \underset{|z|\to+\infty}{=} O(|z^M\exp(f(z))|) \ .$$

Since $g_4$ is complete, for every divergent path $\gamma$, the $g_4$-length $\int_\gamma |\varphi_4(z)|^{1/4}|\mathrm{d}z|$ diverges, hence also $\int_\gamma |z^M\exp(f(z))||\mathrm{d}z|$. Now take

$$h(z) := \int\limits_{[0,z]} \zeta^M\exp(f(\zeta))\mathrm{d}\zeta \ .$$

The function $h$ has a zero of order $M+1$ in $z=0$, hence

$$h(z) = z^{M+1}(b_0 + b_1 z + \ldots) \ ,$$

with $b_0 = a_0/(M+1) \neq 0$. Define another parameter $w$ by a single valued branch of $h^{1/(M+1)}$:

$$w(z) := h^{1/(M+1)}(z) = z(b_0 + b_1 z + \ldots)^{1/(M+1)} \ .$$

The parameter $w$ is well defined in a neighborhood of 0, and since $w'(0) = b_0^{1/(M+1)} \neq 0$ it has an inverse $z = \psi(w)$ where $\psi$ is defined in a neighborhood of zero. We show now that $\psi$ is defined globally. Take a path $\gamma(t) = tw_0$ where $w_0$ is in $\mathbb{C}^*$ and $t$ in $[0, 1)$. The composition $\psi \circ \gamma$ is a path in the $z$-plane. Since $\mathrm{d}(w^{M+1})(z) = \mathrm{d}h(z) = z^M \exp(f(z))$, we have

$$\int_\gamma |z^M \exp(f(z))||\mathrm{d}z| = \int_\gamma |\mathrm{d}(w^{M+1})| = |w_0|^{M+1} .$$

Since $\psi \circ \gamma$ has finite $|z^M \exp(f(z))||\mathrm{d}z|$-length, hence finite $g_4$-length, it is not divergent in the $z$-plane. So $\psi \circ \gamma(t)$ has some adherence value when $t \to 1$ and since $(\psi \circ \gamma)'(1) = \psi'(w_0) \neq 0$, the path $\psi \circ \gamma$ can be continued in $t = 1$. Hence, $\psi$ is an automorphism of $\mathbb{C}$, so is $w$ as a function of $z$. The function $w(z) = h^{1/(M+1)}(z)$ being linear implies that $h$ is a polynomial in $z$. Hence $h'(z) = z^M \exp(f(z))$ is polynomial and $f$ is constant. □

### *5.3. Compact Tits boundaries and nets*

We study in this section the Tits boundary of Hadamard surfaces, and more precisely when the Tits distance is compact.

#### *5.3.1. Characterization of compact Tits boundaries*

**Theorem 5.3.1.** *Suppose that $M$ is a Hadamard surface or a flat metric with conical singularities that is a Hadamard space. Then the following conditions are equivalent.*

*(i) The metric space $(M(\infty), \mathrm{Td})$ is compact,*
*(ii) The metric space $(M(\infty), \mathrm{Td})$ is isometric to a circle of perimeter $P \geq 2\pi$,*
*(iii) The metric space $(M(\infty), l)$ admits a finite 1-net.*

*Moreover, in the case $(M(\infty), \mathrm{Td})$ is isometric to a circle of perimeter $P$, we have for every point $o$ in $M$*

$$P = \lim_{r \to +\infty} \frac{P_o(r)}{r} ,$$

*where $P_o(r)$ denotes the perimeter of the circle of radius $r$ in $(M, d)$ centered at $o$.*

A *$\varepsilon$-net* in a metric space $(E, d)$, following [13, Definition 2.14], is a subset $N$ of $E$ such that $d(x, N) < \varepsilon$ for every $x$ in $E$. We can replace 1-net by $\delta$-net with any $\delta < 2$ in this Theorem, but the value 1 is more convenient.

**Remark 5.3.2.** The last theorem does not hold in higher dimensions. For example, the Tits boundary of a higher rank symmetric space of noncompact type always admits a finite $\sqrt{2}$-net, but is not compact. Indeed, length spaces of dimension 1 are determined by their length and nothing else, but as soon as the dimension is bigger than 1, the situation is much more complex.

**Proof.** $(i) \Rightarrow (ii)$. Suppose that $(M(\infty), \mathrm{Td})$ is compact. Then the natural map $\varphi : (\mathsf{T}^1_x\Sigma, \measuredangle_x) \to (\Sigma(\infty), \mathrm{Td})$ is a homeomorphism and the Tits boundary is homeomorphic to a circle. Now $\mathrm{Td} : \Sigma(\infty)^2 \to \mathbb{R}$ is continuous, hence bounded, and $(M(\infty), \mathrm{Td})$ has finite diameter $D \geq \pi$. As a consequence, since Td is an inner distance, $(M(\infty), \mathrm{Td})$ has finite perimeter $P = 2D \geq 2\pi$.

$(ii) \Rightarrow (iii)$ is immediate.

In order to prove $(iii) \Rightarrow (i)$, we recall the following fact. By Proposition 4.6.4 and Lemma 4.6.8, if the $l$ distance between two points $\xi, \eta$ of $M(\infty)$ is strictly smaller than 2, then

$$\mathrm{Td}(\xi, \eta) = 2 \arcsin\left(\frac{l(\xi, \eta)}{2}\right) < \pi$$

and by Lemma 4.6.10 the Tits geodesic between $\xi$ and $\eta$ in $M(\infty)$ is unique.

Suppose that $(M(\infty), l)$ admits a finite 1-net $N = \{\xi_1, \ldots, \xi_K\}$. Then $N$ is a $\frac{\pi}{3}$-net in $(M(\infty), \mathrm{Td})$. This implies that $(M(\infty), \mathrm{Td})$ is compact. Indeed take a sequence of points $(\eta_k)_{k \in \mathbb{N}}$ in $M(\infty)$. There is an integer $i$ between 1 and $K$ such that an infinite number of $\eta_k$ is at Tits distance less than $\frac{\pi}{3}$. This provides us a subsequence $(\eta_{\psi(k)})_{k \in \mathbb{N}}$ included in $\{\eta \in M(\infty) \mid \mathrm{Td}(\eta, \xi_i) \leq \frac{\pi}{3}\}$ and this last set is homeomorphic to a segment of length $\frac{2\pi}{3}$, yielding the result.

In order to prove the last assertion, it suffices to use the following facts:

- If $(M(\infty), \mathrm{Td})$ is compact, the sequence of metric spaces $(\mathbf{S}_o(r), d/r)$ converges, in the Gromov-Hausdorff sense, to $(M(\infty), l)$ when $r$ goes to $+\infty$. Here $o$ is a point of $M$ and $\mathbf{S}_o(r)$ is the circle centered at $o$ and of radius $r$,
- a sequence of length spaces converges to a length space ([13, Proposition 3.8]).

Hence the length space obtained from $(\mathbf{S}_o(r), d/r)$, that is a circle of perimeter $P_o(r)/r$, converges to $(M(\infty), \mathrm{Td})$. □

*5.3.2. Comparison of distinct Tits boundaries*

In this paragraph, we compare two Hadamard metrics on the same topological surface. We show that if the two metrics are uniformly asymptotic, then they give the same Tits distances in the ideal boundary.

**Corollary 5.3.3.** *Let $\Sigma$ be diffeomorphic to $\mathbb{R}^2$, and $d_1, d_2$ be two Hadamard distances in $\Sigma$. Suppose that for every positive number $\varepsilon$ the identity map is a $(1+\varepsilon, b)$-quasi isometry between $(\Sigma, d_1)$ and $(\Sigma, d_2)$, for a nonnegative number $b$. Then $(\Sigma(\infty), \mathrm{Td}_1)$ is compact* *if and only if* $(\Sigma(\infty), \mathrm{Td}_2)$ *is compact. Moreover, if* $(\Sigma(\infty), \mathrm{Td}_1)$ *and* $(\Sigma(\infty), \mathrm{Td}_2)$ *are compact, they are isometric.*

We fix some notations. Let us fix $o$ a point in $\Sigma$. We denote by $\mathbf{S}_i(r)$ the circle of radius $r$ centered at $o$ in $(\Sigma, d_i)$. We define a map to compare the circles $\mathbf{S}_1(r)$ and $\mathbf{S}_2(r)$ as follows.

If $x$ is a point of $\mathbf{S}_1(r)$, denote by $c_{ox}$ the $d_2$-geodesic of $\Sigma$ from $o$ to $x$. This geodesic intersects the circle $\mathbf{S}_2(r)$ at a unique point $y$, and we define $f_r(x) = y$. More precisely,

$$\begin{aligned} f_r : \mathbf{S}_1(r) &\to \mathbf{S}_2(r) \\ x &\mapsto c_{ox}\left(\frac{r}{d_2(o,x)}\right) . \end{aligned}$$

For topological reasons, the map $f_r$ is surjective.

**Lemma 5.3.4.** *For every positive number* $\varepsilon$*, there is a positive number* $r_0$ *and a function* $\delta : [r_0, +\infty) \to \mathbb{R}_+$ *satisfying* $\lim\limits_{r\to+\infty} \delta(r)/r = 0$ *such that for every* $r$ *bigger than or equal to* $r_0$ *the map* $f_r$ *is a* $(1+\varepsilon, \delta(r))$*-quasi isometry from* $(\mathbf{S}_1(r), d_1)$ *to* $(\mathbf{S}_2(r), d_2)$*.*

**Proof.** Fix a positive number $\varepsilon$. By assumption, there exists a positive number $b$ such that the distances $d_1$ and $d_2$ are $(1+\varepsilon, b)$-quasi isometric through the identity map. Take a positive number $r_0 > (1+\varepsilon)b$, and $r \geq r_0$. Let $x$ be a point in $\mathbf{S}_1(r)$. Since $d_1(o,x) = r$, we have

$$(1+\varepsilon)^{-1}r - b \leq d_2(o,x) \leq (1+\varepsilon)r + b .$$

This implies

$$d_2(f_r(x), x) = \frac{|r - d_2(o,x)|}{d_2(o,x)} \leq \frac{\varepsilon r + b}{(1+\varepsilon)^{-1}r - b} . \tag{11}$$

Denote by $\delta(r)$ the number $\delta(r) := b + 2\dfrac{\varepsilon r + b}{(1+\varepsilon)^{-1}r - b}$. By triangular inequality, for every $x, y$ in $\mathbf{S}_1(r)$,

$$\begin{aligned} d_2(x,y) - d_2(f_r(x), x) - d_2(f_r(y), y) &\leq d_2(f_r(x), f_r(y)) \\ &\leq d_2(f_r(x), x) + d_2(x,y) + d_2(y, f_r(y)) . \end{aligned}$$

Hence, using the quasi isometric inequality and Equation (11), we obtain

$$(1+\varepsilon)^{-1}d_1(x,y) - \delta(r) \leq d_2(f_r(x), f_r(y)) \leq (1+\varepsilon)d_1(x,y) + \delta(r) ,$$

for every $x, y$ in $\mathbf{S}_1(r)$. Hence $f_r$ is a $((1+\varepsilon), \delta(r))$-quasi isometry from $(\mathbf{S}_1(r), d_1)$ to $(\mathbf{S}_2(r), d_2)$, and $\frac{\delta(r)}{r}$ tends to 0 when $r$ goes to $+\infty$. $\square$

We are now ready to prove the Corollary.

**Proof of Corollary 5.3.3.** Suppose that $(\Sigma(\infty), \mathrm{Td}_1)$ is compact. Then the sequence of metric spaces $(\mathbf{S}_1(r), d_1/r)$ converges to $(\Sigma(\infty), l_1)$. Hence there is an integer $K$ and a positive number $r_1$ such that for every $r > r_1$, the metric space $(\mathbf{S}_1(r), d_1/r)$ admits a 1-net with $K$ elements. We can now map this net to $(\mathbf{S}_2(r), d_2/r)$ via the map $f_r$. We obtain that for every $r > \max(r_0, r_1)$, the metric space $(\mathbf{S}_2(r), d_2/r)$ admits a $(1+\varepsilon(r))$-net $\{\xi_1(r), \ldots, \xi_K(r)\}$, where $\varepsilon(r)$ tends to 0 when $r$ goes to $+\infty$.

Now, up to extract, and making $r$ go to $+\infty$, we obtain a 1-net $\{\xi_1, \ldots, \xi_K\}$ in $(\Sigma(\infty), l_2)$. By Theorem 5.3.1 the metric space $(\Sigma(\infty), \mathrm{Td}_2)$ is compact.

Moreover, if we suppose that the two Tits boundaries $(\Sigma(\infty), \mathrm{Td}_1)$ and $(\Sigma(\infty), \mathrm{Td}_2)$ are compact, then looking at nets equidistributed along the circles $(\mathbf{S}_1(r), d_1/r)$ we can show that for every positive number $\varepsilon$ the perimeters $P_1$ and $P_2$ of $(\Sigma(\infty), \mathrm{Td}_1)$ and $(\Sigma(\infty), \mathrm{Td}_2)$ respectively, are related by

$$(1+\varepsilon)^{-1} P_1 \leq P_2 \leq (1+\varepsilon) P_1 \ .$$

Hence $(\Sigma(\infty), \mathrm{Td}_1)$ and $(\Sigma(\infty), \mathrm{Td}_2)$ are isometric, being isometric to a circle of same perimeter. □

**Remark 5.3.5.** Corollary 5.3.3 gives a partial converse to [22, Theorem 2], where Ohtsuka shows that given two Hadamard manifolds with isometric Tits boundaries and same growth rate, they are $(1+\varepsilon, b_\varepsilon)$-quasi isometries between these two manifolds for every positive number $\varepsilon$.

## *5.4. The quartic boundary*

Let $\Sigma$ be a maximal surface with quartic differential $\varphi_4$ and associated quartic metric $g_4$. The metric $g_4$ is a flat metric with conical singularities, and these singularities are in fact negatively curved as we explain now.

### *5.4.1. Curvature near a zero*

Let $z$ be a coordinate centered at a zero of $\varphi_4$ where $\varphi_4$ is written $z^k \mathrm{d}z^4$. Then, as we said in the previous paragraph, $g_4$ has a conical singularity of angle $\theta = \frac{\pi}{2}k + 2\pi$ at the point $z = 0$.

The angle of the cone being bigger than $2\pi$, the surface behaves locally like a non-positively curved manifold.

As a consequence the distance function associated with the metric $g_4$ is Hadamard and the theory of Hadamard spaces applies (see [1] for a reference on Hadamard spaces).

**Remark 5.4.1.** Flat surfaces with conical singularities are very close to Hadamard surfaces. They are manifolds, limits of Hadamard surfaces, and the metric is singular only in a discrete set (finite in the particular case of polynomial quartic differentials).

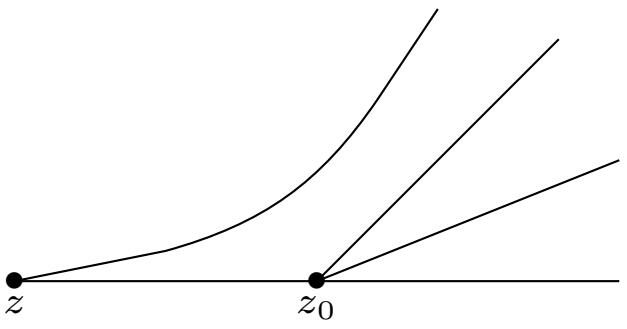

**Fig. 11.** Picture of several non-equivalent rays starting from the same point $z$ and being equal between $z$ and a zero $z_0$ of $g_4$.

*5.4.2. Total curvature*

Suppose now that $\varphi_4$ is a polynomial quartic differential, of degree $N$. Then the total curvature of $(\Sigma, g_4)$ equals

$$\sum_{z \text{ zero}} -\frac{\pi}{2}\mathrm{ord}(z) = -\frac{\pi}{2}N \ .$$

*5.4.3. Ideal boundary*

We define the ideal boundary of $(\mathbb{C}, g_4)$ in the same way as for Hadamard manifolds.

**Definition 5.4.2.** The ideal boundary of $(\mathbb{C}, g_4)$ is the set of equivalence classes of rays parametrized by arclength where two rays are equivalent if they stay at bounded distance.

The ideal boundary of $(\mathbb{C}, g_4)$ has the topological structure of a circle, but this does not come from an identification with the unit tangent space at a point $z \in \mathbb{C}$. Indeed, there can be several non-equivalent rays starting from the same point and being equal along an interval $[0, t]$, because of the zeros of the metric $g_4$ (Fig. 11).

*5.4.4. The l-distance*

In the same way as for Hadamard surfaces, the $l$-distance on the ideal boundary is defined by, given $\xi_1, \xi_2$ two points on the ideal boundary of $(\mathbb{C}, g_4)$ and $c_1, c_2$ being two rays representing the $\xi_i$'s and starting from the same point

$$l(\xi_1, \xi_2) = \lim_{t \to +\infty} \frac{d(c_1(t), c_2(t))}{t} \ .$$

This definition does not depend on the chosen starting point.

*5.4.5. The Tits distance*

As for a Hadamard surface, the Tits distance Td on the ideal boundary of $(\Sigma, g_4)$ is defined as being the inner distance associated to the $l$-distance. It can be computed exactly as the Tits distance of a Hadamard surface (see [1, Section 4]). In particular, we can compute the perimeter of $(\Sigma, \mathrm{Td})$ as being the limit of

$$\frac{P_o(r)}{r} \ ,$$

where $P_o(r)$ denotes the perimeter of the circle of radius $r$ centered at a base point $o$. See more details in Section 5.3 about the computation of this perimeter.

*5.4.6. Computations in a particular case*

Suppose that the quartic differential $\varphi$ equals $z^N \mathrm{d}z^4$. We will compute the Tits distance on the ideal boundary of $(\mathbb{C}, g_0)$, where $g_0 = |z|^{N/2}|\mathrm{d}z|^2$. For $\theta$ in $[0, 2\pi]$ fixed, the curves

$$\begin{aligned} c_\theta : [0, +\infty[ &\to \mathbb{C} \\ t &\mapsto \exp(i\theta)\left(\frac{N+4}{4}t\right)^{\frac{4}{N+4}} \end{aligned}$$

are geodesic rays parametrized by arc length.

The $g_0$-geodesic circle of radius $r$ is then parametrized by:

$$\begin{aligned} \gamma_r : [0, 2\pi] &\to \mathbb{C} \\ \theta &\mapsto \exp(i\theta)\left(\frac{N+4}{4}r\right)^{\frac{4}{N+4}} \quad , \end{aligned}$$

and has perimeter

$$\int_0^{2\pi} g_{0\,\gamma_r(\theta)}(\dot{\gamma_r}, \dot{\gamma_r})^{1/2}\mathrm{d}\theta = \frac{\pi}{2}(N+4)r \ .$$

We obtain easily that the Tits boundary of $(\mathbb{C}, g_0)$ is a circle of perimeter $\frac{\pi}{2}(N+4)$.

*5.4.7. Computations in the general case*

Suppose now that the quartic differential $\varphi_4$ is a general polynomial of degree $N$, equal to

$$\varphi_4(z)\mathrm{d}z^4 = (z^N + a_{N-1}z^{N-1} + \ldots + a_1 z + a_0)\mathrm{d}z^4$$

in the $z$-coordinate. Denote by $g_4$ the induced metric with conical singularities.

**Proposition 5.4.3.** *For every positive number $\varepsilon$ there is a nonnegative number $b$ such that the identity is a $(1+\varepsilon, b)$-quasi-isometry between $g_4$ and $|z|^{N/2}|\mathrm{d}z|^2$.*

**Proof.** Let $\varepsilon$ be a positive number, and take a positive number $R$ such that

$$\frac{|a_{N-1}|}{R} + \ldots + \frac{|a_0|}{R^N} \leq \varepsilon \ .$$

Denote by $K$ the set of points $z$ such that $|z| \leq R$. By triangular inequality, we have for every $z$ in $\mathbb{C} \setminus K$

$$|z|^N (1-\varepsilon) \leq |z^N + a_{N-1}z^{N-1} + \ldots + a_0| \leq |z|^N (1+\varepsilon) \ .$$

Hence $g_4$ and $|z|^{N/2}|\mathrm{d}z|^2$ are $(1+\varepsilon)$-biLipschitz in $\mathbb{C} \setminus K$. Denote by $d_4$ the distance function associated with $g_4$ and by $d_0$ the distance function associated with $|z|^{N/2}|\mathrm{d}z|^2$. Define

$$b := \sup_{(x,y)\in K^2} \{|d_4(x,y) - d_0(x,y)|\} \ .$$

Then $d_4$ and $d_0$ are $(1+\varepsilon, b)$-quasi-isometric on $\mathbb{C}$. □

This implies that the Tits boundaries of $(\mathbb{C}, g_4)$ and of $(\mathbb{C}, |z|^{N+2}|\mathrm{d}z|^2)$ are isometric by Corollary 5.3.3. We deduce then the following.

**Proposition 5.4.4.** *Let $\mathbb{C}$ be endowed with a polynomial quartic differential $\varphi_4$ of degree $N$, and denote by $g_4$ the associated metric with conical singularities. Then the Tits boundary of $(\mathbb{C}, g_4)$ is a circle of perimeter $\frac{\pi}{2}(N+4)$.*

### 5.5. Comparison between the space boundary, the ideal boundary and the quartic boundary

Let $\Sigma$ be an asymptotically flat maximal surface of $\mathbf{H}^{2,n}$. By Lemma 5.2.3 and Theorem 5.2.4, denoting by $g$ the induced metric on $\Sigma$, the Riemann surface $(\Sigma, [g])$ is biholomorphic to the complex plane and has polynomial quartic differential $\varphi_4$. As a consequence, by Proposition 5.4.4, the metric $g_4$ associated with $\varphi_4$ has an ideal boundary of perimeter $(N+4)\pi/2$ for the Tits distance.

In this section, using a careful comparison between the induced metric $g$ and the quartic metric $g_4$, we prove the following.

**Theorem 5.5.1.** *An asymptotically flat maximal surface of $\mathbf{H}^{2,n}$ with polynomial quartic differential of degree $N$ is a polygonal surface with $N+4$ vertices.*

To this aim, we compare the Tits distances on the ideal boundaries of $(\Sigma, g)$ and of $(\Sigma, g_4)$.

#### 5.5.1. Comparison between $g$ and $g_4$

Denote by $d$ and $d_4$ the distance functions associated with $g$ and $g_4$, respectively.

**Lemma 5.5.2.** *Let $\Sigma$ be an asymptotically flat maximal surface of $\mathbf{H}^{2,n}$ with induced metric $g$ and quartic metric $g_4$, and corresponding distances $d$ and $d_4$. For every positive* *number $\varepsilon$, there exists a nonnegative number $b$ such that the identity map is a $(1+\varepsilon, b)$ quasi isometry between $(\Sigma, d)$ and $(\Sigma, d_4)$.*

**Proof.** Fix a positive number $\varepsilon$. We know, by Proposition 5.2.1 that $g$ and $g_4$ are $(1+\varepsilon)$-biLipschitz outside a compact set $K$ in $\Sigma$.

We define

$$b := \sup_{(x,y)\in K^2} \{|d(x,y) - d_4(x,y)|\} \ .$$

Take $x, y$ in $\Sigma$, and denote by $c$ the $g$-geodesic between $x$ and $y$ of length $L = d(x,y)$. If $c([0,L]) \cap K$ is empty, then the $g_4$ length of $c$ is smaller than $(1+\varepsilon)L$, hence

$$d_4(x,y) \leq (1+\varepsilon) d(x,y) \ .$$

If $c([0,L]) \cap K$ is not empty, then denote by $t_1$ the smallest real number such that $c(t_1) \in K$ and by $t_2$ the biggest real number such that $c(t_2) \in K$. Then by triangular inequality

$$d_4(x,y) \leq d_4(x, c(t_1)) + d_4(c(t_1), c(t_2)) + d_4(c(t_2), y) \ .$$

Now, by definition of $b$, and since $c([0,t_1))$ and $c((t_2, L])$ do not intersect $K$, we have

$$d_4(x,y) \leq (1+\varepsilon)t_1 + (t_2 - t_1) + b + (1+\varepsilon)(L - t_2) \leq (1+\varepsilon)L + b \ .$$

Using the same arguments, we obtain

$$d(x,y) \leq (1+\varepsilon) d_4(x,y) + b \ ,$$

so $d$ and $d_4$ are $(1+\varepsilon, b)$-quasi isometric. □

We can now use the results of Section 5.3 to prove that $(\Sigma, g)$ and $(\Sigma, g_4)$ have isometric Tits boundaries.

**Lemma 5.5.3.** *Let $\Sigma$ be an asymptotically flat maximal surface of $\mathbf{H}^{2,n}$ with induced metric $g$ and polynomial quartic differential of degree $N$. The Tits boundary of $(\Sigma, g)$ is isometric to a circle of perimeter $\frac{\pi}{2}(N+4)$.*

**Proof.** By Proposition 5.4.4, $(\Sigma, g_4)$ has a Tits boundary isometric to a circle of perimeter $(N+4)\frac{\pi}{2}$. Applying the previous Lemma and Corollary 5.3.3 gives the result. □

**Proof of Theorem 5.5.1.** By Theorem 4.8.5 and the previous lemma, $\Sigma$ is a polygonal surface with $N+4$ vertices. □

## 6. Main theorem

The Main Theorem is a consequence of Theorem 4.8.5, Theorem 5.2.4 and Theorem 5.5.1. We recall:

**Theorem 6.0.1** *(Main Theorem). Let $\Sigma$ be a maximal surface in $\mathbf{H}^{2,n}$, and denote by $g$ its induced metric. The following are equivalent:*

*(i) The space boundary of $\Sigma$ is a lightlike polygon with finitely many vertices,*
*(ii) The total curvature of $(\Sigma, g)$ is finite,*
*(iii) The surface $(\Sigma, g)$ is asymptotically flat.*

*Moreover, if the conditions are satisfied, and denoting by $N+4$ the number of vertices of $\partial^s\Sigma$, the Riemann surface $(\Sigma, [g])$ is of parabolic type, has polynomial quartic differential of degree $N$, the total curvature of $\Sigma$ equals $-\frac{\pi}{2}N$ and the perimeter of circles of radius $r$ is asymptotic to $(2\pi + N\frac{\pi}{2})r$.*

**Proof.** Let $\Sigma$ be a maximal surface in $\mathbf{H}^{2,n}$. By Theorem 4.8.5, the two first items are equivalent, and if $\Sigma$ is polygonal with $N+4$ vertices then its total curvature equals $-\frac{\pi}{2}N$. By Corollary 4.3.10 a polygonal surface is asymptotically flat. By Theorem 5.2.4, an asymptotically flat maximal surface has parabolic type and polynomial quartic differential. At the end, if $\Sigma$ is asymptotically flat, of parabolic type and has polynomial quartic differential of degree $N$, then by Theorem 5.5.1 it is polygonal with $N+4$ vertices. In addition, a polygonal surface has Tits boundary of perimeter $2\pi + N\frac{\pi}{2}$, and Theorem 5.3.1 proves the last assertion. □